\documentclass[final]{siamltexmm}
\usepackage{verbatim,amsmath,amsfonts,epsfig,graphics,setspace,amssymb,
float,stmaryrd}
\usepackage{epsfig}
\usepackage{algorithm}
\usepackage{algorithmicx}
\usepackage{algpseudocode}
\usepackage{enumerate,color}
\newcounter{rmrk}[section]

\numberwithin{equation}{section}
\newtheorem{cor}{Corollary}[section]
\newtheorem{thm}{Theorem}[section]
\newtheorem{lem}{Lemma}[section]

\newtheorem{defn}{Definition}[section]
\newtheorem{assume}{Assumption}[section]

\numberwithin{equation}{section}
\newcommand{\norm}[1]{\left\Vert#1\right\Vert}
\newcommand{\abs}[1]{\left\vert#1\right\vert}
\newcommand{\set}[1]{\left\{#1\right\}}

\newcommand{\mbf}[1]{\mathbf{#1}}

\title{Solving Stochastic Inverse Problems using Sigma-Algebras on Contour Maps}
\date{\today}
\author{T.~Butler\thanks{Department of Mathematical and Statistical Sciences, University of Colorado Denver, Denver, CO 80202 ({\tt
Troy.Butler@ucdenver.edu}).}
\and D.~Estep\thanks{Department of Statistics, Colorado State University, Fort Collins,
CO 80523 ({\tt estep@stat.colostate.edu}).}
\and S.~Tavener\thanks{Department of Mathematics, Colorado State University, Fort Collins,
CO 80523 ({\tt tavener@math.colostate.edu})}
\and T.~Wildey \thanks{Sandia National Labs, Albuquerque, NM 87185 ({\tt tmwilde@sandia.gov}).
Sandia is a multiprogram laboratory operated by Sandia Corporation, a wholly owned subsidiary of Lockheed Martin company, for the United States Department of Energy's National Nuclear Security Administration under Contract DE-AC04-94-AL85000.}
\and C.~Dawson\thanks{Institute for Computational Engineering and Sciences (ICES), University of Texas at Austin, Austin, Texas 78712 ({\tt
clint@ices.utexas.edu}).}
\and L.~Graham \thanks{Institute for Computational Engineering and Sciences (ICES), University of Texas at Austin, Austin, Texas 78712 ({\tt lgraham@ices.utexas.edu}).}
}

\begin{document}
\maketitle
\newcommand{\slugmaster}{%
\slugger{juq}{xxxx}{xx}{x}{x--x}}

\begin{abstract}
We compute approximate solutions to inverse problems for determining parameters in differential equation models with stochastic data on output quantities. The formulation of the problem and modeling framework define a solution as a probability measure on the parameter domain for a given $\sigma-$algebra. In the case where the number of output quantities is less than the number of parameters, the inverse of the map from parameters to data defines a type of generalized contour map.  The approximate contour maps define a geometric structure on events in the $\sigma-$algebra for the parameter domain.
We develop and analyze an inherently non-intrusive method of sampling the parameter domain and events in the given $\sigma-$algebra to approximate the probability measure. We use results from stochastic geometry for point processes to prove convergence of a random sample based approximation method. We define a numerical $\sigma-$algebra on which we compute probabilities and derive computable estimates for the error in the probability measure. We present numerical results to illustrate the various sources of error for a model of fluid flow past a cylinder.
\end{abstract}


\pagestyle{myheadings} \thispagestyle{plain} \markboth{T.~Butler et.~al.}{Sigma-Algebras on Contour Maps}

%%%%%%%%%%%%%%%%%%%%%%%%%%%%%%%%%%%%%%%%%%%%%%%%%%%%%%%%%%%%%%%%%%%%
%%%%%%%%%%%%%%%%%%%%%%%%%%%%%%%%%%%%%%%%%%%%%%%%%%%%%%%%%%%%%%%%%%%%
\section{Introduction}\label{Sec:Intro}
%%%%%%%%%%%%%%%%%%%%%%%%%%%%%%%%%%%%%%%%%%%%%%%%%%%%%%%%%%%%%%%%%%%%
%%%%%%%%%%%%%%%%%%%%%%%%%%%%%%%%%%%%%%%%%%%%%%%%%%%%%%%%%%%%%%%%%%%%

There are a seemingly unlimited number of proposed methods for solving forward and inverse problems of interest in the scientific community. The computational model considered and the framework assumed in formulating the forward and inverse problems play a critical role in defining the solution methodology. In this work, we consider a set-approximation method for approximating solutions to a stochastic inverse problem for deterministic models formulated within a measure-theoretic framework \cite{Butler2014a}. The solutions we seek are non-parametrically defined probability measures on the model parameter domain defined by uncertain model inputs, e.g.~viscosity of a fluid, initial or boundary conditions, and/or domain parameters such as locations of holes.

The models we consider are deterministic physics-based models for which a limited number of physical observations of the solution are available. We call the observations quantities of interest (QoI) and are particularly interested in the situation where there are a smaller number of QoI than model parameters. In this case, solution of the deterministic inverse problem is complicated by the fact that the map defined by inverting the QoI is a set-valued map. This is often referred to as ill-posedness. Common methods such as regularization effectively alter the QoI map so that it has a well-posed deterministic inverse albeit with an altered geometric structure.  Introducing stochasticity in the form of probability measures on inputs and/or outputs also fails to address this fundamental geometric issue\footnote{Statistical inverse problems where the QoI map is replaced by a statistical map constitute an entirely different framework for formulating inverse problems, see \cite{Butler2014a, Breidt2011} for more detailed relations to statistical and Bayesian inference problems.}. The method we use is based on a computational measure-theoretic approach that fully exploits the set-valued inverses directly preserving specific geometric information contained in the map from inputs to outputs \cite{Breidt2011, Butler2012a, Butler2013a, Butler2014a}. Approximate solutions are defined using simple function approximations to the density of the unique probability measure on the parameter domain and given $\sigma-$algebra.

There are two common steps in computing any simple function approximation to a measurable function such as the density (i.e.~Radon-Nikodym derivative) of a probability measure:
\begin{enumerate}[(1)]
	\item identify measurable sets in a given $\sigma-$algebra partitioning the domain of the function;
	\item assign a nominal function value on each of these sets.
\end{enumerate}
These steps are explicitly used in classical examples applying the Lebesgue Monotone Convergence Theorem. First, a specific sequence of partitions on the range of a given measurable real-valued non-negative function is identified. The preimages of the sets in each partition are used to identify the measurable sets in the $\sigma$-algebra on the domain of the function. This completes the first step and requires evaluation of the inverse map to identify inverse sets. The second step is trivial in this case where the infimum of each output set is assigned to its preimage to define the simple function approximations for each partition.

A probability measure with density on the parameter domain and original $\sigma-$algebra is a solution to the stochastic inverse problem if the QoI maps this parameter density to the output density. Unique solutions to the stochastic inverse problem exist in a \emph{contour} $\sigma-$algebra \cite{Butler2014a} embedded in the original $\sigma-$algebra. Moreover, combining the Disintegration Theorem with an Ansatz describing probabilities on the contour events proves unique solutions to the stochastic inverse problem exist in the original parameter $\sigma-$algebra \cite{Butler2014a}. As shown in \cite{Breidt2011, Butler2014a}, the key step in constructing any simple function approximation to the density solving the stochastic inverse problem is the approximation of the QoI contour map in the given $\sigma-$algebra. The contour map contains all the geometric information available by the QoI map.
It is then straightforward to assign probability values to any approximate partitioning sets of the contour events in the original $\sigma-$algebra \cite{Breidt2011, Butler2012a, Butler2013a, Butler2014a} such that a sequence converges using either the Monotone Convergence Theorem \cite{Breidt2011} or the Lebesgue Dominated Convergence Theorem \cite{Butler2014a}.

The formulations of the forward and inverse problems within this measure-theoretic framework, the existence and uniqueness of solutions, and the approximations of solutions by deterministic approximation techniques have been recently studied \cite{Breidt2011, Butler2012a, Butler2013a, Butler2014a}. The focus of this work is the development and analysis of sampling techniques within a parameter domain to implicitly define approximating sets of both a contour $\sigma-$algebra and the original $\sigma-$algebra. In this way we create a new computational algorithm for the probability measure based on the samples and their approximation properties. We prove that a counting measure approximation converges to the exact probability measure solving the stochastic inverse problem as the number of samples increases. Finally, for the computed probabilities of events, we derive computable error estimates and bounds for the effects of using a finite number of samples and numerical errors in the evaluation of the model.

The outline of this paper is as follows. In Section~\ref{S:Formulation}, we summarize the measure-theoretic framework and formulation of the stochastic inverse problem. In Section~\ref{S:Algorithms}, we summarize the theoretical results involving $\sigma-$algebras on contour maps and solutions to the stochastic inverse problem. We then present computational algorithms including a point-sample based algorithm. We describe the implicit set-approximation properties defined by the sampling of points in parameter space and define a counting measure approximation to the inverse probability measure. In Section~\ref{S:Theory}, we prove that the counting measure approximation converges to the correct probability measure for arbitrary events in the original $\sigma-$algebra. In Section~\ref{S:Errors}, we summarize the types and sources of error in the computed probabilities and how we estimate and/or bound these errors. In Section~\ref{S:Numerics}, we provide numerical results for a Navier-Stokes flow past a cylinder with uncertain viscosity and cylinder location.

\section{Formulating stochastic inverse problems in a measure-theoretic framework}\label{S:Formulation}
%{\bf Don and I will work on cutting this section down and probably moving parts to the introduction. We need a very concise and clear explanation of this so we don't have to repeat it for every paper.}

Let ${\Lambda}\subset\mathbb{R}^n$ denote the parameter domain for the model. We assume $\Lambda$ is a metric space whose metric is specified as part of the model. Thus, there is an induced topology on $\Lambda$ and Borel $\sigma-$algebra $\mathcal{B}_{\Lambda}$. We define the measure space $({\Lambda}, \mathcal{B}_{{\Lambda}}, \mu_{{\Lambda}})$ using the measure $\mu_{{\Lambda}}$ induced by the metric \cite{Friedman_Book}. Let $Q:\Lambda\to\mathcal{D}$ denote the vector-valued QoI map where $\mathcal{D}\subset\mathbb{R}^m$, $m \leq n$, is defined by functionals of solution to the model.
\begin{defn}\label{defn:nonlin_indepdence}
We say that the component maps of $m$-dimensional locally differentiable vector-valued map $Q(\lambda)$ are {\bf geometrically distinct} (GD) if the Jacobian of $Q$ has full rank at every point in ${\Lambda}$.
\end{defn}

In practice, since we deal with measurable events of non-zero measure, we can weaken local differentiabillity of $Q$ and/or the full rank of the Jacobian of $Q$ to hold at a.e.~point in $\Lambda$\footnote{Since we use probability densities, which are $L^1$ functions, to compute probabilities of events, any set of zero measure has no affect on the final solutions and can simply be ignored or removed from $\Lambda$ in the formulation of the inverse problem for convenience.}. Assuming $Q$ is locally differentiable with GD components, $\mu_{{\Lambda}}$ defines a ``push-forward'' measure  $\mu_{\mathcal{D}}$ on $(\mathcal{D},\mathcal{B}_{\mathcal{D}})$, where $\mathcal{B}_{\mathcal{D}}$ is the Borel $\sigma-$algebra on $\mathcal{D}$, and for any $A\in \mathcal{B}_{\mathcal{D}}$, $\mu_{\mathcal{D}}(A) = \mu_{{\Lambda}}(Q^{-1}(A))$.
%\begin{equation*}
%\mu_{\mathcal{D}}(A) = \int_A d\mu_{\mathcal{D}} = \int_{Q^{-1}(A)}  d\mu_{{\Lambda}} = \mu_{{\Lambda}}(Q^{-1}(A)),
%\end{equation*}
This yields the measure space $(\mathcal{D}, \mathcal{B}_{\mathcal{D}}, \mu_{\mathcal{D}})$.  The measures $\mu_{\Lambda}$ and $\mu_{\mathcal{D}}$ are volume measures not probability measures. We often assume the probability measures are absolutely continuous with respect to the volume measures in which cases we typically refer to the associated probability density functions (i.e.~the Radon-Nikodym derivatives).%These measures determine the volumes, not probabilities, of measurable sets (i.e.~events). % and are used within the theoretical and algorithmic framework to define a unique solution to the stochastic inverse problem.

%\paragraph*{A forward stochastic sensitivity analysis problem} We assume the probability space $(\Lambda, \mathcal{B}_{{\Lambda}}, P_{{\Lambda}})$ is given where probability measure $P_{\Lambda}$ is absolutely continuous with respect to $\mu_{\Lambda}$, i.e.~$\lambda$ is a continuous random variable on $\Lambda$ with density $\rho_{\Lambda}$ defined by the Radon-Nikodym derivative $dP_{\Lambda}/d\mu_{\Lambda}$.
%The goal of the stochastic sensitivity analysis  is to compute the  probability density $\rho_{\mathcal{D}}$ with respect to $\mu_{\mathcal{D}}$ induced by $\rho_{{\Lambda}}$, i.e. for event $A\in \mathcal{B}_{\mathcal{D}}$,
%\begin{equation*}
%P_{\mathcal{D}}(A) = \int_A \rho_{\mathcal{D}} d\mu_{\mathcal{D}} = \int_{Q^{-1}(A)} \rho_{{\Lambda}} d\mu_{{\Lambda}} = P_{{\Lambda}} (Q^{-1}(A)).
%\end{equation*}
%This is a familiar problem in uncertainty quantification, and a simple solution method is the classic Monte Carlo method, but more sophisticated methods using the sensitivity of $Q$ exist, e.g.~see \cite{EstepNeckelsPapers, MalquistEstepPapers}.

%\paragraph*{An inverse stochastic sensitivity analysis problem}

%Description of the inverse problem is complicated by the fact that the deterministic inverse problem is most often defined by a set-valued function.

The inverse sensitivity problem of interest is the direct inversion of the forward stochastic sensitivity analysis problem\footnote{The forward problem is standard: given probability measure $P_{\Lambda}$  on measurable space $(\Lambda,\mathcal{B}_{\Lambda})$ compute the  probability measure $P_{\mathcal{D}}$ on measurable space $(\mathcal{D},\mathcal{B}_{\mathcal{D}})$.}.
In other words, the stochastic inverse problem is to determine probability measures on measurable space $(\Lambda,\mathcal{B}_{\Lambda})$ such that the push-forward probability measures match given probability measures on $(\mathcal{D},\mathcal{B}_{\mathcal{D}})$. Fundamentally, this requires a description of the mapping between sets in the various $\sigma-$algebras $\mathcal{B}_{\mathcal{D}}$ and $\mathcal{B}_{\Lambda}$.

%Given probability space $(\mathcal{D},\mathcal{B}_{\mathcal{D}},P_{\mathcal{D}})$ where probability measure $P_{\mathcal{D}}$ has density (i.e.~Radon-Nikodym derivative) $\rho_{\mathcal{D}}$ with respect to $\mu_{\mathcal{D}}$, the goal is to compute probability measure $P_{\Lambda}$ on measurable space $(\Lambda,\mathcal{B}_{\Lambda})$ given probability measure $P_{\mathcal{D}}$ on measurable space $(\mathcal{D},\mathcal{B}_{\mathcal{D}})$. We assume that the probability measures are all absolutely continuous with respect to the underlying ``volume'' measures and often write the corresponding densities as Radon-Nikodym derivatives.

When $m<n$,  %by the assumption that $Q$ has GD component maps,
the Implicit Function Theorem guarantees that for any fixed point in $\mathcal{D}$, there exists some union of locally continuous $(n-m)$-dimensional manifolds defining the set-valued inverse. We use the notation and definitions of \cite{Butler2014a, Breidt2011, Butler2012a, Butler2013a} and define a {\em generalized contour} as any such set-valued inverse of the map $Q$.
We assume that $\mathcal{D}$ is defined by $Q(\Lambda)$, i.e. $\mathcal{D}$ is the range of the QoI map containing all possible physically observable data that can be mapped to from $\Lambda$. Thus, we can decompose ${\Lambda}$ into a union of generalized contours in 1-1 correspondence with the points in $\mathcal{D}$. In other words, the map $Q$ defines a type of generalized contour map on $\Lambda$.
There exists a (possibly piecewise-defined) continuous $m$-dimensional indexing manifold, called a {\em transverse parameterization}, in $\Lambda$ defining a bijection between $\mathcal{D}$ and the generalized contours, see \cite{Butler2014a} for more details. As way of analogy, a specific transverse parameterization is like a particular path of ascent up a mountain that a hiker plots out using a contour map of elevations.

We note that the contour map obtained from $Q^{-1}$ defines an equivalence class relation on ${\Lambda}$. Denote this space of equivalence classes as $\mathcal{L}$ where each point in $\mathcal{L}$ corresponds to a set of points in ${\Lambda}$. We obtain a measure space $(\mathcal{L}, \mathcal{B}_{\mathcal{L}}, \mu_{\mathcal{L}})$, where the $\sigma-$algebra  $\mathcal{B}_{\mathcal{L}}$ on $\mathcal{L}$ can be generated using inverse images of a collection of Borel sets in $\mathcal{B}_{\mathcal{D}}$ and  the volume measure $\mu_{\mathcal{L}}$ induced by $\mu_{\mathcal{D}}$.

We exploit the equivalence relation to define the induced $\sigma-$algebra $\mathcal{C}_{{\Lambda}}$ on $\Lambda$ that can be generated from the set of equivalence classes for a set of generating events in $\mathcal{B}_{\mathcal{L}}$. For $m<n$, $\mathcal{C}_{{\Lambda}}$ is a proper subset of $\mathcal{B}_{{\Lambda}}$. We define $\mathcal{C}_{\Lambda}$ as the {\em contour $\sigma-$algebra on $\Lambda$} and call events in $\mathcal{C}_{{\Lambda}}$ {\em contour events}. The geometric structure of $Q$ on $\Lambda$ is fully exploited to define the contour events in the contour $\sigma$-algebra $\mathcal{C}_{\Lambda}$. We emphasize that events in $\mathcal{C}_{\Lambda}$ can be uniquely determined by events in $\mathcal{B}_{\mathcal{D}}$.

The goal is to define solutions to the stochastic inverse problem on $(\Lambda,\mathcal{B}_{\Lambda})$, e.g.~in terms of a density defined on {\em points in $\Lambda$ not points on a contour map}. Below, we prove there exists a $\sigma-$algebra defined on the generalized contour map on $\Lambda$ equivalent to the $\sigma-$algebra $\mathcal{B}_{\Lambda}$ on $\Lambda$. We then use this $\sigma-$algebra equivalency in the Disintegration Theorems that follow to prove the existence and uniqueness of solutions to the stochastic inverse problem.

 % Below, we review the theoretical and algorithmic framework for the existence, uniqueness, and approximation of solutions to this inverse sensitivity problem.

\section{Solving the stochastic inverse problem using sigma-algebras on contour maps}\label{S:Algorithms}

A common starting point for constructing a measure on a $\sigma-$algebra is to first define an algebra used to generate the $\sigma-$algebra of interest. The next steps can vary, but typically we define a premeasure on the algebra which induces an outer-measure in a natural way. Finally, employ Carath\'{e}odory's Theorem to extend the outer-measure uniquely to a complete measure on the generated $\sigma-$algebra. We can define an algebra with elementary families of elementary sets for which a clear notion of measure is defined, e.g.~by using a partitioning of a space where sets are either disjoint or possibly intersect only at boundaries such as hypercubes or generalized rectangles in $\mathbb{R}^n$. We will exploit such generating sets in much of the theory below.

In Section~\ref{S:SetDisintegration}, we relate $\sigma-$algebras on $\Lambda$ to $\sigma-$algebras on the generalized contour maps. In Section~\ref{S:TheoryReview}, we summarize the theory of the existence and uniqueness of solutions to the stochastic inverse problem with respect to these $\sigma-$algebras. In Section~\ref{S:ComputingDetails}, we describe the approximation of solutions.

\subsection{Sigma-Algebras on COntour Maps (SACOM)}\label{S:SetDisintegration}

The generalized contours define a type of contour map on $\Lambda$. Given the geometry of these generalized contours and a topology, we can define $\sigma-$algebras on contour maps (SACOM). In Section~\ref{S:TheoryReview} below, we show that the solution to the stochastic inverse problem can be solved uniquely using SACOM. The goal is to relate the $\sigma-$algebras on $\Lambda$ (e.g.~$\mathcal{C}_{\Lambda}$ or $\mathcal{B}_{\Lambda}$) to SACOM.
%The advantage of relating $\sigma-$algebras on $\Lambda$ to SACOM is made clear by the Disintegration Theorem in Section~\ref{S:TheoryReview} which defines solutions of the stochastic inverse problem with respect to SACOM.

To motivate what follows, consider the case where the map $Q$ is linear. In this case, the $(n-m)$-dimensional generalized contours are hyperplanes. For simplicity, let $\Lambda = \mathbb{R}^n$, and assume we use the typical Euclidean distance metric and induced topology giving the Borel $\sigma-$algebra on the product space $\Lambda$. Let $\pi_{\ell}(\lambda)=\ell$ denote the projection map from $\Lambda$ to the equivalence class $\ell\in\mathcal{L}$ containing $\lambda$. For arbitrary $\ell\in\mathcal{L}$, let $C_{\ell} := \pi^{-1}(\ell)$ denote the associated generalized contour as a subset of $\Lambda$.
Changing coordinates with respect to the directions orthogonal and parallel to the hyperplanes  gives $\mathcal{B}_{\Lambda}=\mathcal{B}_{\mathcal{L}}\otimes\mathcal{B}_{C_{\ell}}$\footnote{See Propositions 1.5 in \cite{Folland_Book} identifying $\Lambda = \mathbb{R}^n$ and from the change of coordinates $\mathcal{L}=\mathbb{R}^m$ and $C_{\ell} = \mathbb{R}^{(n-m)}$.}. Identifying such a product decomposition is useful for applying Fubini's theorem in order to integrate certain functions using iterated integrals. The key is that the domain of integration, which is a set in the original $\sigma-$algebra, can be ``cut'' into products of lower-dimensional sets in the component measure space $\sigma-$algebras. In this case, such sets are ``cut'' into products of sets both along and transverse to the generalized contours.  Below, we extend this for the general case of nonlinear generalized contours indexed by a nonlinear manifold defining a transverse parameterization.

Let $\mathcal{F}_{C_{\ell}}$ denote a $\sigma-$algebra on a given $C_{\ell}$. Below, we describe two natural choices for $\mathcal{F}_{C_{\ell}}$ for each $C_{\ell}$ defining a family of measurable spaces $\set{C_{\ell}, \mathcal{F}_{C_{\ell}}}$ indexed by $\ell\in\mathcal{L}$.

Any space contains the so-called trivial $\sigma-$algebra consisting of the empty set or the entire space. Therefore, there exists at least one $\sigma-$algebra on each $C_{\ell}$, given by $\mathcal{T}_{C_{\ell}}=\set{C_{\ell},\emptyset}$. Alternatively, using the induced topology on each $C_{\ell}$, we define Borel $\sigma-$algebras $\mathcal{B}_{C_{\ell}}$. These are equivalently defined by $\set{A\cap C_{\ell} \, | \, A\in\mathcal{B}_{\Lambda}}$\footnote{See Lemma 6.2.4 of \cite{Bogachev_vol2}.}. In other words, the Borel $\sigma-$algebra on $C_{\ell}$ can be defined by the restriction of Borel measurable sets in $\mathcal{B}_{\Lambda}$ to the given generalized contour.   These two particular choices for $\mathcal{F}_{C_{\ell}}$ prove useful in the Disintegration Theorems of Section~\ref{S:TheoryReview} for showing existence and uniqueness of solutions to the stochastic inverse problem.

As in many of the classical proofs of theorems involving product $\sigma-$algebras, we make explicit use of the generating sets for each $\sigma-$algebra. Since each measurable space we consider is metrizable, we assume the generating sets are taken from the implied topology. For each $\ell\in\mathcal{L}$, let $F_{C_{\ell}}$ denote any family of subsets of $C_{\ell}$ generating $\mathcal{F}_{C_{\ell}}$. A standard measure-theory result states that $\mathcal{F}_{C_{\ell}}$ is the unique smallest $\sigma-$algebra generated by $F_{C_{\ell}}$ for each $\ell\in\mathcal{L}$. Similarly, we let $F_{\mathcal{L}}$ denote a generating set for the Borel $\sigma-$algebra $\mathcal{B}_{\mathcal{L}}$.

\begin{defn}\label{d:tp_algebra}
We define the \textbf{transverse product $\sigma-$algebra} as the smallest $\sigma-$algebra on $\Lambda$ generated by $\bigcup_{\ell\in A} B_{C_{\ell}}$ %for all countable sets\footnote{Note that we have assumed $\Lambda$ is a compact metric space so is separable.} of
where $A\in F_{\mathcal{L}}$ and $B_{C_{\ell}}\in F_{C_{\ell}}$ for all $\ell$. We denote this $\sigma-$algebra by $ \bigotimes_{\mathcal{B}_{\mathcal{L}}}\set{\mathcal{F}_{C_{\ell}}}$.
\end{defn}

\begin{thm}\label{thm:ContourSigmaAlgebraDecomp}
$\mathcal{C}_{\Lambda}=\bigotimes_{\mathcal{B}_{\mathcal{L}}} \set{\mathcal{T}_{C_{\ell}}}$.
\end{thm}

{\em Proof:}
This follows immediately from the definition of contour $\sigma-$algebra $\mathcal{C}_{\Lambda}$ using the equivalence relation determined by $Q^{-1}$ and the generating sets for $\mathcal{B}_{\mathcal{L}}$ and $\set{\mathcal{T}_{C_{\ell}}}$. $\Box$

%From the ass $\mathcal{E}_A\in\mathcal{B}_{\mathcal{L}}$. Similarly,
%For any $A\in\mathcal{B}_{\Lambda}$ and $\ell\in\mathcal{L}$, let $A_{\ell}:= \pi_{\ell}^{-1}(\ell)\cap A$. Note that $A_{\ell}$ is a subset of $\Lambda$ with an induced topology. Let $\mathcal{B}_{\ell}$ denote the Borel $\sigma-$algebra given by $\set{A_{\ell} \, | \, A\in\mathcal{B}_{\Lambda}}$.

\begin{thm}\label{thm:OriginalSigmaAlgebraDecomp}
$\mathcal{B}_{\Lambda}=\bigotimes_{\mathcal{B}_{\mathcal{L}}} \set{\mathcal{B}_{C_{\ell}}}$.
%If $\mathcal{L}$ is separable, then $\mathcal{B}_{\Lambda}=\bigotimes_{\mathcal{B}_{\mathcal{L}}} \set{\mathcal{B}_{C_{\ell}}}$
%\footnote{Note that the inclusion $\mathcal{B}_{\mathcal{L}} \vec{\otimes}\set{\mathcal{B}_{C_{\ell}}}\subset\mathcal{B}_{\Lambda}$ is not at all obvious (and is possibly untrue) since we do not define $\mathcal{B}_{\mathcal{L}} \vec{\otimes}\set{\mathcal{B}_{C_{\ell}}}\subset\mathcal{B}_{\Lambda}$ as generated by countable unions of $\ell\in\mathcal{L}$ yet $\mathcal{B}_{\Lambda}$ is a countably generated $\sigma-$algebra. Fortunately, it is also not necessary.}.
\end{thm}

Before we prove this theorem, we note that the inclusion $\bigotimes_{\mathcal{B}_{\mathcal{L}}} \set{\mathcal{B}_{C_{\ell}}}\subset\mathcal{B}_{\Lambda}$ is not at all obvious since $A\in F_{\mathcal{L}}$ may be uncountable and $\sigma-$algebras are closed under countable unions.  Since we assume $Q$ is locally differentiable, there exists a bijection between $\mathcal{L}$ and a transverse parameterization in $\Lambda$ \cite{Breidt2011, Butler2014a}. Since the transverse parameterization is a piecewise-continuous $m$-dimensional manifold in separable space $\Lambda$, it follows that $\mathcal{L}$ is separable. This is used below. Also, we let $\mathcal{E}_A := \pi_{\ell}(A)$ for any $A\in\mathcal{B}_{\Lambda}$. We use this notation elsewhere as convenient.
%This inclusion is also not required in any of the results that follow.

{\em Proof:}
Let $F_{\Lambda}$ denote any Borel generating set of $\mathcal{B}_{\Lambda}$ and consider any $A\in F_{\Lambda}$. Restrict the map $Q$ to $A$, then by the assumption of GD component maps, there exists an $(n-m)$-dimensional (piecewise) continuous manifold defining a transverse parameterization on $A$ \cite{Butler2014a}. Moreover, the transverse parameterization is a Borel set in $\Lambda$ by assumption of the local differentiability of $Q$ and is in 1-1 correspondence with $\mathcal{E}_A$. Thus, we have\footnote{See Theorem 6.9.7~\cite{Bogachev_vol2}.} that the Borel set $A$ defines the unique {\em Borel set} $\mathcal{E}_A\subset\mathcal{L}$. Restricting any generating sets of $\mathcal{B}_{\Lambda}$ to a generalized contour $C_{\ell}$ defines a generating set for $\mathcal{B}_{C_{\ell}}$. From Definition~\ref{d:tp_algebra}, we have that $A\in \bigotimes_{\mathcal{B}_{\mathcal{L}}} \set{\mathcal{B}_{C_{\ell}}}$. It follows that $\mathcal{B}_{\Lambda}\subset \otimes_{\mathcal{B}_{\mathcal{L}}} \set{\mathcal{B}_{C_{\ell}}}$\footnote{See Lemma 1.1 of \cite{Folland_Book}.}.

Suppose we are given Borel generating set $F_{\mathcal{L}}$ and family of Borel generating sets $\set{F_{C_{\ell}}}$ defined by all the open sets on these spaces. For all $\ell\in\mathcal{L}$, any Borel set on $C_{\ell}$ is also Borel in $\Lambda$. By assumption, $Q^{-1}(A)$ is Borel in $\Lambda$ for any $A\in F_{\mathcal{L}}$. If $A$ is countable, then $\bigcup_{\ell\in A} B_{\ell} \in \mathcal{B}_{\Lambda}$ for any $B_{\ell}\in F_{\ell}$. Suppose $A$ is uncountable. Let $C_{\mathcal{L}}$ denote a countable dense set in $\mathcal{L}$ and $E_{\mathcal{L}}$ the collection of open balls in $\mathcal{L}$ with rational radius and center in $C_{\mathcal{L}}$. Then every open set in $\mathcal{L}$ is a countable union of open balls in $E_{\mathcal{L}}$. In other words, $A$ is countably generated by $E_{\mathcal{L}}$, so the set $\bigcup_{\ell\in A} B_{C_{\ell}}$ is countably generated by unions of Borel sets and is itself Borel. Thus, $\bigcup_{\ell\in A} B_{C_{\ell}}\in\mathcal{B}_{\Lambda}$, so $\bigotimes_{\mathcal{B}_{\mathcal{L}}} \set{\mathcal{B}_{C_{\ell}}} \subset \mathcal{B}_{\Lambda}$. $\Box$

\subsection{Existence and uniqueness of solutions}\label{S:TheoryReview}

To solve the stochastic inverse problem, we use a form of the Disintegration Theorem \cite{Butler2014a, Chang_Pollard, Dellacherie_Meyer}, which is a powerful theoretical tool for rigorous definition of conditional probabilities.  We focus on convenient forms of the theorem written for probabilities and direct the interested reader to \cite{Butler2014a} for a more thorough presentation. Using Theorem~\ref{thm:ContourSigmaAlgebraDecomp} and following the steps of \cite{Butler2014a}, we have,
\begin{thm}[Disintegration of Contour Map Probabilities]\label{thm:ProbabililtyDisintegrationTrivial}
Let $({\Lambda}, \mathcal{C}_{{\Lambda}})$ be a measurable space. Assume that $P_{{\Lambda}}$ is a probability measure on $({\Lambda}, \mathcal{C}_{{\Lambda}})$. There exists a family of conditional probability measures $\set{P_{\ell}}$ on $\set{(C_{\ell},\mathcal{T}_{\Lambda})}$
giving the disintegration,
\begin{equation}
	P_{{\Lambda}}(A)  = \int_{\pi_{\mathcal{L}}(A)} \bigg( \int_{\pi^{-1}_{\mathcal{L}}(\ell) \cap A} \, dP_{\ell}(\lambda) \bigg) dP_{\mathcal{L}}(\ell), \ \forall A\in\mathcal{C}_{\Lambda}.
\end{equation}
\end{thm}

%\begin{thm}[The Disintegration Theorem]\label{thm:disintegration}
%Let $(\Lambda, \mathcal{B}_{\Lambda})$ be a measurable space and  $Q:{\Lambda}\to\mathcal{D}$ be a measurable map with GD component maps, and assume that $\Psi$ is a measure on $(\Lambda, \mathcal{B}_{\Lambda})$. There is a family of measures $\{\Psi_{\ell}\}$ on $({\Lambda},\mathcal{B}_{\Lambda})$ defined for a.e. $\ell \in \mathcal{L}$ such that
%\begin{equation*}
%\Psi_\ell (\lambda) = 0, \quad \lambda \in {\Lambda} \setminus \pi^{-1}_{\mathcal{L}}(\ell) , \quad \text{a.e. } \ell \in \mathcal{L},
%\end{equation*}
%i.e., $\Psi_\ell (A) =  \Psi_\ell ( \pi^{-1}_{\mathcal{L}}(\ell) \cap A ) $ for all $A \in \mathcal{B}_{\Lambda}$, and which gives the  following { disintegration} for $\Psi$,
%\begin{equation}\label{eq:disintegration1}
%	\Psi(A)  = \int_{\mathcal{E}_A} \Psi_{\ell}(A) \, d\mu_{\mathcal{L}}(\ell) = \int_{\mathcal{E}_A} \bigg( \int_{\pi^{-1}_{\mathcal{L}}(\ell) \cap A} \, d\Psi _{\ell} (\lambda) \bigg) \, d\mu_{\mathcal{L}}(\ell),
%\end{equation}
%for $A \in \mathcal{B}_{\Lambda}$.
%\end{thm}

It is clear we can compute the induced probability measure $P_{\Lambda}$ on $(\Lambda,\mathcal{C}_{\Lambda})$ defined by $P_{\Lambda}(A) = P_{\mathcal{L}}(\pi_{\mathcal{L}}(A))= P_{\mathcal{D}}(Q(A))$, and $P_{\mathcal{L}}$ is defined by a probability density $\rho_{\mathcal{L}}$ on $(\mathcal{L},\mathcal{B}_{\mathcal{L}})$ with respect to $\mu_{\mathcal{L}}$,
\begin{equation*}
P_{\mathcal{L}}(A) = \int_{\pi_{\mathcal{L}}(A)} \rho_{\mathcal{L}} d\mu_{\mathcal{L}} = \int_{Q(A)}  \rho_{\mathcal{D}} d\mu_{\mathcal{D}} = P_{\mathcal{D}}(Q(A)), \ \forall A\in\mathcal{B}_{\mathcal{L}}.
\end{equation*}
From Theorem~\ref{thm:ProbabililtyDisintegrationTrivial}, the conditional probability measures for $P_{\Lambda,\mathcal{C}_{\Lambda}}$ are given by $P_{\ell}(A) = 1$ if $\pi^{-1}(\ell)\subset A\in\mathcal{C}_{\Lambda}$ and $P_{\ell}(A)=0$ otherwise. This proves the following
\begin{thm}
The stochastic inverse problem has a unique solution on $(\Lambda,\mathcal{C}_{\Lambda})$.
\end{thm}

The primary goal is not a probability measure on $(\Lambda,\mathcal{C}_{\Lambda})$, but a probability measure $P_{\Lambda}$ on $(\Lambda,\mathcal{B}_{\Lambda})$. Since we will often work with densities given as Radon-Nikodym derivatives (i.e.~$dP_{\Lambda}/d\mu_{\Lambda}$), we find useful the following corollary of the Disintegration Theorem \cite{Butler2014a},
\begin{cor}[Disintegration of the Volume Measure]\label{cor:disintegrate_volume}
There exists a family of volume measures $\set{\mu_{C_{\ell}}}$ on $\set{(C_{\ell},\mathcal{B}_{C_{\ell}})}$ such that for any $A\in\mathcal{B}_{\mathbf{\Lambda}}$,
\[
	\mu_{\mathbf{\Lambda}}(A) = \int_{\mathcal{E}_A} \mu_{C_{\ell}}(\pi_{\mathcal{L}}^{-1}(\ell) \cap A) \, d\mu_{\mathcal{L}}(\ell) = \int_{\mathcal{E}_A} \int_{\pi^{-1}_{\mathcal{L}}(\ell) \cap A} \, d\mu_{C_{\ell}}(\lambda)\, d\mu_{\mathcal{L}}(\ell).
\]
\end{cor}

Using Theorem~\ref{thm:OriginalSigmaAlgebraDecomp}, we arrive at the more common form of the Disintegration Theorem for probability measures given by
\begin{thm}\label{thm:ProbabililtyDisintegration}
Let $({\Lambda}, \mathcal{B}_{{\Lambda}})$ be a measurable space. Assume that $P_{{\Lambda}}$ is a probability measure on $({\Lambda}, \mathcal{B}_{{\Lambda}})$. There exists a family of conditional probability measures $\set{P_{\ell}}$ on $\set{(C_{\ell},\mathcal{B}_{C_\ell})}$
giving the disintegration,
\begin{equation}\label{eq:finaldisintegration}
	P_{{\Lambda}}(A)  = \int_{\pi_{\mathcal{L}}(A)} \bigg( \int_{\pi^{-1}_{\mathcal{L}}(\ell) \cap A} \, dP_{\ell}(\lambda) \bigg) dP_{\mathcal{L}}(\ell), \ \forall A\in\mathcal{B}_{\Lambda}.
\end{equation}
\end{thm}

In the above Disintegration Theorems, the conditional probability measures can be extended to $(\Lambda,\mathcal{C}_{\Lambda})$ or $(\Lambda,\mathcal{B}_{\Lambda})$ by extending the measures $P_{\ell}(A)=0$ for all $A\subset \Lambda \backslash C_{\ell}$ for each $\ell\in\mathcal{L}$.  Theorem~\ref{thm:ProbabililtyDisintegration} guarantees that any probability measure on $(\Lambda,\mathcal{B}_{\Lambda})$ can be decomposed into a form involving a probability measure on $(\mathcal{L},\mathcal{B}_{\mathcal{L}})$ uniquely defined by $P_{\mathcal{D}}$ and probability measures on each measurable generalized contour space $(C_{\ell},\mathcal{B}_{C_{\ell}})$ defined by the conditional probabilities $P_{\ell}$. It is not obvious what these conditional probability measures are in Theorem~\ref{thm:ProbabililtyDisintegration}.
Clearly, any conditional probability measures on $\set{(C_{\ell},\mathcal{B}_{C_{\ell}})}$ can not be determined by observations of $Q(\lambda)\in\mathcal{D}$. This motivates the adoption of an {\em Ansatz} in which the probability measures along $\set{(C_{\ell},\mathcal{B}_{C_{\ell}})}$ are specified.
\smallskip

\begin{center}
\fbox{
\parbox{.9\linewidth}{
{\bf Ansatz: \ } For all $\ell \in \mathcal{L}$, assume a probability measure $P_{\ell}(\cdot)$ is given on $(C_{\ell},\mathcal{B}_{C_{\ell}})$.
}}
\end{center}

\smallskip

\begin{thm}\label{thm:Lambda_inverse}
Under the Ansatz, the stochastic inverse problem has a unique solution on $(\Lambda,\mathcal{B}_{\Lambda})$.
\end{thm}

We prefer a specific Ansatz that can be interpreted as prescribing a ``non-probabilistic'' or ``non-preferential'' weighting determined by the disintegration of volume measure to compute probabilities of events inside of a contour event\footnote{This is valid when $\mu_{\Lambda}(\Lambda)<\infty$, which is generally satisfied in practice using compact domains. If this is not the case, implying certain physical parameters are unbounded, then we may employ techniques similar to those used in Bayesian analysis for ``non-informative'' priors.}. However, the approximation method and resulting Algorithm~\ref{Alg_old} (summarized below) can be easily modified for any Ansatz (see \cite{Butler2014a} for more details).
\smallskip
\begin{center}
\fbox{
\parbox{.9\linewidth}{
{\bf Standard Choice for Ansatz: \ }
$P_{\ell} = \mu_{C_{\ell}}/\mu_{C_{\ell}}(C_{\ell}), \ \forall\ell\in\mathcal{L}. $
}}
\end{center}

If we do not specify an Ansatz, the stochastic inverse problem can only be solved on $(\Lambda,\mathcal{C}_{\Lambda})$ where there is only one possible probability measure on each $(C_{\ell},\mathcal{T}_{C_{\ell}})$. From Theorems~\ref{thm:ContourSigmaAlgebraDecomp} and \ref{thm:OriginalSigmaAlgebraDecomp}, it is evident that the choice of $\sigma-$algebra on $\Lambda$ dictates the SACOM and whether or not we must adopt an Ansatz for the solution.
%This is equivalent to substituting the measure space $(\Lambda,\mathcal{C}_{\Lambda})$ for $(\Lambda,\mathcal{B}_{\Lambda})$ in the Disintegration Theorem, which then describes a decomposition of a probability measure on $(\mathcal{L},\mathcal{B}_{\mathcal{L}})$ and the family of measurable contours $\set{(C_{\ell},\mathcal{T}_{C_{\ell}})}$.

\subsection{Approximating solutions}\label{S:ComputingDetails}

In \cite{Butler2014a}, we describe several approximation issues that need to be addressed in any practical computation of $P_{\mathbf{\Lambda}}$.  As described below, the fundamental approximation issues of any measure are the approximation of events in the various $\sigma-$algebras. The numerical approximation issues involve error in numerical evaluation of the model and computation of probability measures on some collection of events. Here, we discuss the event approximations. In Section~\ref{S:Errors}, we analyze the effect of errors on the approximate probability measures computed from the approximating sets.

%The solutions we approximate are probability measures absolutely continuous with respect to the underlying volume measures. Thus, there exist Radon-Nikodym derivatives, which in the context of probability theory are referred to as probability densities. We assume any probability density on $\mathcal{D}$ is continuous $\mu_{\mathcal{D}}$-a.e.

%In this paper, we use samples in $\Lambda$ to define the approximating events for the $\sigma-$algebras of $\Lambda$. In Section~\ref{S:StochGeom}, we use results from stochastic geometry/point processes to prove convergence to the exact probability measure $P_{\Lambda}$ in the absence of numerical error. However, numerical error can pollute the computations by  incorrectly identifying a contour event such that refinement of approximations of events in $\mathcal{B}_{\Lambda}$ fails to improve results. In Section~\ref{S:Errors}, we decompose the sources of error in the probability measure to account for the effect of finite samples and numerical accuracy separately.

\paragraph*{Brief review of approximations}

%Assuming that some algebra of elementary sets used to generate a $\sigma-$algebra can be defined such that the volume measure of any event from the $\sigma-$algebra  can be computed arbitrarily well from the outer-measure using finite unions of elements from this algebra,
Repeated application of the Lebesgue Dominated Convergence Theorem yields the following,
\begin{thm}\label{thm:original_approx}
Given probability measure $P_{\mathcal{D}}$, absolutely continuous with respect to $\mu_{\mathcal{D}}$, on $(\mathcal{D},\mathcal{B}_{\mathcal{D}})$ with density $\rho_{\mathcal{D}}$ and event $A\in\mathcal{B}_{\Lambda}$, there exists a sequence of approximations $P_{\Lambda,N}(A)$ using simple function approximations to probability densities $\rho_{{\Lambda},N}$ and $\rho_{\mathcal{D},M}$ requiring only calculations of volumes in $\Lambda$ that converges to $P_{\Lambda}(A)$ as $N,M\to\infty$.
\end{thm}

The proof of Theorem~\ref{thm:original_approx} details approximating probability densities $\rho_{\Lambda}$ and $\rho_{\mathcal{D}}$ (i.e.~the Radon-Nikodym derivatives of the probability measures) in order to apply the Lebesgue Dominated Convergence Theorem and outlines the computational measure-theoretic Algorithm~\ref{Alg_old} (see \cite{Butler2014a} for more details).  Unsurprisingly, the first step is the generation of partitions $\set{I_i}_{i=1}^M\subset\mathcal{D}$ and $\set{b_j}_{j=1}^N\subset\Lambda$ such that (1) arbitrary events in $\mathcal{B}_{\mathcal{D}}$ and $\mathcal{B}_{\Lambda}$ are approximated by unions of sets from these partitions, and (2) simple function approximations to the densities can be computed on these partitions. These partitions can be chosen from an algebra of elementary sets.

Since densities are $L^1$ functions, we can use convex sets with continuous boundaries to define a finite partition $\set{I_i}_{i=1}^M\subset\mathcal{D}$ and associated simple function approximation $\rho_{\mathcal{D},M}$ of sufficient accuracy in the $L^1$-norm\footnote{See Theorem 2.41 in \cite{Folland_Book}.}. The local differentiability of $Q$ implies that $Q^{-1}(I_i)$ is a measurable event in both $\mathcal{C}_{\Lambda}$ and $\mathcal{B}_{\Lambda}$ with boundary of zero $\mu_{\Lambda}$-measure. Such events can have their measures approximated within any prescribed tolerance by a finite set of hypercubes\footnote{See Theorem 2.40 of \cite{Folland_Book}.}. This defines an obvious choice for the cells $\set{b_j}_{j=1}^N$ partitioning $\Lambda$.
Once the probabilities of cells $\set{b_j}\subset\Lambda$ have been approximated by Algorithm~\ref{Alg_old}, we may estimate $P_{\Lambda}(A)$ for arbitrary event $A\in\mathcal{B}_{{\Lambda}}$ in any of the usual measure-theoretic ways such as using inner or outer sums, averages of inner and outer sums, or direct integration of the simple function approximation $\rho_{\Lambda,N}$.
%\footnote{If $\mu_{\Lambda}(\partial A)\neq 0$, then only direct integration of the simple function approximation $\rho_{\Lambda,N}$ will yield an approximation to $P_{\Lambda}(A)$ that converges as $N,M\to 0$.}.

Suppose instead that we wish to solve the stochastic inverse problem on $(\Lambda,\mathcal{C}_{\Lambda})$ instead of $(\Lambda,\mathcal{B}_{\Lambda})$. In this case,  we can approximate the solution using the same approximate solution on $(\Lambda,\mathcal{B}_{\Lambda})$ computed using Algorithm~\ref{Alg_old} since $\mathcal{C}_{\Lambda}\subset\mathcal{B}_{\Lambda}$. In other words, we may use the same approximate finite generating set $\set{b_j}$ for $\mathcal{B}_{\Lambda}$ as an approximate finite generating set for $\mathcal{C}_{\Lambda}$.

Unless otherwise stated, we let $\set{I_i}_{i=1}^M$ denote a partition of $\mathcal{D}$ used to approximate events in $\mathcal{B}_{\mathcal{D}}$ and define the simple function approximation $\rho_{\mathcal{D},M} = \sum_{i=1}^M p_i \mathbf{1}_{I_i}(q)$, where $p_i = P_{\mathcal{D}}(I_i)$. In other words, the probability on each set $I_i\subset\mathcal{D}$ is defined by the expected value of $\mathbf{1}_{I_i}(q)$ on $\mathcal{D}$. This is a natural choice equivalent to using the Integral Mean Value Theorem to define a normalized approximation, i.e.~$\rho_{\mathcal{D},M}$ is a density with this choice of $p_i$.

\begin{algorithm}\caption{Approximation of the Inverse Density}\label{Alg_old}
\begin{algorithmic}
\State Generate partitions $\set{I_i}_{i=1}^M\subset\mathcal{D}$ and $\set{b_j}_{j=1}^N\subset\Lambda$
\State Fix and normalize the simple function approximation $\rho_{\mathcal{D},M} =\sum_{i=1}^M p_i \mbf{1}_{I_i}(q)$
\State Let $\set{A_i}_{i=1}^M\subset\Lambda$ denote the induced regions of generalized contours partitioning ${\Lambda}$
\For{$j=1,\ldots,N$}
	\For{$i=1,\ldots,M$}
		\State Compute $\mu_{\Lambda}(A_i\cap b_j)$ and store as $ij$-component in matrix $V$.
	\EndFor
	
\EndFor
\For{$j=1,\ldots,N$}
	\State Set $P_{\Lambda,N}(b_j)$ to $\sum_{i=1}^{M} p_i(V_{ij}/\sum_{j=1}^N V_{ij})$
\EndFor
\end{algorithmic}
\end{algorithm}

%\begin{algorithm}\caption{Approximation to Exact Parameter Probability Distribution}\label{Alg_old2}
%\begin{algorithmic}
%\State Given map $q:{\Lambda}\to\mathcal{D}$ and family of probability laws $P_{\bar{Q}}$
%\State Fix and normalize simple function approximation, $\rho_{\mathcal{D},M}(\bar{Q})=\sum_{i=1}^M p_i \mbf{1}_{I_i}(\bar{Q})$, to $\rho_{\mathcal{D}}(q)$
%
%\State Let $\set{A_i}_{i=1}^M\subset\Lambda$ denote induced region of generalized contours partitioning ${\Lambda}$
%
%\State Fix simple function approximations $\sum_{k=1}^{K_i} u_{ik}\mbf{1}_{A_i\cap {\cal E}_{k}}(\lambda)$, to conditional densities defined by $P_{\bar{Q}}$ on each $A_i$.
%
%\State Partition ${\Lambda}$ with geometrically simple cells $\set{b_j}_{j=1}^{N}$
%\For{$j=1,\ldots,N$}
%	\For{$i=1,\ldots,M$}
%		\For{$k=1,\ldots,K_i$}
%			\State Compute $\mu_{\Lambda}(b_j\cap {\cal E}_k)$ and store in the $jk$-component of matrix $V^{(i)}$
%		\EndFor
%	\EndFor
%\EndFor
%\For{$j=1,\ldots,N$}
%		\State Set $P(b_j)$ to $\sum_{i=1}^M p_i\left[\sum_{k=1}^{K_i} u_{ik} \left(V_{jk}^{(i)}/\sum_{j=1}^N V_{jk}^{(i)}\right)\right]$
%\EndFor
%
%\end{algorithmic}
%\end{algorithm}

\paragraph*{A sample based approximation and counting measure}

Previous implementations of Algorithm~\ref{Alg_old} used regular grids so that $\set{b_j}_{j=1}^N$ was defined as a set of generalized rectangles or hypercubes \cite{Breidt2011, Butler2012a, Butler2013a}. While such sets may be refined and the measure of unions of these sets made to approximate any Borel set arbitrarily well \cite{Folland_Book}, there are obvious practical difficulties with using such an approximating set in high dimensions. Here, we consider an alternative where the sets $\set{b_j}_{j=1}^N$ are defined implicitly by a finite collection of samples in $\Lambda$ satisfying some particular properties detailed below.

\begin{defn}
For a fixed number of samples $\set{\lambda^{(j)}}_{j=1}^N\subset\Lambda$, there is a \textbf{Voronoi tessellation} of ${\Lambda}$ denoted by $\set{\mathcal{V}_j}_{j=1}^{N}\subset\Lambda$ defined by
\[
	\mathcal{V}(\lambda^{(j)}) :=  \set{\lambda\in{\Lambda} \, : \, d_v(\lambda^{(j)},\lambda) \leq d_v(\lambda^{(i)},\lambda), \forall i=1,\ldots,N}.
\]
Here, $d_v(\cdot,\cdot)$ denotes a metric on ${\Lambda}$ used to define the Voronoi cells\footnote{The metric $d_v(\cdot,\cdot)$ is possibly different from the metric that induces the volume measure $\mu_{{\Lambda}}$ and Borel $\sigma$-algebra $\mathcal{B}_{{\Lambda}}$. Common choices for $d_v(\cdot,\cdot)$ are the standard Euclidean 2-norm or 1-norm.}.
\end{defn}

Clearly the goal is to define a set of samples $\set{\lambda^{(j)}}_{j=1}^N$ implicitly defining a tessellation of $\Lambda$ that is useful for approximating events, i.e.~measurable sets, in $\mathcal{B}_{\Lambda}$. We often approximate events in a given $\sigma-$algebra by a Voronoi coverage.

\begin{defn}
We say that $A_N$ is the \textbf{Voronoi coverage} of $A\in\mathcal{B}_{\Lambda}$ and $\mu_k(A)$ its volume defined by
\[
	A_N:=\bigcup_{\lambda^{(j)}\in A, 1\leq j\leq N} \mathcal{V}(\lambda^{(j)}), \ \text{and} \ \mu_k(A) = \mu_{\Lambda}(A_k) = \sum_{j=1}^{N} \mu_{\Lambda}(\mathcal{V}(\lambda^{(j)}))\mathbf{1}_{\lambda^{(j)}\in A}.
\]
\end{defn}

\begin{defn}
A rule for defining any $N$ samples $\set{\lambda^{(j)}}_{j=1}^N\subset\Lambda$ is called \textbf{$\mathcal{B}_{\Lambda}$-consistent} if
\[
	\mu_{\Lambda}(A\triangle A_k) \to 0 \text{ as }N\to\infty, \ \forall A\in\mathcal{B}_{\Lambda}, \text{s.t. } \mu_{\Lambda}(\partial A) = 0.
\]
\end{defn}

Clearly any $\mathcal{B}_{\Lambda}$-consistent rule implies that any generalized rectangle or finite union of generalized rectangles may be approximated arbitrarily well by a Voronoi coverage defined by a finite number of samples. Replacing $\set{b_j}_{j=1}^N$ with a Voronoi tesselation $\set{\mathcal{V}_j}_{j=1}^N$ generated from a $\mathcal{B}_{\Lambda}$-consistent rule in Algorithm~\ref{Alg_old}, we define the following,
\begin{defn}
Let $\tilde{P}_{{\Lambda},N}$ denote the \textbf{counting probability measure} (or simply \textbf{counting measure} on $({\Lambda},\mathcal{B}_{{\Lambda}})$ defined by
\[
	\tilde{P}_{{\Lambda},N}(A) := \sum_{j=1}^N P_{\Lambda,N}(\mathcal{V}_j) \mbf{1}_{\lambda^{(j)}}(A).
\]
\end{defn}

Note that $\tilde{P}_{{\Lambda},N}$ is simply a measure of point masses at each sample from $\set{\lambda^{(j)}}_{j=1}^N$. A technical point is that we may use a counting measure to estimate the probability of {\em any} $A\in\mathcal{B}_{\Lambda}$. However, we only prove convergence to the exact probability when $\mu_{\Lambda}(\partial A) = 0$.

It is possible to have compact Borel sets with boundaries of non-zero measure. We describe in Section~\ref{S:non-intrusive-converge} how to approximate the probability of such sets.
%On the other hand, we can show that for any $\epsilon>0$, there exists a finite collection of generalized rectangles $\set{R_j}_{j=1}^N$ with subset $\set{R_{n_j}}$ such that $P_{\Lambda,N}(A\triangle \union_{n_j} R_{n_j}) <\epsilon$. For this
In practical computations, we are not concerned with such sets with fractal boundaries. Furthermore, such sets do not enter into any of the computational algorithms. This is due to the assumption of local differentiability of the map $Q$, which implies that all generalized contours, and subsequently, all induced regions of generalized contours defined by $Q^{-1}(I)$ for any Borel set $I\in\mathcal{B}_{\mathcal{D}}$ are sets with boundaries of zero $\mu_{\Lambda}$-measure.

\paragraph*{A Monte Carlo approach}
There are many possible $\mathcal{B}_{\Lambda}$-consistent rules that we may choose. For example, we may define a sequence of uniformly refined grids of which the grid points are numbered and sequentially sampled in a serpentine manner. This choice produces Voronoi cells that are generalized rectangles (for certain choices of $N$) and has been used previously, e.g.~see \cite{Breidt2011,Butler2012a,Butler2013a}. While this produces small Voronoi cells everywhere in the domain, the approximation properties of the cells have negative dependence on the dimension of the parameter space in terms of requiring large $N$ to achieve reasonable $\mu_{\Lambda}$-approximations by the Voronoi coverage.

We may choose to normalize the volume measure (if possible) and generate N i.i.d.~samples from this distribution corresponding to a ``uniform'' sampling density (see Section~\ref{S:StochGeom} below). This amounts to a standard Monte Carlo (MC) approach for sampling on $(\Lambda,\mathcal{B}_{\Lambda},\mu_{\Lambda})$.
% and permits a straightforward approximation to $P_{\Lambda,N}$ where $P_{\Lambda,N}(b_j)$ is replaced by $P_{\Lambda,N}(\mathcal{V}_j)$ in Algorithm~\ref{Alg_MC} below, where the standard MC approximation that all
Furthermore, we may use the standard MC approximation that all Voronoi cells have the same volume which greatly reduces the computational demands of the algorithm, see Algorithm~\ref{Alg_MC}.
% in the last for-loop where we divide by the number of cells identified as approximating an induced region of generalized contours.

\begin{algorithm}\caption{A Monte Carlo Approximation of the Inverse Density}\label{Alg_MC}
\begin{algorithmic}
\State Let $\set{\lambda^{(j)}}_{j=1}^N\subset{\Lambda}$ denote uniform i.i.d.~random samples from $\mathcal{B}_{\Lambda}$-consistent rule, and $\set{\mathcal{V}_j}_{j=1}^{N}\subset\Lambda$ denote the associated Voronoi tessellation of ${\Lambda}$.
\For{$j=1,\ldots,N$}
\State Assign a nominal value of $Q_j$ to $\mathcal{V}_j$, e.g.~in the continuous case use $Q_j=Q(\lambda^{(j)})$.
\EndFor
\State Generate partition $\set{I_i}_{i=1}^M\subset\mathcal{D}$.
\State Fix  and normalize the simple function approximation $\rho_{\mathcal{D},M} =\sum_{i=1}^M p_i \mbf{1}_{I_i}(q)$.
\State Initialize $M\times 1$ counting vector $\mathbf{c}$ and $N\times 1$ pointer vector $\mathbf{i_o}$ to zeros.
\For{$j=1,\ldots,N$}
	\State Set $i=1$ and \textbf{flag}$=0$
	\While{$i\leq M$ and \textbf{flag}$=0$}
		\If{ $Q_j\in I_i$} $\mathbf{c}(i) = \mathbf{c}(i)+1$, $\mathbf{i_o}(j)=i$, \textbf{flag}$=1$.
		\Else{ $i=i+1$}.
		\EndIf
	\EndWhile
	
\EndFor
\For{$j=1,\ldots,N$}
	\State Set $P_{\Lambda,N}(\mathcal{V}_j)$ to $p_{\mathbf{i_o}(j)}/\mathbf{c}(\mathbf{i_o}(j))$
\EndFor
\end{algorithmic}
\end{algorithm}

%Another choice is taking a (non-random) uniform grid of $N$ points as in \cite{Breidt2011,Butler2012a}. The Voronoi tessellation is a uniform rectangular grid in this case. That provides small cells everywhere in the domain, but the approximation properties of the cells have negative dependence on the dimension of the parameter space.

The counting measure definition is consistent with the MC approach of Algorithm~\ref{Alg_MC} to estimate the probabilities of individual Voronoi cells $\mathcal{V}_j$ inside of a particular induced region of generalized contours by counting the number of samples within this event. Specifically, in the algorithm, $\mathbf{c}(\mathbf{i_o}(j))$ counts the number of Voronoi cells we associate within the induced generalized region of contours defined by $Q^{-1}(I_{\mathbf{i_o}(j)})$.

Standard MC sampling algorithms tend to produce random samples that are clustered resulting in large variations in the sizes of the corresponding Voronoi cells. This may not provide good set-approximation properties for small sample sizes. Sampling points equidistant along a space filling curve or using Latin hypercube sampling are other options. Or, we may simply choose the samples deterministically based on other information we know or we want to impose a certain resolution on the problem, e.g.~we may wish to obtain a fine resolution of some particular subset of $\Lambda$ for some design purposes. Some pseudo-MCMC sampling might be best. The underlying theory only requires that the samples be generated from a $\mathcal{B}_{\Lambda}$-consistent rule. Exploring all possibilities is beyond the scope of this paper and we limit presentation to the basic theory of convergence of counting measures for $\mathcal{B}_{\Lambda}$-consistent rules with special attention paid to random sampling rules satisfying criterion as discussed in Section~\ref{S:StochGeom} below.

\section{Convergence of counting measures}\label{S:Theory}

\subsection{General theory}\label{S:non-intrusive-converge}

%The proof of convergence uses the almost sure convergence in measure of Voronoi coverages and explicitly uses generalized rectangles (i.e. hypercubes) intersecting the domain as in the proof of Lemma~\ref{lem:SLLN}. This is a classical measure-theoretic argument based upon simple function approximations defined using some approximating generating set to a Borel $\sigma-$algebra. We use a triangle inequality to separate the error in probability into two parts involving the approximation properties of the Voronoi coverages to hypercubes and the approximation properties of the hypercubes to more general Borel-measurable sets.

The proof of convergence is based upon classical measure-theoretic arguments using simple function approximations defined using some approximate generating set to a Borel $\sigma-$algebra. We use a triangle inequality to separate the error in probability into two parts involving the approximation properties of the Voronoi coverages to generalized rectangles and the approximation properties of the generalized rectangles to other Borel-measurable sets.

\begin{lemma}\label{lem:help}
For any $A\in\mathcal{B}_{\Lambda}$ and $\epsilon>0$, there exists simple function approximation $\rho_{\mathcal{D},M}$ and $K<\infty$ generalized rectangles in $\mathbb{R}^n$ partitioning compact ${\Lambda}$ with $\abs{P_{\Lambda,K}(A)-P_{\Lambda}(A)}<\epsilon$.
\end{lemma}
%The family of all generalized rectangles in $\mathbb{R}^n$ partitioning compact $\Lambda$ forms a generating set of the Borel $\sigma-$algebra $\mathcal{B}_{\Lambda}$. Furthermore, any event in $\mathcal{B}_{\Lambda}$ can have its $\mu_{\Lambda}$-measure approximated arbitrarily well by taking sufficiently fine generalized rectangles.

{\em Proof:}
Let $\epsilon>0$ be given.
There exists a partition of $\mathcal{D}$ into $M$ generalized rectangles, $\set{I_i}_{i=1}^M$, such that the $L^1$ error of $\rho_{\mathcal{D},M}$ and $\rho_{\mathcal{D}}$ is less than $\epsilon/3$. For each $A_i := Q^{-1}(I_i)\in\mathcal{B}_{\Lambda}$, which have boundaries of zero $\mu_{\Lambda}$-measure, there is a finite number $K_i$ of generalized rectangles partitioning $\Lambda$ such that $\abs{P_{\Lambda,K_i}(A_i)-P_{\Lambda}(A)}<\epsilon/(3M)$. Moreover, for any $A\in\mathcal{B}_{\Lambda}$ there is a finite number $K_A$ of generalized rectangles, $\set{R_j}_{j=1}^{K_A}$ in $\Lambda$ such that $\mu_{\Lambda}(A\triangle \bigcup_{j} R_j) <\epsilon/3$. The conclusion follows from Theorem~\ref{thm:original_approx} where $P_{\Lambda,K}$ is constructed as in Algorithm~\ref{Alg_old} with $\set{b_j}_{j=1}^K$ given by $K$ sufficiently fine generalized rectangles partitioning $\Lambda$ from which the finite number of generalized rectangles used above may be constructed by finite unions. $\Box$

% Since the hypercubes in $\mathbb{R}^n$ intersecting $\Lambda$ can be used to form a generating set for the Borel $\sigma-$algebra $\mathcal{B}_{\Lambda}$, we can choose the $N$ hypercubes partitioning $\Lambda$, denoted by $\set{b_j}_{j=1}^N$, such that $\mu_{\Lambda}(A_k\delta A) < \eta_1$ for any $\eta_1>0$. The Lebesgue Dominated Convergence Theorem implies there exists simple function $\rho_{\mathcal{D},M}$ of sufficient accurac and using Theorem~\ref{thm:original_approx}

%\begin{thm}\label{thm:counting_converges}
%Suppose $Q:{\Lambda}\to\mathcal{D}$ is sufficiently smooth with GD component maps. Given probability measure $P_{\mathcal{D}}$ absolutely continuous with respect to $\mu_{\mathcal{D}}$ with density $\rho_{\mathcal{D}}$, sampling distribution $F$ absolutely continuous with respect to $\mu_{{\Lambda}}$ with density $f(\lambda)>0$ for almost every $\lambda\in{\Lambda}$, then for any event $A\in\mathcal{B}_{{\Lambda}}$, there exists a sequence of approximations of simple functions $\rho_{\mathcal{D},M}$ and counting measures $\tilde{P}_{\Lambda,N}(A)$ such that $\tilde{P}_{{\Lambda},N}(A)\to P_{{\Lambda}}(A)$ almost surely as $N,M\to\infty$.
%\end{thm}

\begin{thm}\label{thm:counting_converges}
Given probability measure $P_{\mathcal{D}}$, absolutely continuous with respect to $\mu_{\mathcal{D}}$, on $(\mathcal{D},\mathcal{B}_{\Lambda})$ with density $\rho_{\mathcal{D}}$ and a $\mathcal{B}_{\Lambda}$-consistent rule for generating samples, then for any event $A\in\mathcal{B}_{{\Lambda}}$ with $\mu_{\Lambda}(\partial A)=0$, there exists a sequence of approximations of simple functions $\rho_{\mathcal{D},M}$ and counting measures $\tilde{P}_{\Lambda,N}(A)$ such that $\tilde{P}_{{\Lambda},N}(A)\to P_{{\Lambda}}(A)$ a.s.~as $N,M\to\infty$.
\end{thm}

{\em Proof:}
Let $\epsilon>0$ be given and $\lambda^{(1)},\lambda^{(2)},\ldots$ denote a sequence of samples in ${\Lambda}$ computed from the $\mathcal{B}_{\Lambda}$-consistent rule. We have by the assumptions of $Q$ that for any fixed $I_i\in\mathcal{B}_{\mathcal{D}}$ with continuous boundary that the induced region of generalized contours $A_i:=Q^{-1}(I_i)\in\mathcal{B}_{\Lambda}$ has boundary with zero $\mu_{{\Lambda}}$-measure.
By assumption of $\rho_{\mathcal{D}}$, there exists a sequence of partitions $\set{I}_{i=1}^M\subset\mathcal{D}$ such that for any fixed but arbitrary $M$ and $I_i\in\mathcal{B}_{\mathcal{D}}$ we have that $\mu_{{\Lambda}}(A_{i,N} \triangle A_i) \to 0$ as $N\to\infty$ where $A_{i,N}$ denotes the Voronoi coverage of $A_i$. By construction of $P_{\Lambda,N}(\mathcal{V}_i)$, the definition of $\tilde{P}_{{\Lambda},N}$, and the fact that the Borel sets on $\Lambda$ can be generated by the collection of all generalized rectangles, there is a sufficiently fine partition of $\Lambda$ by $K_1$ generalized rectangles and an $N_1$ such that $\abs{\tilde{P}_{\Lambda,N}(A_i)-P_{\Lambda,K_1}(A_i)}<\epsilon/2$ for all $N>N_1$.

By Lemma~\ref{lem:help}, there is a sufficiently fine partition of $\Lambda$ by $K_2$ generalized rectangles such that $\abs{P_{\Lambda,K_2}(A_i)-P_{\Lambda}(A_i)}<\epsilon/2$. Let $K$ denote the number of generalized rectangles needed to generate all of the above generalized rectangles by finite unions and intersections. Then by a standard triangle inequality  and Lemma~\ref{lem:help}, we have that for partition $\set{I_i}_{i=1}^M$ as above defining simple function approximation $\rho_{\mathcal{D},M}$ to $\rho_{\mathcal{D}}$ that $\abs{\tilde{P}_{\Lambda,N}(A_i)-P_{\Lambda}(A_i)}<\epsilon$ for all $N>N_1$ for each $1\leq i\leq M$. The result extends to any $A\in\mathcal{B}_{\Lambda}$ with $\mu_{\Lambda}(\partial A)=0$ by applying the Lebesgue Dominated Convergence Theorem to a sequence of $\rho_{\mathcal{D},M}$ constructed from a sequence of  partitions $\set{I_i}_{i=1}^M\in\mathcal{B}_{\mathcal{D}}$ of generalized rectangles in $\mathcal{D}$ and using Theorem~\ref{thm:Lambda_inverse}. $\Box$

We can approximate the probability of general $A\in\mathcal{B}_{\Lambda}$ using the counting measure with error less than any $\epsilon>0$. For such an $\epsilon>0$, we use Lemma~\ref{lem:help} with $\epsilon/2$ to define a finite union of generalized rectangles in $\Lambda$  approximating $A$. Let $B$ denote this finite union of generalized rectangles. Then $B$ has the property that $\mu_{\Lambda}(\partial B) = 0$ and Theorem~\ref{thm:counting_converges} applies to the set $B$. Let $B_N$ denote the Voronoi coverage of $B$ with $N$ samples in $\Lambda$ satisfying $\abs{\tilde{P}_{N,\Lambda}(B)-P_{\Lambda}(B)}<\epsilon/2$ for any $\epsilon>0$. Then using $(A\triangle B_N) = (A\triangle B) \triangle (B\triangle B_N)$ and a triangle inequality, the result follows. This proves the following
\begin{thm}\label{thm:counting_approx}
Given probability measure $P_{\mathcal{D}}$, absolutely continuous with respect to $\mu_{\mathcal{D}}$, on $(\mathcal{D},\mathcal{B}_{\Lambda})$ with density $\rho_{\mathcal{D}}$ and a $\mathcal{B}_{\Lambda}$-consistent rule for generating samples. For any event $A\in\mathcal{B}_{{\Lambda}}$ and $\epsilon>0$, there exists simple function $\rho_{\mathcal{D},M}$, finite $N$ samples defining counting measure $\tilde{P}_{\Lambda,N}$, and a set $B\in\mathcal{B}_{\Lambda}$ with $\mu_{\Lambda}(\partial B)=0$ such that $\tilde{P}_{\Lambda,N}(B\triangle A) <\epsilon$.
\end{thm}

\subsection{Convergence of the Monte Carlo approximation}\label{S:StochGeom}

%The theory of convergence of the Monte Carlo approximation to the measure-theoretic inverse makes explicit use of the classic Implicit Function Theorem for multivariate functions and a result from stochastic geometry using Poisson processes that can be
We present a result that can be viewed as a strong law of large numbers for Voronoi tessellations. %We summarize the terminology and required theory of Poisson point processes below and
We direct the interested reader to \cite{Penrose, Khmaladze} for a more thorough exposition on the subject within the context of stochastic geometry and point processes. % of which we provide a more complete proof below.

%\subsubsection{Results from stochastic geometry}

\begin{defn}
A \textbf{sampling distribution}, denoted by $F$, is any distribution absolutely continuous with respect to $\mu_{{\Lambda}}$ such that the corresponding \textbf{sampling density} $f(\lambda)>0$ for almost every $\lambda\in{\Lambda}$.
\end{defn}

%Let $F$ denote a distribution absolutely continuous with respect to $\mu_{{\Lambda}}$ such that the corresponding density $f(\lambda)>0$ for all $\lambda\in{\Lambda}$. The necessity of $f(\lambda)>0$ is made clear below. Suppose $\lambda^{(1)},\lambda^{(2)},\ldots$ are i.i.d.~random variables in ${\Lambda}$ with distribution $F$. Let $N_k$ for $k=1,2,\ldots$ denote a sequence of Poisson$(k)$ random variables independent of the $\lambda^{(i)}$ samples\footnote{In practice, we simple use a monotonic sequence for the number of samples. However, restricting  the number of samples to increase monotonically is not required for convergence so we use the weaker assumption.}. We say that
%\[
%	\xi_k(A) = \sum_{j=1}^{N_k} \mathbf{1}_{\set{\lambda^{(j)}\in A}}, \ k= 1,2,\ldots
%\]
%%is a Poisson process on the Borel $\sigma$-algebra $\mathcal{B}_{{\Lambda}}$ of measurable subsets of ${\Lambda}$ with $E(\xi_k(A)) = kF(A)$ for measurable $A\in\mathcal{B}_{{\Lambda}}$. Note that a Poisson process is defined for each Poisson random variable $N_k$, that $N_k\to\infty$ almost surely as $k\to\infty$, and that for $N_{k}>N_{l}$, $\xi_k(A)\geq \xi_l(A)$.
%The quantity $nf(\lambda)$ defines the intensity of the Poisson process.
% For a fixed $N_k$, we have a corresponding Voronoi tessellation of ${\Lambda}$ associated with the Poisson process. We use this tessellation to define approximating events to any event in $\mathcal{B}_{\Lambda}$.

The following Lemma from \cite{Khmaladze} has been modified to conform to our notation. To our knowledge, \cite{Khmaladze} was the first to summarize and explicitly prove the fundamental elements of Lemma~\ref{lem:SLLN}. Below, we modify and expand the proof presented in \cite{Khmaladze} to highlight details of the convergence and uniqueness that are of practical use, e.g., when considering design of adaptive sampling procedures to determine more accurate Voronoi coverages of specific events \cite{Butler2014b} or in understanding issues related to a posteriori error estimates described below in Section~\ref{S:Errors}.

\begin{lem}\label{lem:SLLN}
Given sampling density $f(\lambda)>0$ almost everywhere on ${\Lambda}$, if $A\in\mathcal{B}_{{\Lambda}}$ such that $\mu_{{\Lambda}}(\partial A) = 0$, then almost surely
\[
	\mu_{{\Lambda}}(A_N \triangle A) \to 0, \ N\to\infty.
\]
\end{lem}

{\em Proof.} Suppose $A\in\mathcal{B}_{{\Lambda}}$, $\mu_{{\Lambda}}(A)>0$ and $\mu_{{\Lambda}}(\partial A) = 0$. Let $A^{\delta} := (A+\delta B(0,1))\cap {\Lambda}$ denote the Minkowski sum of $A$ and $\delta B(0,1)$ restricted to domain ${\Lambda}$, where $B(0,1)$ is the unit ball in $\mathbb{R}^n$. Let $A_{\delta} = ((A^c)^{\delta})^c$.  Since $\mu_{{\Lambda}}(\partial A)=0$, for all $\epsilon>0$ there exists $\delta(\epsilon)>0$ such that $\partial A\subset A^{\delta(\epsilon)}\backslash A_{\delta(\epsilon)}$ and $\mu_{{\Lambda}}(A^{\delta(\epsilon)}\backslash A_{\delta(\epsilon)})<\epsilon$.

Let $r(\lambda^{(j)}):=\max_{\lambda\in\mathcal{V}(\lambda^{(j)})} d_v(\lambda,\lambda^{(j)})$ and $r_N:=\max_{1\leq j\leq N} r(\lambda^{(j)})$.  We now prove that $r_N\to 0$ almost surely. For any $h\in\mathbb{N}$, let $\mathcal{H}_h$ denote the finite set of hypercubes in $\mathbb{R}^n$ with edge-length $1/h$ partitioning compact ${\Lambda}$.
% i.e.~$H\in\mathcal{H}_h$ is defined by $\left\{\Pi_{i=1}^n[j_i/h,(j_i+1)/h]\right\}\cap {\Lambda}$ for some $j_i\in\mathbb{N}$ for each $i$.
%Since $h>0$ is finite,
By Kninchine's strong law of large numbers,
\[
	\max_{H\in \mathcal{H}_h} \abs{\frac{\sum_{j=1}^{N} \mathbf{1}_{\set{\lambda^{(j)}\in H}}}{N} - \int_H f(\lambda)\, d\mu_{{\Lambda}}(\lambda)} \to 0 \text{ almost surely as } N\to\infty
\]
and
\[
	\min_{H\in\mathcal{H}_h} \int_H f(\lambda)\, d\mu_{{\Lambda}}(\lambda)>0 \text{ for all } h.
\]
Suppose $r_N$ does not converge to zero almost surely, then there exists $\epsilon>0$ such that $r_N\geq \epsilon$ infinitely often. By construction, $r_N\geq r_K$ for any $N \geq K$, and there exists a fixed $j$ such that $r(\lambda^{(j)})\geq\epsilon$ for all $k\geq j$. Choosing a set of hypercubes intersecting $\Lambda$ with edge-length sufficiently small (and positive) such that the maximum distance in the $d_v$-metric between any two points in the hypercubes is less than $\epsilon/4$ yields a contradiction since  $\abs{\frac{\sum_{j=1}^{N} \mathbf{1}_{\set{\lambda^{(j)}\in H}}}{N} - \int_H f(\lambda)\, d\mu_{{\Lambda}}(\lambda)} \to 0$ almost surely~for each of these hypercubes, so there must exist a sample $\lambda^{(i)}$ with $i\neq j$ such that its distance is less than $\epsilon/2$ from $\lambda^{(j)}$.  Thus $r_N\to 0$ almost surely. It follows that, for all $\epsilon>0$, there exists a Voronoi coverage $A_N$ of $A$ such that $A_{\delta(\epsilon)}\subset A_N\subset A^{\delta(\epsilon)}$ for all $N>N(\epsilon)$ almost surely, which implies $\mu_{{\Lambda}}(A_N\triangle A)\to 0$ almost surely. $\Box$

%We note that Lemma~\ref{lem:SLLN} holds if we replace the Poisson process $\xi_k$ with random variable $\zeta_N$, where $\zeta_N = N$ with probability 1 and $N=k$.

This implies that any random sampling scheme with sampling density $f(\lambda)>0$ produces a $\mathcal{B}_{\Lambda}$-consistent rule with probability $1$. Theorem~\ref{thm:counting_converges} then applies to counting measures computed using random sampling as long as the sampling density is positive.

%$f(\lambda)>0$ to guarantee convergence and may actually improve numerical results.
\subsection{Choosing sampling densities}

If $f(\lambda)=0$ on any set $A\in\mathcal{B}$ such that $\mu_{{\Lambda}}(A)>0$, then any sequence of counting measures evaluated on $A$ will result in a sequence of zeros. If $A$ shares any volume with a region of induced generalized contours of non-zero probability, then according to Theorem~\ref{thm:ProbabililtyDisintegration}\footnote{With any Ansatz such that the generalized contours in $A$ have non-zero probability.}, $P_{\Lambda}(A)>0$ and the sequence of approximate probabilities fails to converge. This highlights the importance of the sampling density being strictly positive on $\Lambda$.

In Algorithm~\ref{Alg_MC}, a natural choice is to let $f(\lambda)=1$, i.e.~we sample from ${\Lambda}$ according to its volume measure. However, we may choose a non-uniform intensity to improve resolution or accuracy of the approximate probability measure for certain induced regions of generalized contours, e.g.~in locations where $\norm{\nabla Q}$ is largest. Such ideas of computational accuracy with a finite number of samples and adaptive sampling are the subjects of future work.

We note that in Algorithm~\ref{Alg_MC}, we use a standard MC approximation that can be interpreted as approximating the volumes of Voronoi cells to be equal. We may instead opt to more accurately estimate the volumes for each $\mathcal{V}_j$, in which case the algorithm more closely resembles Algorithm~\ref{Alg_old} where we explicitly use volumes of the partition of $\Lambda$ to obtain $P_{\Lambda,N}$. Similarly, if there exists $A\in \mathcal{B}_{\Lambda}$ such that $\mu_{\Lambda}(A)>0$ and $f(\lambda)\neq 1$ for $\lambda\in A$, then the standard MC approximation that all $\mu_{\Lambda}(\mathcal{V}_j)$ are equal for all $j$ no longer applies. The necessary modification is that  approximation to the ratio of volumes of $\mu_{\Lambda}(\mathcal{V}_j)/\sum_{\{k\, :\, \mathbf{i_o}(k)=\mathbf{i_o}(j)\}} \mu_{\Lambda}(\mathcal{V}_k)$ replace division by $\mathbf{c}(\mathbf{i_o}(j))$ in the last for-loop\footnote{In this case, the counting vector $\mathbf{c}$ can be entirely removed from Algorithm~\ref{Alg_MC}.}. Similar modifications are required for a non-standard choice of the Ansatz. Below, we assume such modifications to Algorithm~\ref{Alg_MC} are made as necessary to correctly compute the counting measure $\tilde{P}_{\Lambda,N}$.
\section{Sources of error in the counting measure}\label{S:Errors}

There are two types of error in approximating probability measure $P_{\Lambda}$ using Algorithm~\ref{Alg_MC}: stochastic and deterministic. Stochastic error arises from using $N$ {\em finite} samples $\set{\lambda^{(j)}}_{j=1}^N$ to implicitly define a Voronoi tessellation $\set{\mathcal{V}_j}$ of $\Lambda$. Deterministic error arises from the numerical evaluation of the map $Q$ for each of the $N$ samples. Each of these errors affects $\tilde{P}_{\Lambda,N}$ in different ways. The deterministic error may lead to misidentification of induced regions of generalized contours (possibly leading to incorrect component values of the vectors $\mathbf{c}$ and $\mathbf{i_o}$ in Algorithm~\ref{Alg_MC}). We let $\tilde{P}_{\Lambda,N,h}$ denote the computed counting measure using numerical computations $Q_h$ approximating the map $Q$, where $h$ denotes some numerical discretization parameter. The choice of samples defines the set-approximation properties of events in $\mathcal{B}_{\Lambda}$ in terms of unions and complements of Voronoi cells (the first step in Algorithm~\ref{Alg_MC}). In fact, the ability to estimate induced regions of generalized contours is limited by the choice of samples, so it is the stochastic source of error that we examine first.

Below, $P_{\Lambda}$ denotes the exact probability measure on $(\Lambda,\mathcal{B}_{\Lambda})$ for the standard Ansatz\footnote{The error analysis can be altered for a non-standard Ansatz where volume measures and counting are replaced by the probabilistic weighting of the contour events.} {\em given} a fixed simple-function approximation to $P_{\mathcal{D}}$. In practice, when repeated experiments and the subsequent measurements are used to empirically determine $P_{\mathcal{D}}$, such an approximation would be determined by the binning of the relevant QoI. In general, we decompose the error in the  computed probability measure $\tilde{P}_{\Lambda,N,h}$ as
\begin{equation}\label{eq:P_error_decompose}
	P_{\Lambda}(A) - \tilde{P}_{\Lambda,N,h}(A)  = \underbrace{\left(P_{\Lambda}(A)-\tilde{P}_{\Lambda,N}(A)\right)}_{\mathbf{I}} - \underbrace{\left(\tilde{P}_{\Lambda,N,h}(A)-\tilde{P}_{\Lambda,N}(A)\right)}_{\mathbf{II}}, \ \forall A\in \tilde{\mathcal{B}}_{\Lambda}\subset\mathcal{B}_{\Lambda}.
\end{equation}
Here, term $\mathbf{I}$ is the error in approximating the probability of event $A$ by the counting measure $\tilde{P}_{\Lambda,N}$.  Term $\mathbf{II}$ is the error in using the numerical map $Q_h$ to identify induced regions of generalized contours in Algorithm~\ref{Alg_MC}. The events $A$ belong to a $\sigma-$algebra $\tilde{\mathcal{B}}_{\Lambda}$ that is a subset of the original $\sigma-$algebra $\mathcal{B}_{\Lambda}$.

\subsection{Stochastic set-approximation and numerical $\sigma-$algebras}\label{S:stoch_error}

We find convenient the following,
\begin{defn}\label{def:numerical_algebra}
For a given set of samples on $\Lambda$, $\set{\lambda^{(j)}}_{j=1}^N$, let $\mathcal{B}_{\Lambda,N}$ denote the {\bf numerical $\sigma-$algebra} generated by the implicitly defined Voronoi tessellation $\set{\mathcal{V}_j}_{j=1}^N$\footnote{Borel $\sigma-$algebras on $\mathbb{R}$ are commonly generated using generating sets defined, for example, by all half-open intervals. Analogously, $\mathcal{B}_{\Lambda,N}$ is equivalently generated by the set of all Voronoi coverages for all $A\in\mathcal{B}_{\Lambda}$. This particular definition uses the {\em minimal} generating set.}.
\end{defn}

By the definition or Borel sets and Voronoi cells, $\mathcal{B}_{\Lambda,N}$ is a proper subset of $\mathcal{B}_{\Lambda}$ for any finite $N$, i.e.~any event in $\mathcal{B}_{\Lambda,N}$ is $\mu_{\Lambda}$-measurable.

Definition~\ref{def:numerical_algebra} is illustrative of the fundamental types of events we  may compute probabilities for using Algorithm~\ref{Alg_MC} with no set-approximation errors, i.e.~with no term $\mathbf{I}$ error. Specifically, for a fixed number of samples, a fixed approximation to $\rho_{\mathcal{D}}$, and given exact map $Q$ in Algorithm~\ref{Alg_MC}, we may compute $\tilde{P}_{\Lambda,N}(A)$ exactly for all $A\in\mathcal{B}_{\Lambda,N}$. From Lemma~\ref{lem:SLLN}, we have that any $A\in\mathcal{B}_{\Lambda}$ can be sufficiently approximated in $\mu_{\Lambda}$-measure by a Voronoi coverage almost surely, i.e.~we are almost surely guaranteed to find an element of $\mathcal{B}_{\Lambda,N}$ for some $N$ that approximates $A$ to the desired accuracy. This was exploited in the proof of Theorem~\ref{thm:counting_converges} that the counting measures converge. Here, we focus on the events $A\in\mathcal{B}_{\Lambda}$ of non-zero $\tilde{P}_{\Lambda,N}$-measure and the errors in the computed probabilities of such events.

Below, we show that for certain events in $\mathcal{B}_{\Lambda}$ that there are, in a probabilistic sense, optimal and unique approximations to these events in $\mathcal{B}_{\Lambda,N}$.

\begin{lemma}\label{lemma:nonempty}
If $B,C\in\mathcal{B}_{\Lambda,N}$ and $\mu_{\Lambda}(B\triangle C)>0$, then there exists at least one $\lambda^{(k)}$, $k\in\set{1,2,\ldots,N}$), such that $\lambda^{(k)}$ is in one and only one of the measurable sets $B$ or $C$.
\end{lemma}

{\em Proof:}
By definition, the $\sigma-$algebra $\mathcal{B}_{\Lambda,N}$ is constructed by the closure of complements and countable unions of the finite set of individual Voronoi cells $\set{\mathcal{V}_j}_{j=1}^N$. Thus, by construction, any sets of non-zero $\mu_{\Lambda}$-measure in $\mathcal{B}_{\Lambda,N}$ must contain at least the interior of a single Voronoi cell. Since $B\triangle C \in \mathcal{B}_{\Lambda,N}$, the conclusion follows. $\Box$

\begin{thm}\label{thm:stoch_error}
Suppose $\set{\lambda^{(j)}}_{j=1}^N$ is fixed and choose any $A\in\mathcal{B}_{\Lambda}$ such that there exists at least one $\lambda^{(j)}\in A$. There exists unique (up to a set of zero $\mu_{\Lambda}$-measure) $B\in\mathcal{B}_{\Lambda,N}$ such that for any $P_{\mathcal{D}}$ (absolutely continuous with respect to $\mu_{\mathcal{D}}$), for any simple function approximation to $\rho_{\mathcal{D}}$, and for any choice of Ansatz, we have $\tilde{P}_{\Lambda,N}(A) = \tilde{P}_{\Lambda,N}(B)$.
\end{thm}

{\em Proof:}
Fix an arbitrary $A\in\mathcal{B}_{\Lambda}$ containing at least one sample and let $\mathcal{J}$ denote the set of indices such that $\lambda^{(j)}\in A$ for $j\in\mathcal{J}$ and $B:=\cup_{j\in \mathcal{J}} \mathcal{V}_j$. By construction, $B$ is a Voronoi coverage of $A$ and by definition of counting measure $\tilde{P}_{\Lambda,N}$, for any probability measure $P_{\mathcal{D}}$ absolutely continuous with respect to $\mu_{\mathcal{D}}$ and any simple function approximation to $\rho_{\mathcal{D}}$, we have that $\tilde{P}_{\Lambda,N}(A)=\tilde{P}_{\Lambda,N}(B)$. Lemma~\ref{lemma:nonempty} gives uniqueness (see the Appendix for details). $\Box$

Theorem~\ref{thm:stoch_error} guarantees that for any event in the original $\sigma$-algebra $\mathcal{B}_{\Lambda}$ of which we want to compute its probability that there is a ``best'' approximation in the numerical $\sigma$-algebra $\mathcal{B}_{\Lambda,N}$ as long as the original event in $\mathcal{B}_{\Lambda}$ contains at least one sample in $\Lambda$. While we use probabilities to prove this, the probabilities on $\mathcal{D}$ are independent of the choices of samples on $\Lambda$ (the probability measure on outputs exists independently of our choice of approximating events in the parameter domain). In other words, Theorem~\ref{thm:stoch_error} is a statement about set-approximation not the approximation of probabilities since $\tilde{P}_{\Lambda,N}$ is exact.

%The quality of the set-approximation, or its error, is bounded by standard Monte Carlo bounds for a given choice of $A\in\mathcal{B}_{\Lambda}$, e.g.~by manipulating standard inequalities such as the Dvoretzky-Kiefer-Wolfowitz inequality for errors in empirical distribution functions. If it is the case that $A\in\mathcal{B}_{\Lambda}$ but $A\notin \mathcal{B}_{\Lambda,N}$, then it is possible to obtain computable a posteriori {\em probability} bounds for the set-approximation error (assuming the sampling density is uniform on $\Lambda$) \cite{BDW2} such as
%\begin{equation*}
% \mathbf{I} \leq \frac{\sqrt{\log(2)}}{N^{1/4}}
% % \frac{\mu_{\Lambda}(A)}{\mu_{\Lambda}(\Lambda)}
%\end{equation*}
%with probability greater than $1-2^{n+1}e^2N^n4^{-\sqrt{n}}$.

The set-approximation error can be determined {\em a priori} to any computation of the model (and bounded in $\mu_{\Lambda}$-measure independently of any probability measure) either globally or for a certain collection of events of physical importance or interest in $\mathcal{B}_{\Lambda}$. For the sake of simplicity and also for practical reasons, we consider any event of interest to belong to $\mathcal{B}_{\Lambda,N}$, i.e.~in Eq.~(\ref{eq:P_error_decompose}) we set $\tilde{\mathcal{B}}_{\Lambda} = \mathcal{B}_{\Lambda,N}$. In other words, we assume that all the events we wish to compute the probabilities of can be represented exactly (or with negligible error) in $\mathcal{B}_{\Lambda,N}$. This prevents the problem of trying to compute probabilities of complicated events that exist even in standard Borel $\sigma-$algebras, which are usually not of interest in a realistic setting where events are often defined by hypercubes or balls. Moreover, all such events $A\in\mathcal{B}_{\Lambda,N}$ have the property that $\mu_{\Lambda}(\partial A) = 0$.

Recall that $\set{I_i}_{i=1}^M$ denotes the partition on $\mathcal{D}$ used to define the simple-function approximation to the density of $P_{\mathcal{D}}$ as
\begin{equation*}
	\rho_{\mathcal{D},M} = \sum_{i=1}^M P_{\mathcal{D}}(I_i) \mathbf{1}_{I_i}(q).
\end{equation*}
Suppose that every $A_i:=Q^{-1}(I_i)\in \mathcal{B}_{\Lambda,N}$, then it follows that the error given by term $\mathbf{I}$ is \emph{zero} for all $A\in\mathcal{B}_{\Lambda,N}$. However, it is unlikely that the induced regions of generalized contours $A_i$ are exactly represented in $\mathcal{B}_{\Lambda,N}$ and this is the source of error in term $\mathbf{I}$ of Eq.~(\ref{eq:P_error_decompose}). The exact density of $P_{\Lambda}$ given the above density $\rho_{\mathcal{D},M}$ can be written as
\begin{equation}
	\rho_{\Lambda} = \sum_{i=1}^M P_{\mathcal{D}}(I_i) \mathbf{1}_{A_i}(\lambda),
\end{equation}
and the exact density of $P_{\Lambda,N}$ given $\rho_{\mathcal{D},M}$ can be written as
\begin{equation}
	\tilde{\rho}_{\Lambda,N} = \sum_{i=1}^M P_{\mathcal{D}}(I_i) \mathbf{1}_{A_{i,N}}(\lambda),
\end{equation}
where $A_{i,N}$ is the unique Voronoi coverage of $A_i$ guaranteed to exist by Theorem~\ref{thm:stoch_error} for each $i=1,2,\ldots,M$. Thus, for any $A\in\mathcal{B}_{\Lambda,N}$, we have that
\begin{equation}\label{eq:stoch_error_uncomputable}
	P_{\Lambda}(A) - P_{\Lambda,N}(A) = \sum_{i=1}^M P_{\mathcal{D}}(I_i) \left( \frac{\mu_{\Lambda}(A\cap A_i)}{\mu_{\Lambda}(A_i)} - E(A,A_i) \right).
\end{equation}
Here $E(A,A_i)$ is a computable constant determined by the Voronoi cells defining the unique covers of $A$ and $A_{i,N}$. If we use the standard MC implementation of Algorithm~\ref{Alg_MC}, then $E(A,A_i)$ is equal to the number of samples in $A\cap A_i$ over the number of samples in $A_i$. If we do not use the standard MC approximation or use non-uniform sampling so that the variant of Algorithm~\ref{Alg_MC} discussed in Section~\ref{S:non-intrusive-converge} is used, then $E(A,A_i)$ is equal to $\mu_{\Lambda}(A\cap A_{i,N})/\mu_{\Lambda}(A_{i,N})$. We also consider these options in deriving computable estimates of term $\mathbf{II}$ in Section~\ref{S:det_error} below. The error representation of Eq.~\ref{eq:stoch_error_uncomputable} is uncomputable as it requires knowledge of the exact volume of $A_i$.

To derive the computable bounds for term $\mathbf{I}$ described below, we require the following
\begin{assume}\label{A:nonempty}
Assume that for a fixed approximation to $P_{\mathcal{D}}$ that a sufficient number of samples in $\Lambda$ are taken so that there is at least one sample $\lambda^{(j)}$ in $A_{i,N}$ for each $i$ such that the associated Voronoi cell $\mathcal{V}_j$ shares no boundary with a Voronoi cell from any other region of approximate induced generalized contours $A_{k,N}$, $k\neq i$\footnote{It is not uncommon for a sufficient number of samples to be assumed in standard inequalities for empirical distribution functions, e.g.~forms of the multivariate Dvoretzky-Kiefer-Wolfowitz inequality require the square root of the number of samples to be greater than the dimension of the random variable divided by the desired $\alpha$-level \cite{Devroye77, Massart90}.}.
\end{assume}

\begin{thm}
If Assumption~\ref{A:nonempty} holds, then there exists signed computable lower and upper a posteriori bounds on term $\mathbf{I}$ of Eq.~\ref{eq:P_error_decompose} for every $A\in\mathcal{B}_{\Lambda,N}$.
\end{thm}

{\em Proof:}
Let $B_{i,N}$ denote the union of Voronoi cells belonging to $A_{i,N}$ that share no boundary with $A_{k,N}$ for all $i$, $k\neq i$. By Assumption~\ref{A:nonempty}, $B_{i,N}\neq \emptyset$ and $\mu_{\Lambda}(B_{i,N})>0$ for all $i$. Let $C_{i,N}$ denote the union of $A_{k,N}$ and all Voronoi cells sharing boundary with $A_{k,N}$. The sets $B_{i,N}$ and $C_{i,N}$ are identifiable (e.g.~using nearest neighbor searches) and computable since they belong to $\mathcal{B}_{\Lambda,N}$ for all $i$. Furthermore, we have that
\begin{equation}\label{eq:volume_error_bounds}
	\frac{\mu_{\Lambda}(A\cap B_{i,N})}{\mu_{\Lambda}(C_{i,N})} \leq \frac{\mu_{\Lambda}(A\cap A_i)}{\mu_{\Lambda}(A_i)} \leq \frac{\mu_{\Lambda}(A\cap C_{i,N})}{\mu_{\Lambda}(B_{i,N})}.
\end{equation}

Substitution of Eq.~(\ref{eq:volume_error_bounds}) into Eq.~(\ref{eq:stoch_error_uncomputable}) gives the following signed computable lower and upper a posteriori bounds on $\mathbf{I}$,
\begin{eqnarray}
	\mathbf{I} &\leq & \sum_{i=1}^M P_{\mathcal{D}}(I_i) \max\left\{ \left( \frac{\mu_{\Lambda}(A\cap B_{i,N})}{\mu_{\Lambda}(C_{i,N})} - E(A,A_i) \right), \left( \frac{\mu_{\Lambda}(A\cap C_{i,N})}{\mu_{\Lambda}(B_{i,N})} - E(A,A_i)\right) \right\}, \label{eq:stoch_upper_bd}\\
	\mathbf{I} &\geq & \sum_{i=1}^M P_{\mathcal{D}}(I_i) \min\left\{ \left( \frac{\mu_{\Lambda}(A\cap B_{i,N})}{\mu_{\Lambda}(C_{i,N})} - E(A,A_i) \right), \left( \frac{\mu_{\Lambda}(A\cap C_{i,N})}{\mu_{\Lambda}(B_{i,N})} - E(A,A_i)\right) \right\}. \label{eq:stoch_lower_bd}
\end{eqnarray}
This completes the proof $\Box$.

\subsection{Deterministic error and misidentification of inverse sets}\label{S:det_error}

As described above in Section~\ref{S:stoch_error}, we assume $\tilde{\mathcal{B}}_{\Lambda}=\mathcal{B}_{\Lambda,N}$ in Eq.~(\ref{eq:P_error_decompose}). By assumption of local differentiability of map $Q$, we assume the numerical method used to compute $Q_h$ produces a piecewise continuous approximation to $Q$ and there exists a piecewise continuous error function $e$ such that
\[
	Q(\lambda)=Q_h(\lambda)+e(\lambda), \ \forall \lambda\in\Lambda.
\]
We assume there exists computable { a posteriori} error estimates for $e$ denoted by $e_h$ such that
\[
 Q(\lambda) \approx Q_h(\lambda)+e_h(\lambda), \ \forall \lambda\in\Lambda.
\]
Analogous to Algorithm~\ref{Alg_old}, let $A_{i,N}$ and $A_{i,N,h}$ denote the Voronoi coverages of $Q^{-1}(I_i)$ and $Q^{-1}_h(I_i)$, respectively. Following Algorithm~\ref{Alg_MC}, we have that, for any $A\in\mathcal{B}_{\Lambda,N}$, term $\mathbf{II}$ in Eq.~(\ref{eq:P_error_decompose}) can be written as
\begin{equation}\label{eq:det_error}
	\mathbf{II} = \sum_{i=1}^M P_{\mathcal{D}}(I_i)\left( \frac{ \set{\# j\, : \, \lambda^{(j)}\in A_{i,N,h}\cap A} }{ \set{\# j\, : \, \lambda^{(j)} \in A_{i,N,h}} } - \frac{ \set{\# j\, : \, \lambda^{(j)}\in A_{i,N}\cap A} }{ \set{\# j\, : \, \lambda^{(j)} \in A_{i,N}} }   \right).
\end{equation}
Let $J_{i,A}$, and $J_{i}$, denote $\set{\# j\, : \, \lambda^{(j)}\in A_{i,N,h}\cap A}$ and $\set{\# j\, : \, \lambda^{(j)} \in A_{i,N,h}}$, respectively. Note that the numbers $J_{i,A}$ and $J_{i}$ are computable. Terms involving $A_{i,N}$ are not computable since they require the exact map $Q$. We may invert $Q_h+e_h$ to more accurately estimate $A_{i,N}$, and define the computable approximations
\begin{eqnarray}
	\set{\# j\, : \, \lambda^{(j)}\in A_{i,N}\cap A} &\approx & \set{\# j\, : \, \lambda^{(j)}\in A, \text{and } Q_h(\lambda^{(j)})+e_h(\lambda^{(j)}) \in I_i}, \label{eq:approx1} \\
     \set{\# j\, : \, \lambda^{(j)} \in A_{i,N}}   &\approx &  \set{\# j\, : \, Q_h(\lambda^{(j)})+e_h(\lambda^{(j)}) \in I_i}. \label{eq:approx2}
\end{eqnarray}
Let $J_{i,A,e}$ and $J_{i,e}$ denote the computable approximations of Eqs.~(\ref{eq:approx1}) and (\ref{eq:approx2}), respectively. Substituting these approximations into Eq.~(\ref{eq:det_error}) and re-arranging terms gives the computable a posteriori estimate
\begin{equation}
	\mathbf{II} \approx  \sum_{i=1}^M P_{\mathcal{D}}(I_i) \left( \frac{ J_{i,A}J_{i,e} - J_{i,A,e}J_{i} } { J_i J_{i,e} }	\right)
\end{equation}

As discussed in Section~\ref{S:non-intrusive-converge}, if we replace the counting vector in Algorithm~\ref{Alg_MC} with estimates of volumes for the Voronoi cells, then the above summations of number of cells are replaced by the summations of the volumes of the associated Voronoi cells. For example, in this case we set $J_{i,A}$ to be the number determined by the following
\[
	\sum_ {\set{j\, : \, \lambda^{(j)}\in A_{i,N,h}\cap A}} \mu_{\Lambda}(\mathcal{V}_j).
\]

\subsection{Improving computations with error estimates}\label{S:Error_Improve}

Computable error estimates and/or error bounds are useful in adaptive algorithms. This is beyond the scope of this paper, but is the subject of future work on adaptive sampling strategies for computing accurate probabilities of specified events, e.g.~rare events.

 In the numerics below, we consider an alternative use of the a posteriori error estimates of the map $Q_h$. Specifically, we use the a posteriori error estimates to improve the functional evaluations. This is motivated by previous work on improved linear functionals where a posteriori error estimates are used to improve the pointwise accuracy of the QoI map \cite{Butler_2012b} and improve the pointwise accuracy of probability distributions  propagated through likelihoods involving a QoI map \cite{Butler2013b}.

To show the effect of using a posteriori error estimates to improve the accuracy of computed probabilities of a counting measure we first compute a so-called ``reference solution'' that is an accurate approximation to $P_{\Lambda}$ on a fine grid of uniformly spaced rectangles partitioning $\Lambda$. When computing the reference solution, we use a posteriori error estimates on the fine grid to correct for the deterministic error of the map. We follow this by the random sampling algorithm for computing counting measures both with and without the a posteriori error estimate correcting for the deterministic error. Finally, we demonstrate that the errors of the counting measures with respect to the reference solution evaluated on generating events for a given $\tilde{\mathcal{B}}_{\Lambda}$ are reduced globally when the a posteriori error estimate is used to improve the pointwise accuracy of $Q_h$.

\section{Numerical Results}\label{S:Numerics}

\subsection{Improving probabilities with a posteriori error estimates}

\begin{figure}[htb]
	\centering
\includegraphics[height=4.1cm,clip=true]{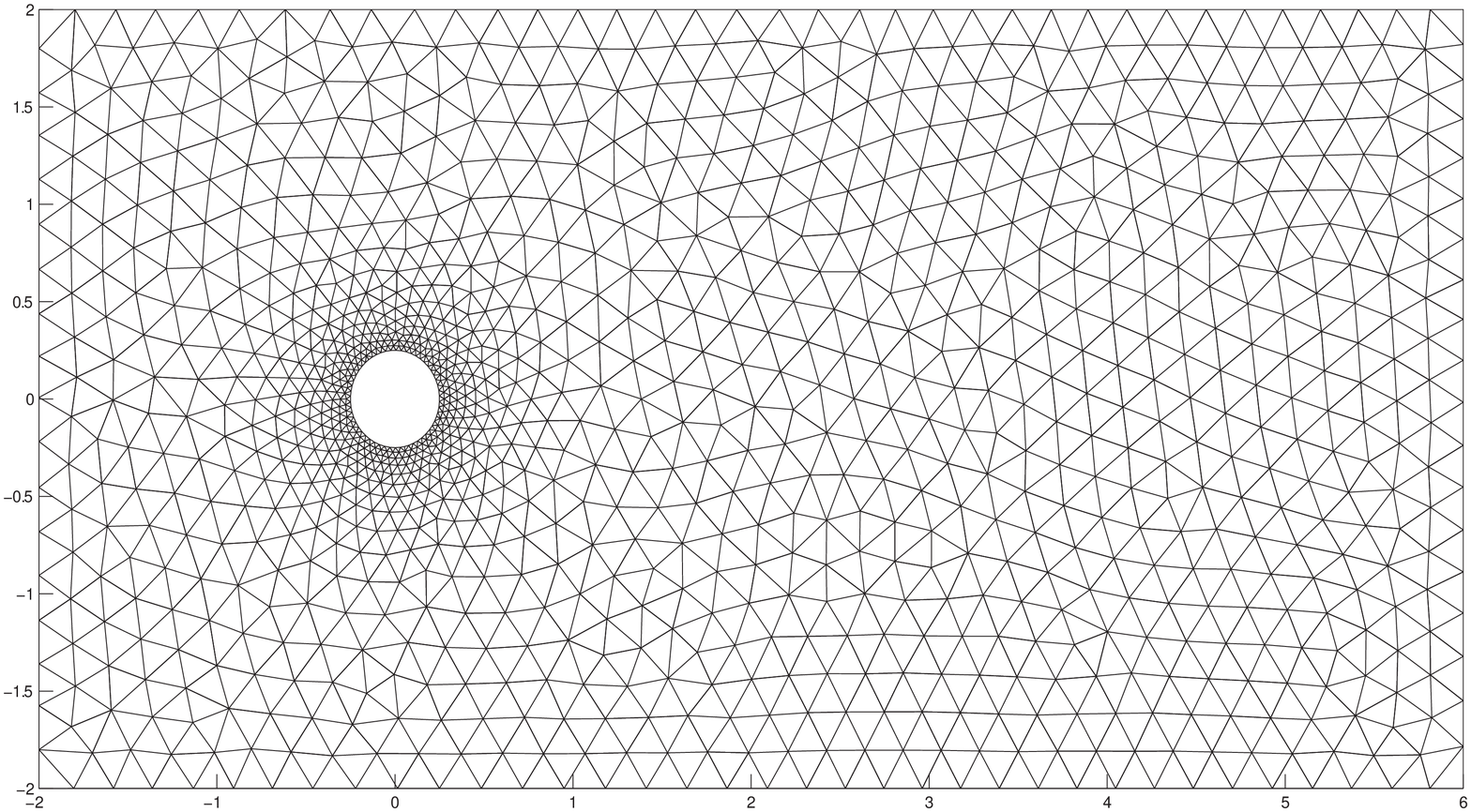} \includegraphics[height=4.3cm,clip=true]{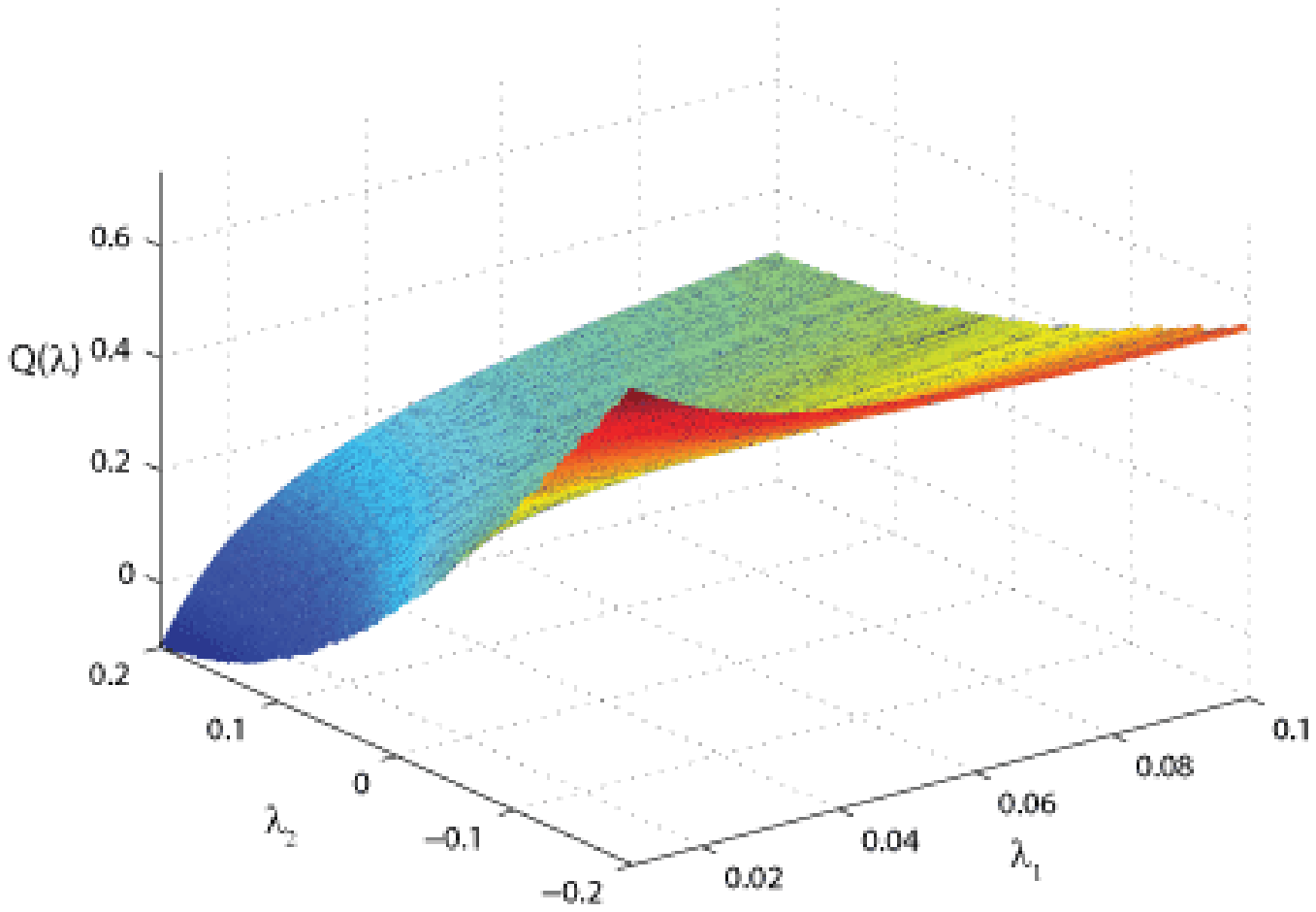}
		\caption{Left: Physical domain $\Omega$ in Eq.~\ref{eq:NS}. Right: The QoI map $Q(\lambda)$ is the first component of velocity at the point $(1,0.2)$ in $\Omega$. The viscosity $\nu$ is taken as $\lambda_1$ in the parameter domain. The vertical displacement of the center of the cylindrical hole is taken as $\lambda_2$ in the parameter domain.}
	\label{fig:NS_domain}
\end{figure}

We take as the model the Navier-Stokes equation
\begin{equation}\label{eq:NS}
	\frac{\partial \mathbf{u}}{\partial t} - \nu \nabla \cdot\nabla \mathbf{u} + \mathbf{u}\cdot\nabla \mathbf{u} + \nabla p = \mathbf{g}_{\mathbf{u}}, \ x\in \Omega, \ t\in (t_0,t_f),
\end{equation}
where $\Omega$ is defined on a rectangular domain with a cylindrical hole as shown in the left plot of Figure~\ref{fig:NS_domain}. We describe the boundary conditions as specified for a similar problem in \cite{Wildey_2014} with fixed parameter values. For the rectangular boundary, the inflow velocity is set to $(1,0)$, non-penetrating boundary conditions ($\mathbf{u}\cdot\mathbf{n}=0$) are set along the top and bottom boundaries, and a natural outflow condition is set on the right boundary. We use no-slip boundary conditions ($\mathbf{u}=0$) on the cylindrical boundary. The diameter of the cylinder is $1/2$ and has uncertain vertical displacement. For this choice of cylinder and possible locations, the Reynolds number is determined from the viscosity as Re$=(2\nu)^{-1}$. We focus on steady-state solutions, so we consider the situation where the fluid viscosity is uncertain but known to be within $[0.01,0.1]$ resulting in Reynolds numbers bounded by $50$. Thus, the parameters for this problem are the viscosity and the vertical displacement of the center of the cylinder. Let $\lambda\in\Lambda = [0.01,0.1]\times[-0.2,0.2]\subset\mathbb{R}^2$ denote the parameter domain with viscosity the first component.

For the numerical solution of the model, we use SUPG-PSPG-LSIC stabilization with piecewise linear finite elements for the forward problem. The QoI we consider is the first component of the velocity at the point $(1,0.2)$ which is just behind the cylinder and off the center-line. This QoI defines a differentiable map shown on the right plot of Figure~\ref{fig:NS_domain}. To estimate the error in the numerically evaluated QoI, we use an a posteriori error estimate using solution to an adjoint problem. We solve the adjoint problem using SUPG-PSPG-LSIC stabilization with piecewise quadratic finite elements. This approach has been demonstrated to produce reliably accurate error estimates for a variety of QoI in \cite{Wildey_2014}.

\begin{figure}[htb]
	\centering
\includegraphics[height=4.5cm,bb = 80 205 534 578 clip=true]{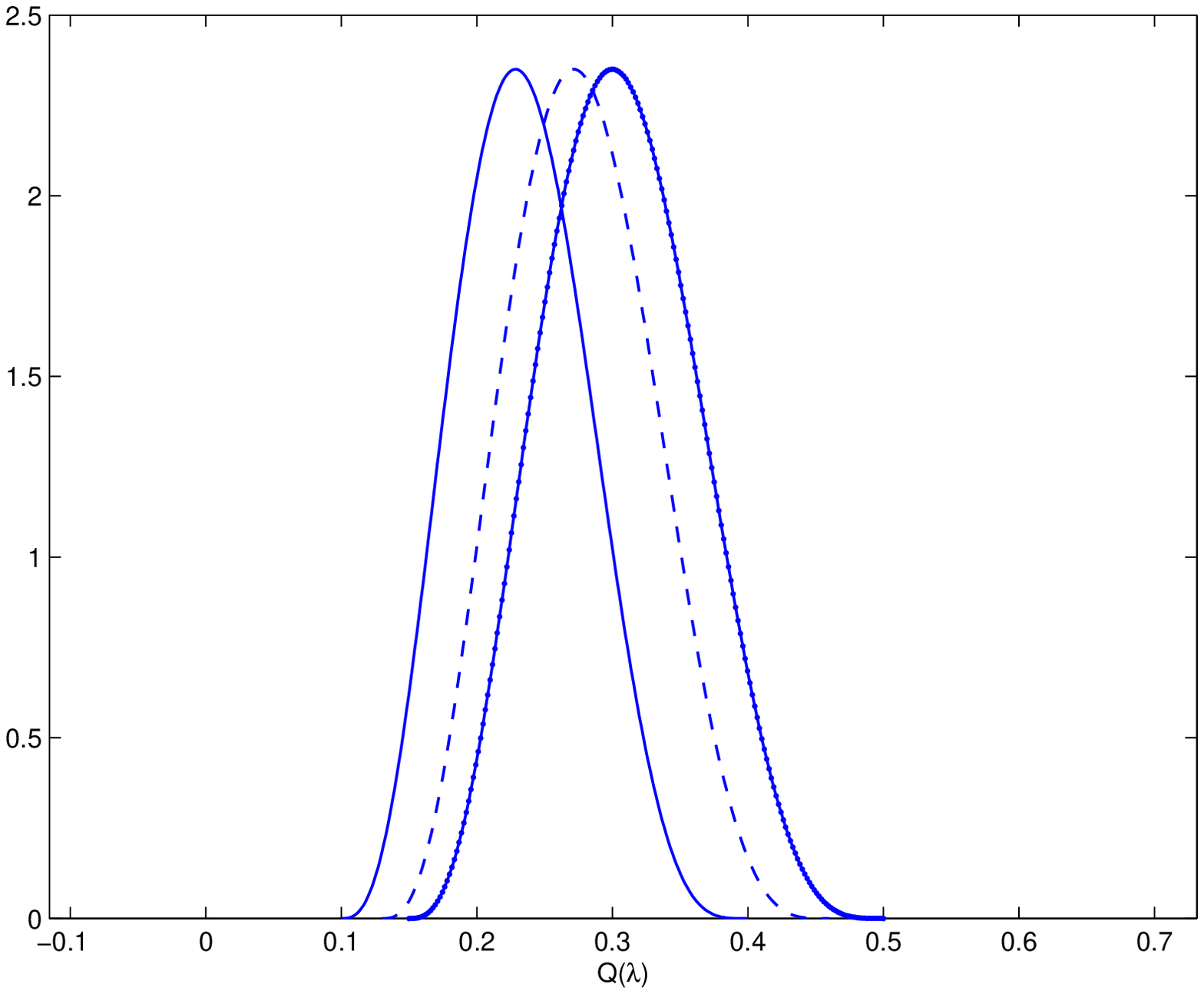} \includegraphics[height=4.5cm,,bb= 165 205 429 593,clip=true]{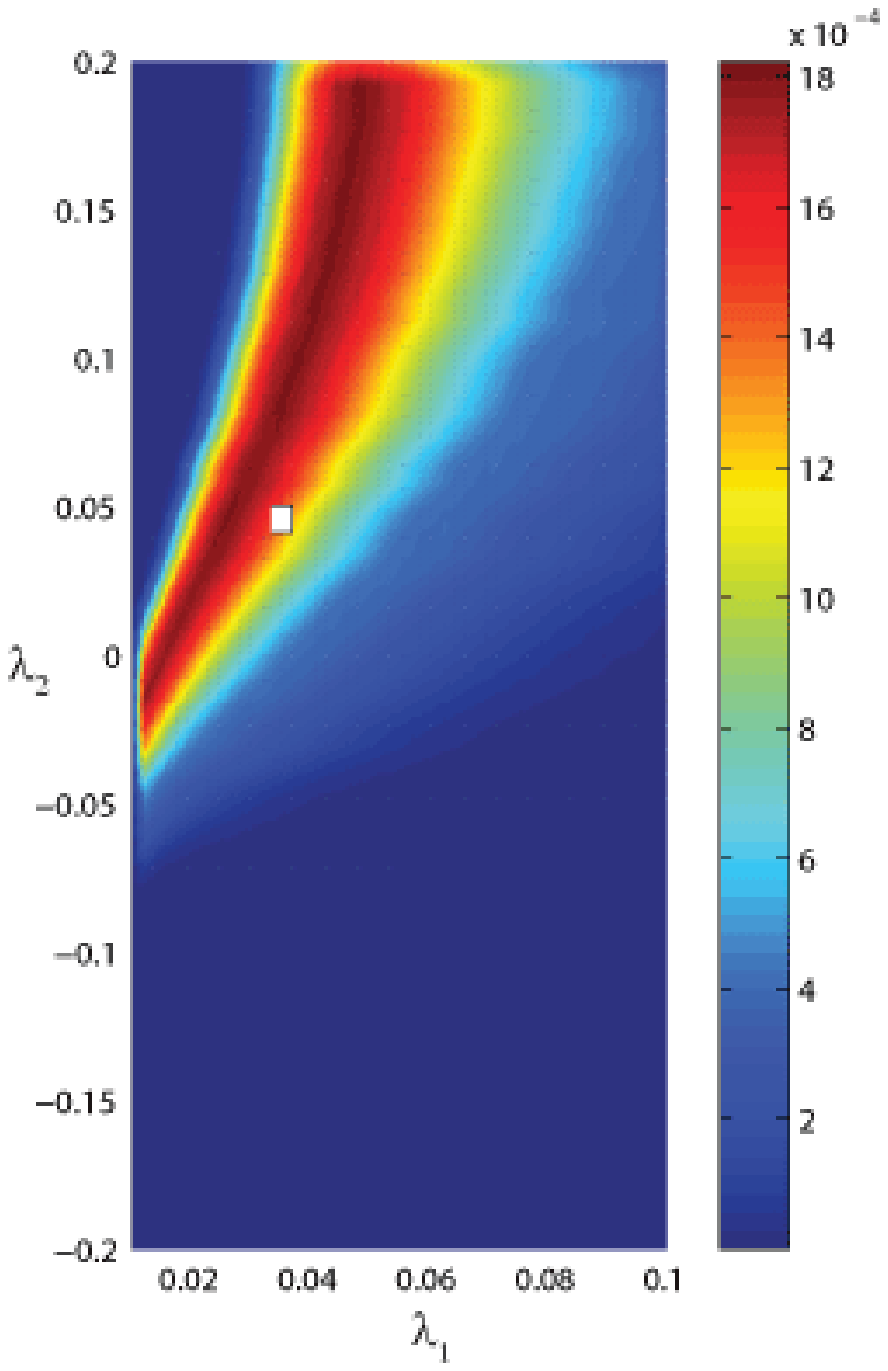}
\includegraphics[height=4.5cm,bb= 165 205 429 593,clip=true]{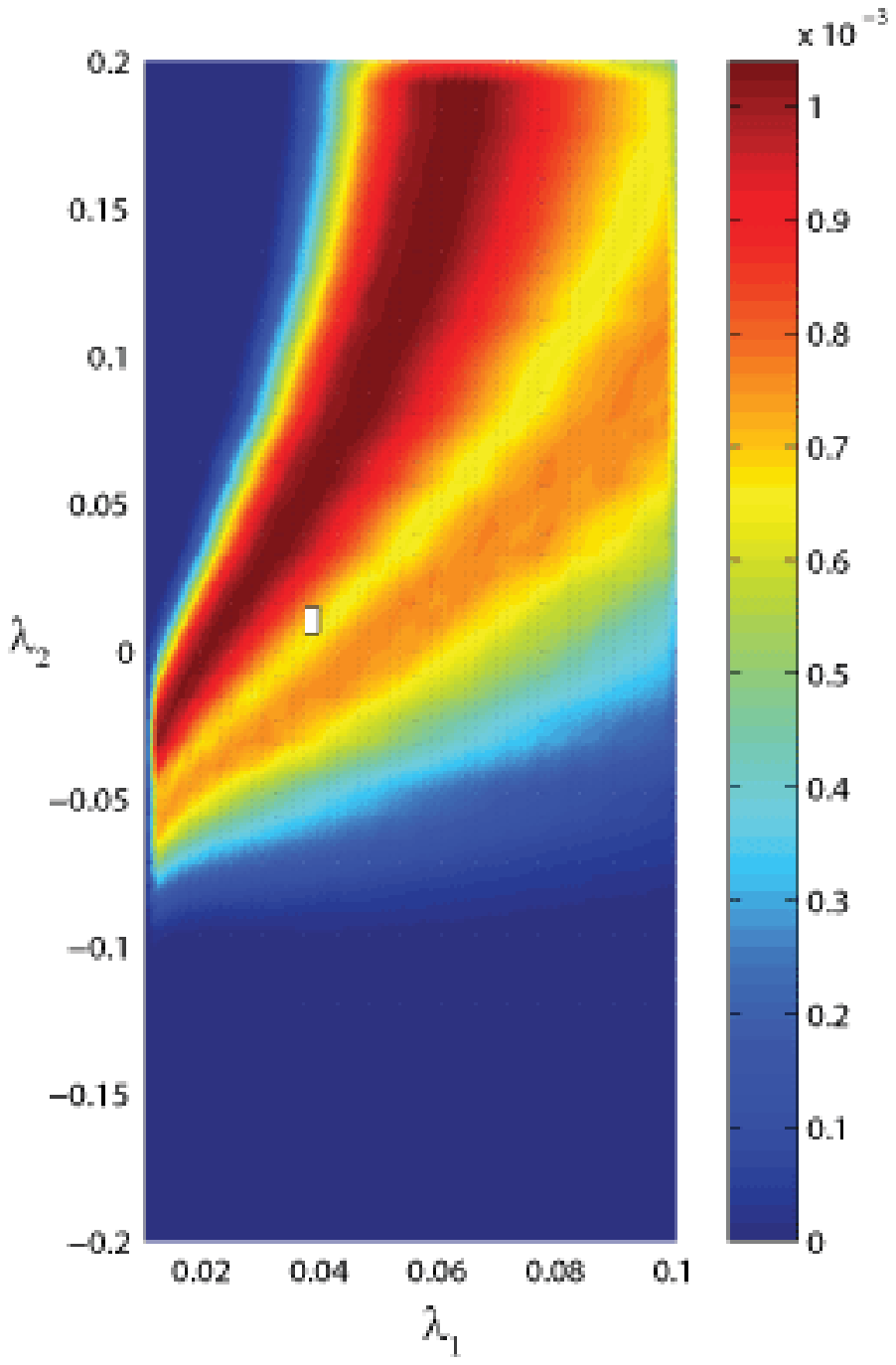}
\includegraphics[height=4.5cm,bb= 165 205 429 593,clip=true]{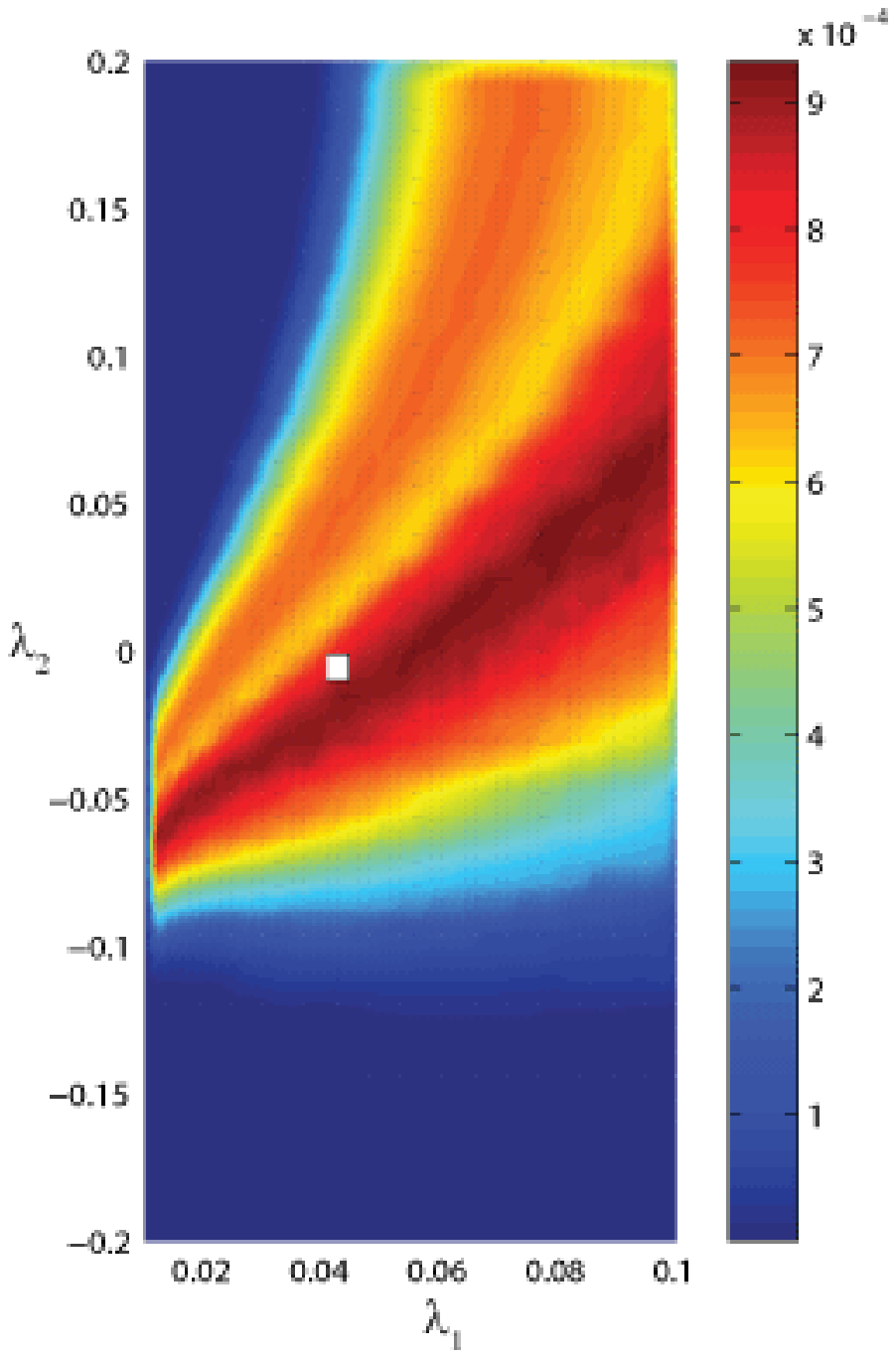}
		\caption{Left Plot: Plots of the Beta$(4,5)$ distributions centered at $Q(0.0346,0.0458)$ (solid line with leftmost peak), $Q(0.0371,0.0104)$ (dashed line with middle peak), and $Q(0.0422,-0.0056)$ (dotted line with rightmost peak). Remaining plots show approximations of $P_{\Lambda}(R_j)$ on $50\times 50$ rectangles $\set{R_j}$ partitioning $\Lambda$ so that moving left-to-right the associated output densities have means at $Q(0.0346,0.0458)$, $Q(0.0371,0.0104)$, and $Q(0.0422,-0.0056)$, respectively. The parameter values $(0.0346,0.0458)$, $(0.0371,0.0104)$, and $(0.0422,-0.0056)$ are denoted by white boxes in these plots.}
	\label{fig:NS_various_distributions}
\end{figure}

\begin{figure}[htb]
	\centering
 \includegraphics[height=5cm,,bb= 135 200 492 580,clip=true]{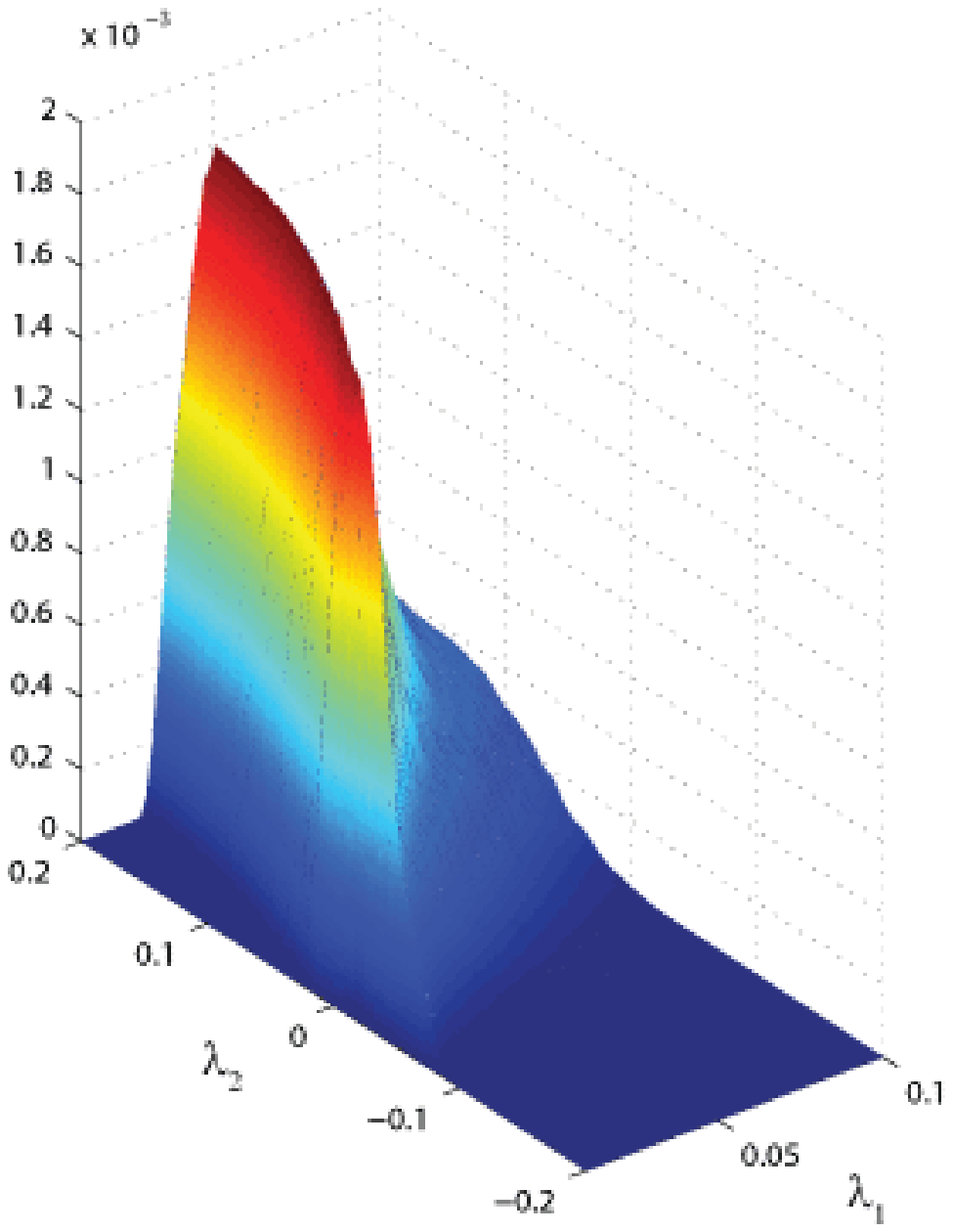}
\includegraphics[height=5cm,bb= 135 200 492 580,clip=true]{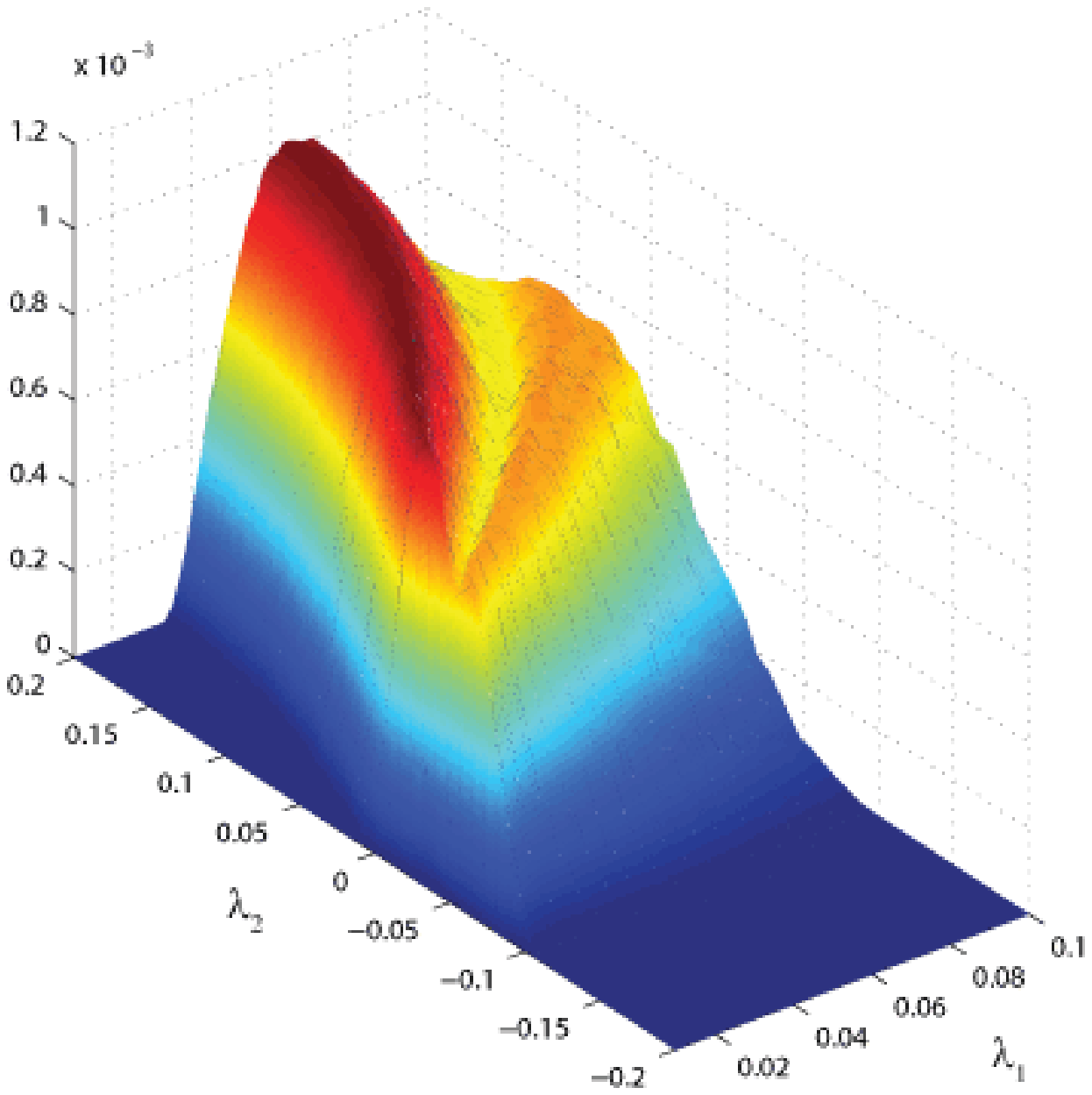}
\includegraphics[height=5cm,bb= 135 200 492 580,clip=true]{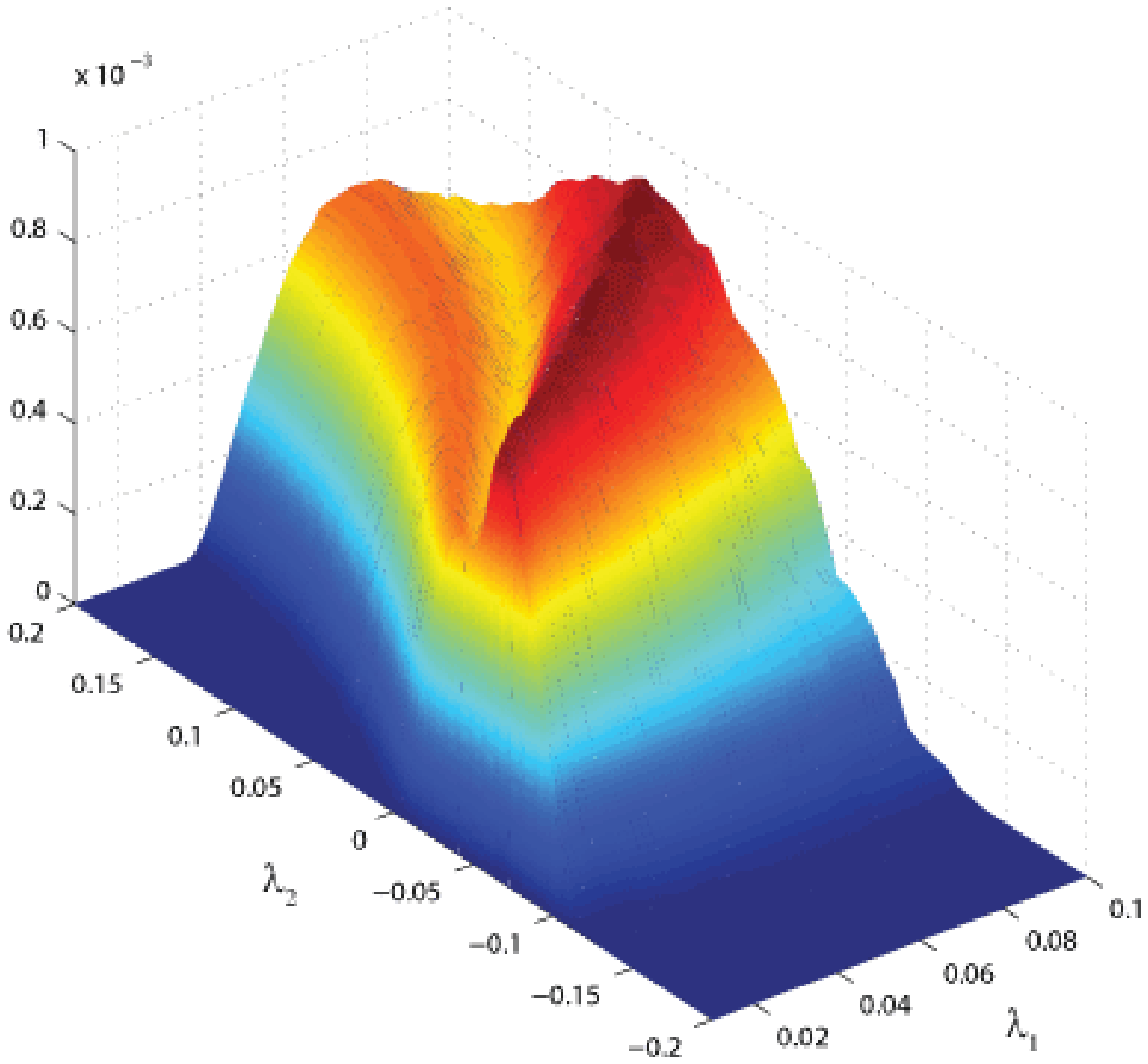}
		\caption{Left: to Right: A 3-D view of approximations of $P_{\Lambda}(R_j)$ on $50\times 50$ rectangles $\set{R_j}$ partitioning $\Lambda$ so that moving left-to-right the associated output Beta densities have means at $Q(0.0346,0.0458)$, $Q(0.0371,0.0104)$, and $Q(0.0422,-0.0056)$, respectively. The parameter values $(0.0346,0.0458)$, $(0.0371,0.0104)$, and $(0.0422,-0.0056)$ are denoted by white boxes in these plots.}
	\label{fig:NS_various_distributions_view3}
\end{figure}

We create an approximate reference solution to $P_{\Lambda}$ in the following way. First, we use 62,500 samples taken from a regularly spaced $250\times250$ set of points in $\Lambda$. For each sampled parameter, we compute both the QoI and the associated a posteriori error estimate. Since $\mathcal{D}$ is compact, we assume the output density is approximated parametrically by a Beta distribution Beta$(\alpha,\beta)$ within the range of the QoI map. The parameters $\alpha$ and $\beta$ are chosen so the mean of the distribution is at a specified ``exact'' observed value and the distribution ``appears Normal'' in shape\footnote{We could also simply truncate Normal distributions and renormalize. We simply avoid this nuisance here.}. For the cases shown here, we set $\alpha=4$ and $\beta=5$ and took the mean QoI values from a more accurate (i.e.~a posteriori error corrected) mapping of known parameter values, see the leftmost plot of Figure~\ref{fig:NS_various_distributions}. We used three separate parameters within $\Lambda$ that map to nearby QoI values to demonstrate the sensitivity of inverse probability measures to small changes in the location of high probability intervals in $\mathcal{D}$. We approximate each $\rho_{D}$ with a simple function $\rho_{D,M}$ defined on a grid of $200$ uniformly spaced points defining $199$ bins partitioning $\mathcal{D}$ and binning $1E+6$ i.i.d.~samples from the associated Beta distributions. We then follow the steps of Algorithm~\ref{Alg_old} using the 62,500 samples of $\Lambda$ evaluated on the error corrected map. In Figures~\ref{fig:NS_various_distributions} and \ref{fig:NS_various_distributions_view3}, we show the resulting approximate plots of $P_{\Lambda}(R_j)$ computed on $50\times 50$ uniformly sized rectangles $\set{R_j}$ partitioning $\Lambda$. In other words, we use the $250\times 250$ partition to define $P_{\Lambda,62500}\approx P_{\Lambda}$ and create plots of probabilities on a $50\times 50$ grid of rectangles, which gives a good visual approximation to the shape of the density function. Below, we generally refer to any computation involving $P_{\Lambda,62500}$ as a reference solution.

We observe in the plots of the reference solution in Figures~\ref{fig:NS_various_distributions} and \ref{fig:NS_various_distributions_view3} that in two cases the associated densities will exhibit a type of ``bi-modality'' across two distinct regions of induced generalized contours that appear to fan out as viscosity (i.e.~the $\lambda_1$ parameter) is increased. This is due in part to the fact that the peaks of the associated output Beta distributions with means at $Q(0.0371,0.0104)$ and $Q(0.0422,-0.0056)$ are near a QoI value for which the surface shown in Figure~\ref{fig:NS_domain} has a saddle point. We focus on the use of a counting measure to compute approximation of probabilities in $\Lambda$ associated with the output Beta$(4,5)$ distribution with mean at $Q(0.0422,-0.0056)$ where this bi-modality feature appears particularly predominant in the reference solution.

% We note that the reference solution for the associated output Beta distribution with mean at $Q(0.0371,0.0104)$ has more probability in events of $\Lambda$ containing higher values of $\lambda_2$ and smaller values of $\lambda_1$ compared to the reference solution associated with output Beta distribution with mean value of $Q(0.0422,-0.0056)$.

\begin{figure}[htb]
	\centering
 \includegraphics[height=5cm,,bb= 165 205 429 593,clip=true]{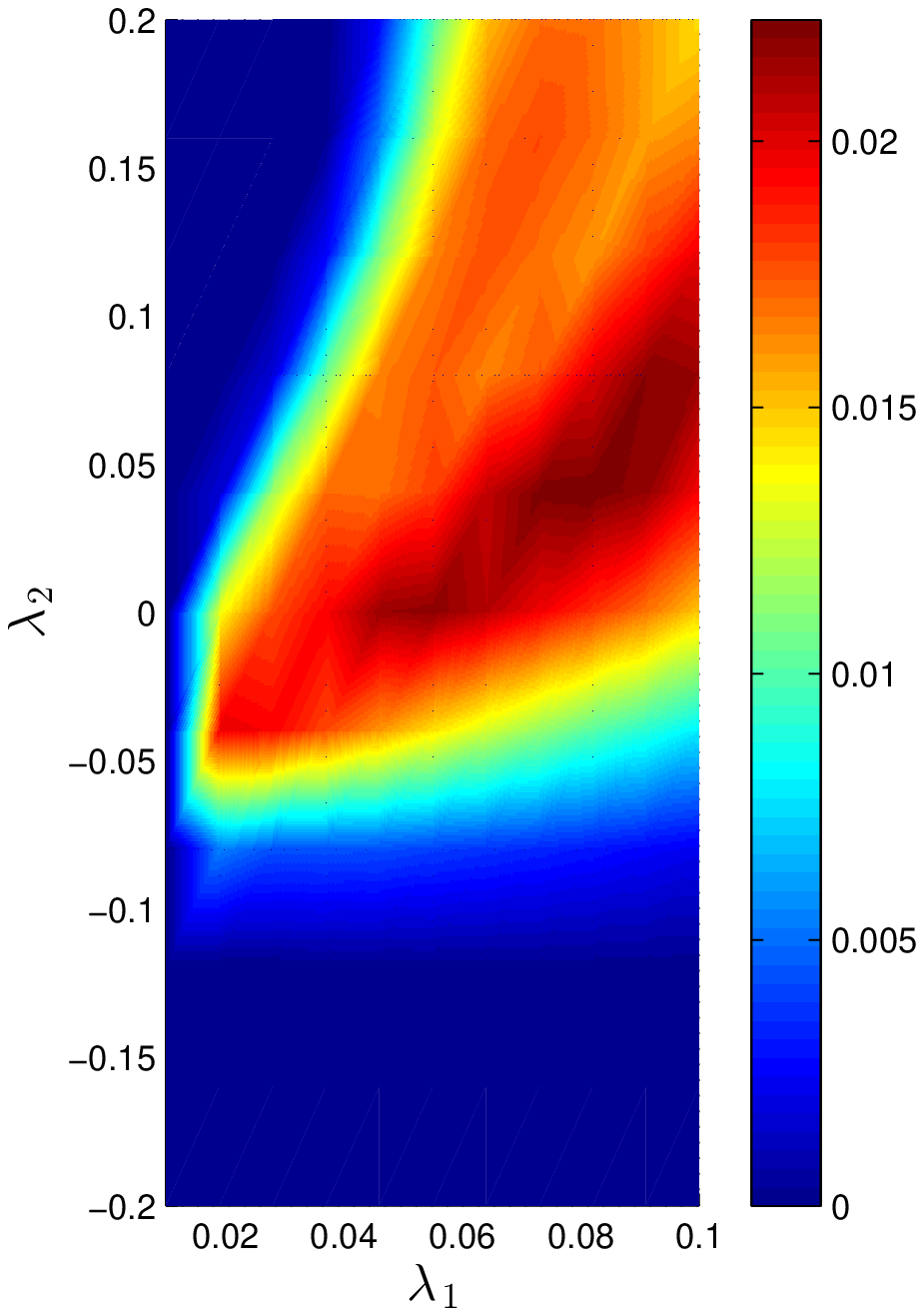}
\includegraphics[height=5cm,bb= 165 205 429 593,clip=true]{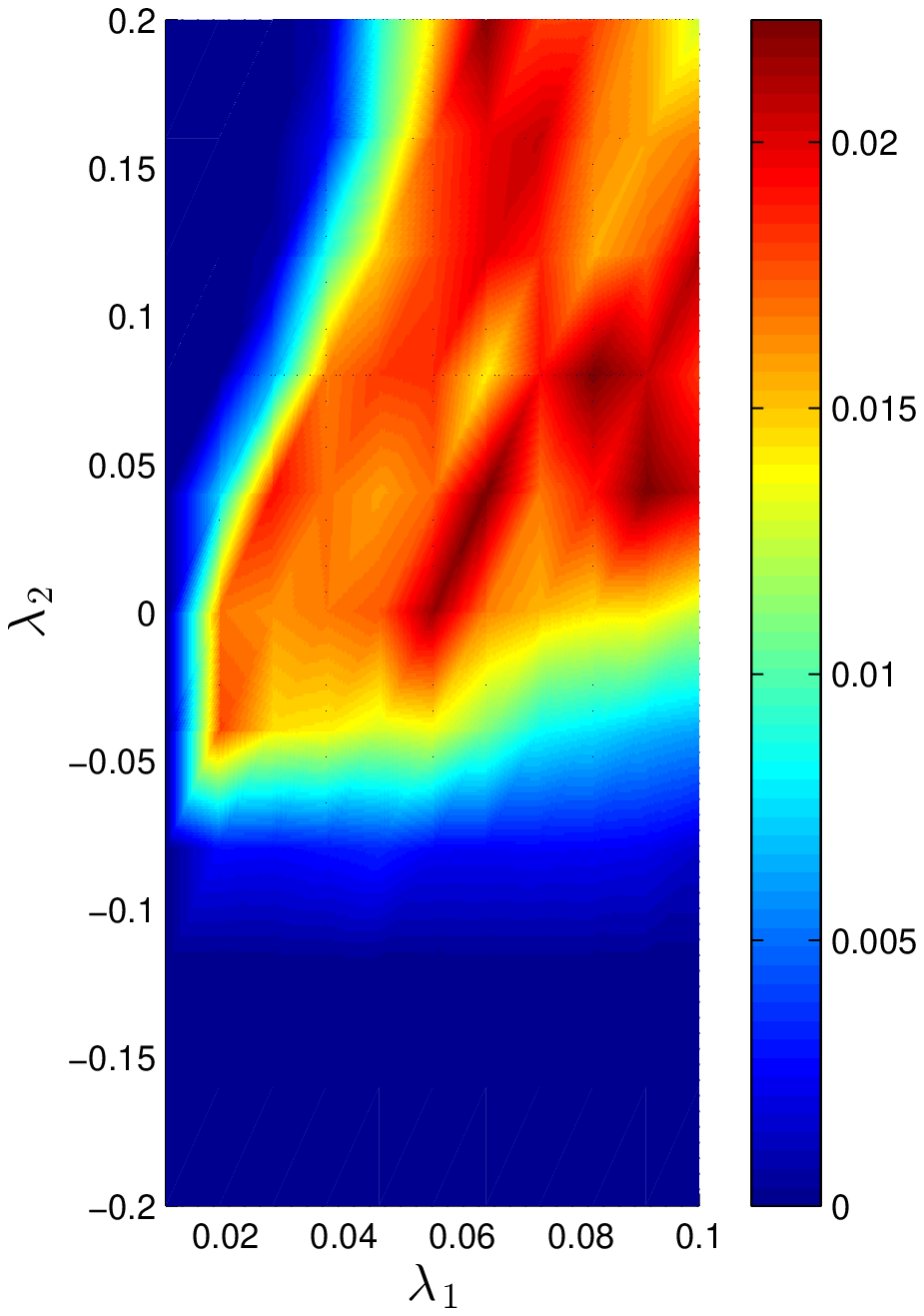}
\includegraphics[height=5cm,bb= 165 205 429 593,clip=true]{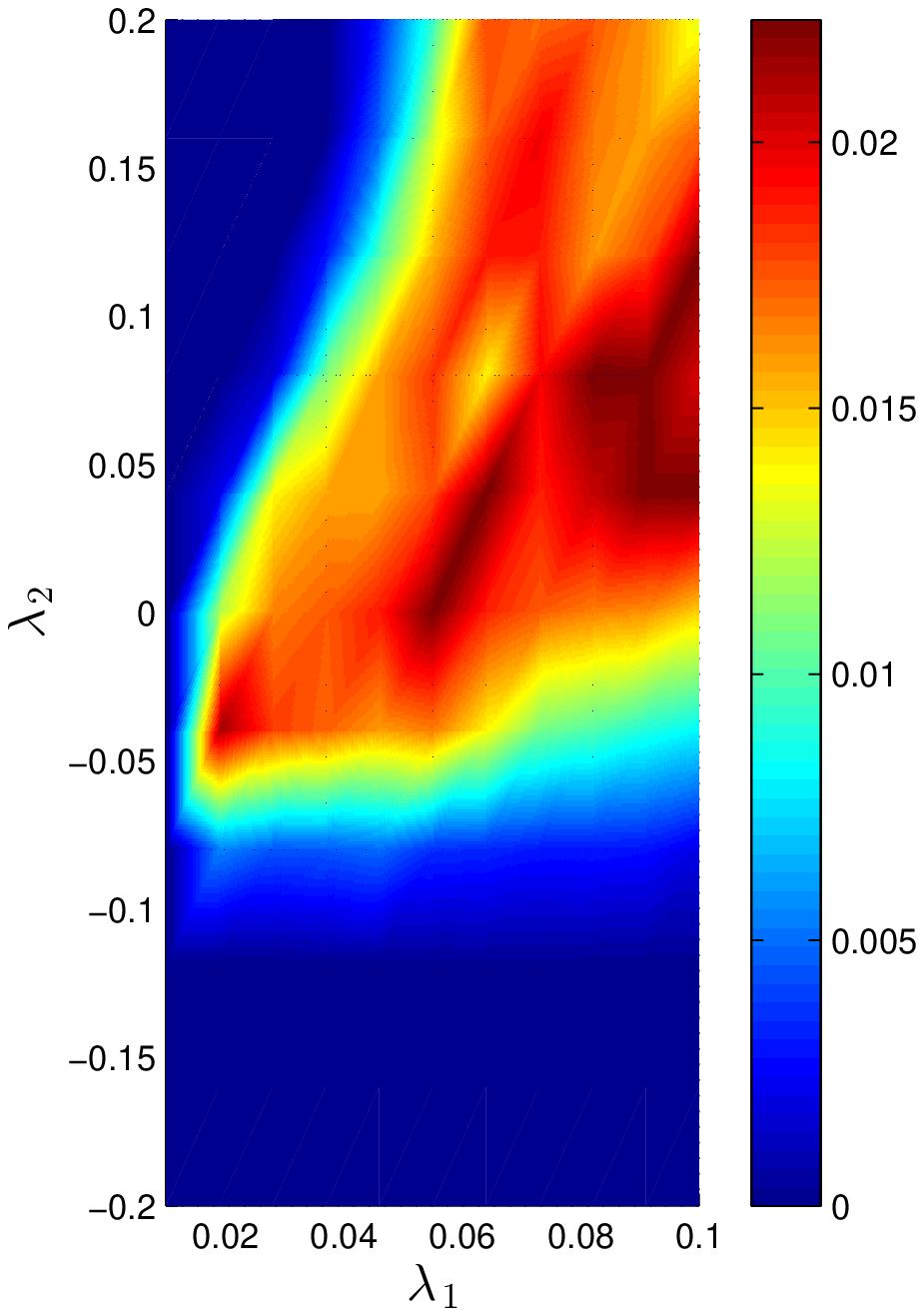}
		\caption{Reference solution (left), $\tilde{P}_{\Lambda,15k,h_1}$ (middle), and  $\tilde{P}_{\Lambda,15k,h_2}$ (right) evaluated on the 100 generating events of $\tilde{\mathcal{B}}_{\Lambda}$}
	\label{fig:NS_LowRes}
\end{figure}

\begin{figure}[htb]
	\centering
\includegraphics[height=5cm,bb= 135 200 492 580,clip=true]{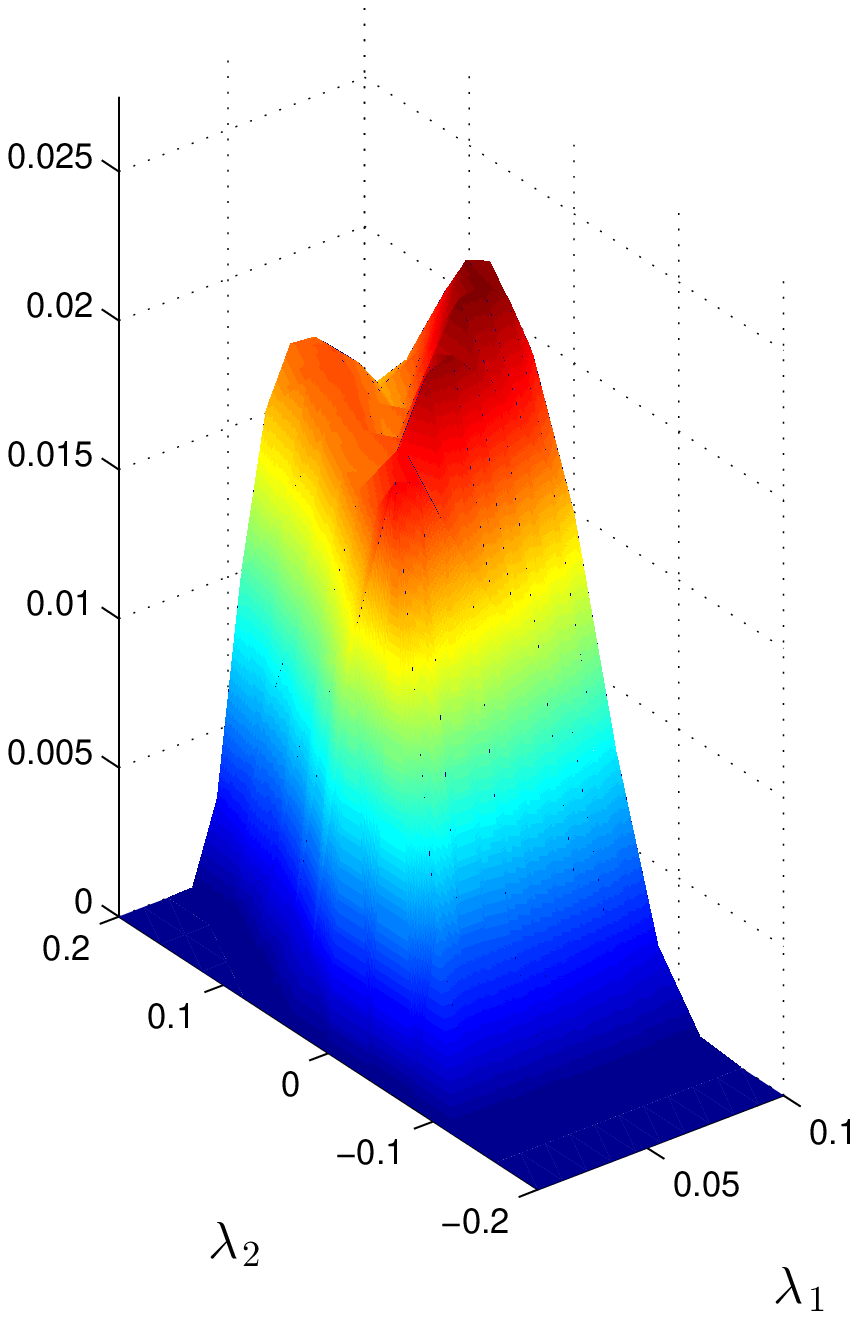}
\includegraphics[height=5cm,bb= 135 200 492 580,clip=true]{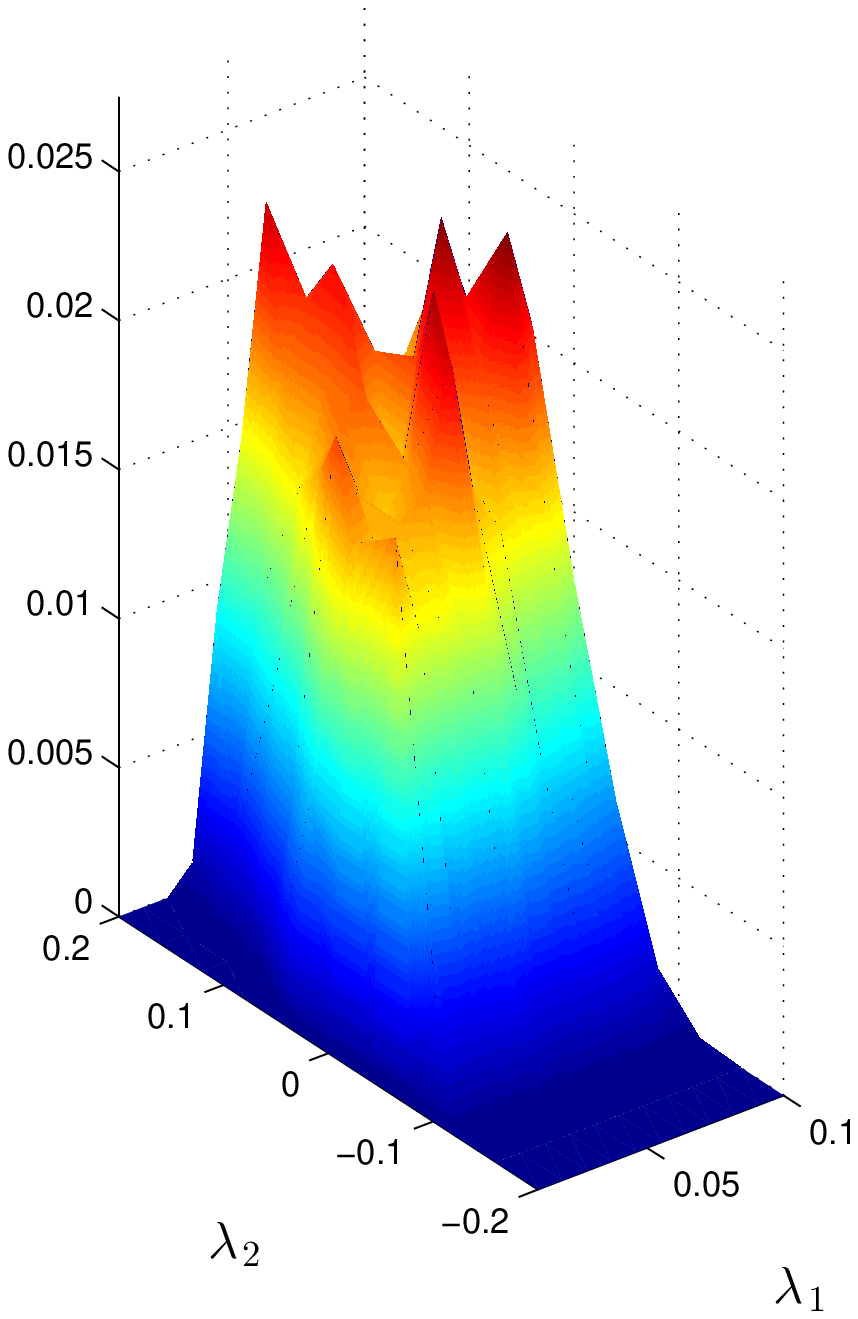}
\includegraphics[height=5cm,bb= 135 200 492 580,clip=true]{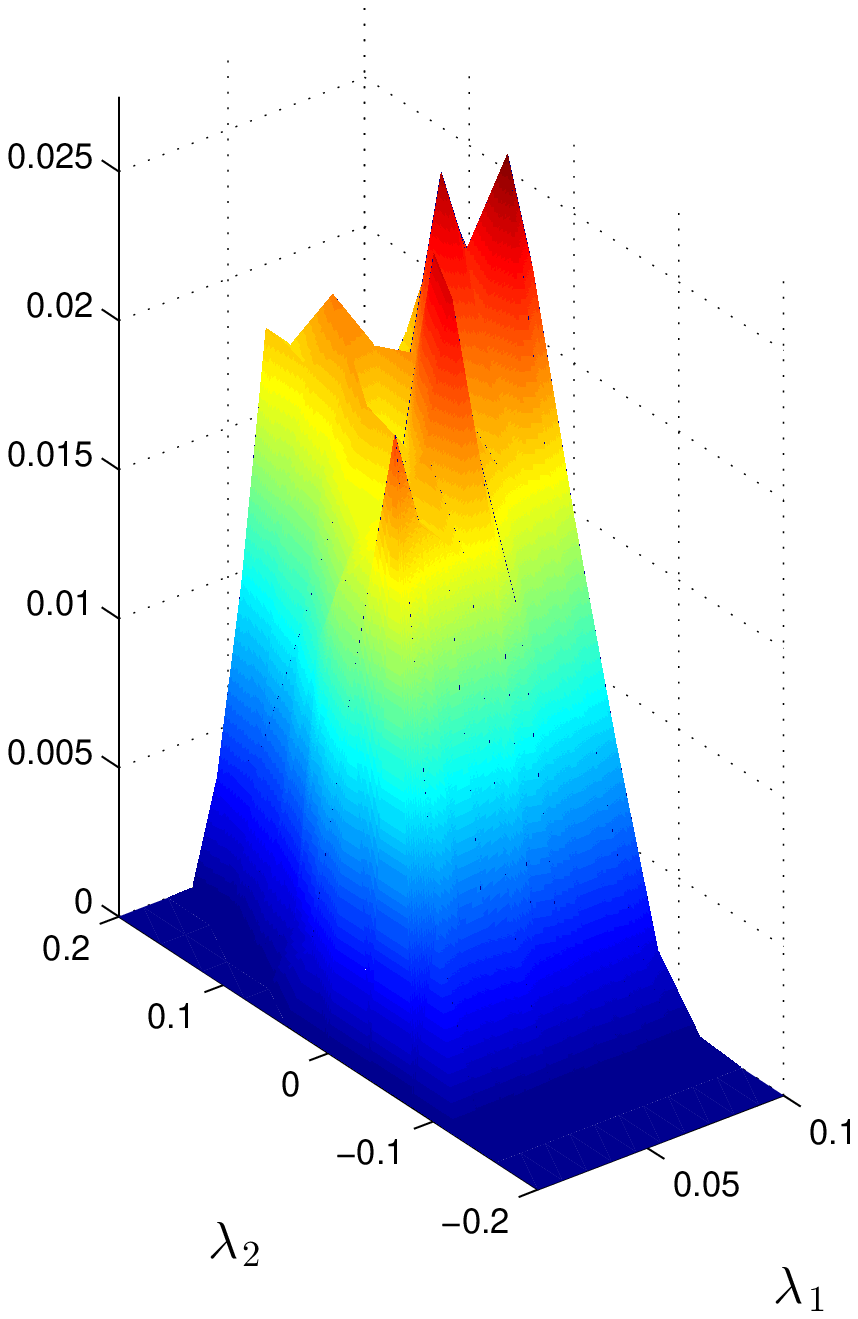}
		\caption{3-D views of reference solution (left), $\tilde{P}_{\Lambda,15k,h_1}$ (middle), and  $\tilde{P}_{\Lambda,15k,h_2}$ (right) evaluated on the 100 generating events of $\tilde{\mathcal{B}}_{\Lambda}$}
	\label{fig:NS_LowRes_view3}
\end{figure}

We define $\tilde{\mathcal{B}}_{\Lambda}$ as the $\sigma-$algebra generated from a $10\times 10$ uniform rectangular grid of $\Lambda$, i.e.~we assume any event we want to measure in probability can be defined by unions of a subset of 100 uniformly sized rectangles in $\Lambda$. We plot the reference solution computed on the generating events of $\tilde{\mathcal{B}}_{\Lambda}$ in the left plots of Figures~\ref{fig:NS_LowRes} and \ref{fig:NS_LowRes_view3}.

We use 15,000 i.i.d.~uniform random samples in $\Lambda$ to compute two separate counting measures denoted by $\tilde{P}_{\Lambda,15k,h_1}$ and $\tilde{P}_{\Lambda,15k,h_2}$. We compute $\tilde{P}_{\Lambda,15k,h_1}$ using $Q_h$ to evaluate the samples in $\Lambda$ and compute $\tilde{P}_{\Lambda,15k,h_2}$ using $Q_h+e_h$ (i.e.~using computed a posteriori error estimates to correct for error in the numerically evaluated QoI) to evaluate the same set of samples in $\Lambda$. We plot  $\tilde{P}_{\Lambda,15k,h_1}$ and  $\tilde{P}_{\Lambda,15k,h_2}$ given the generating events of $\tilde{\mathcal{B}}_{\Lambda}$ in the middle and right plots, respectively, of Figures~\ref{fig:NS_LowRes} and \ref{fig:NS_LowRes_view3}.

\begin{figure}[htb]
	\centering
\includegraphics[height=5cm,bb= 165 205 429 593,clip=true]{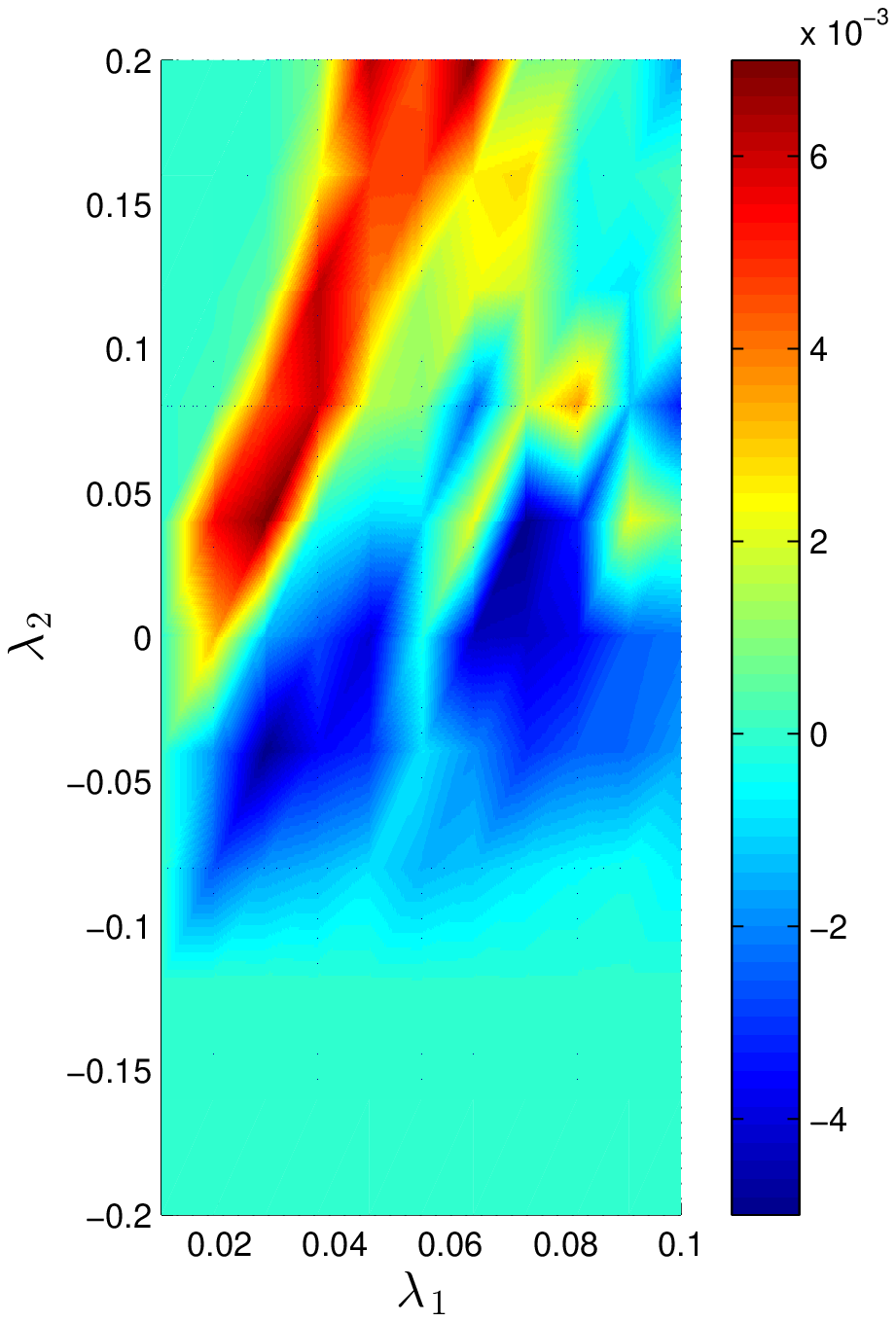}
\includegraphics[height=5cm,bb= 155 200 475 580,clip=true]{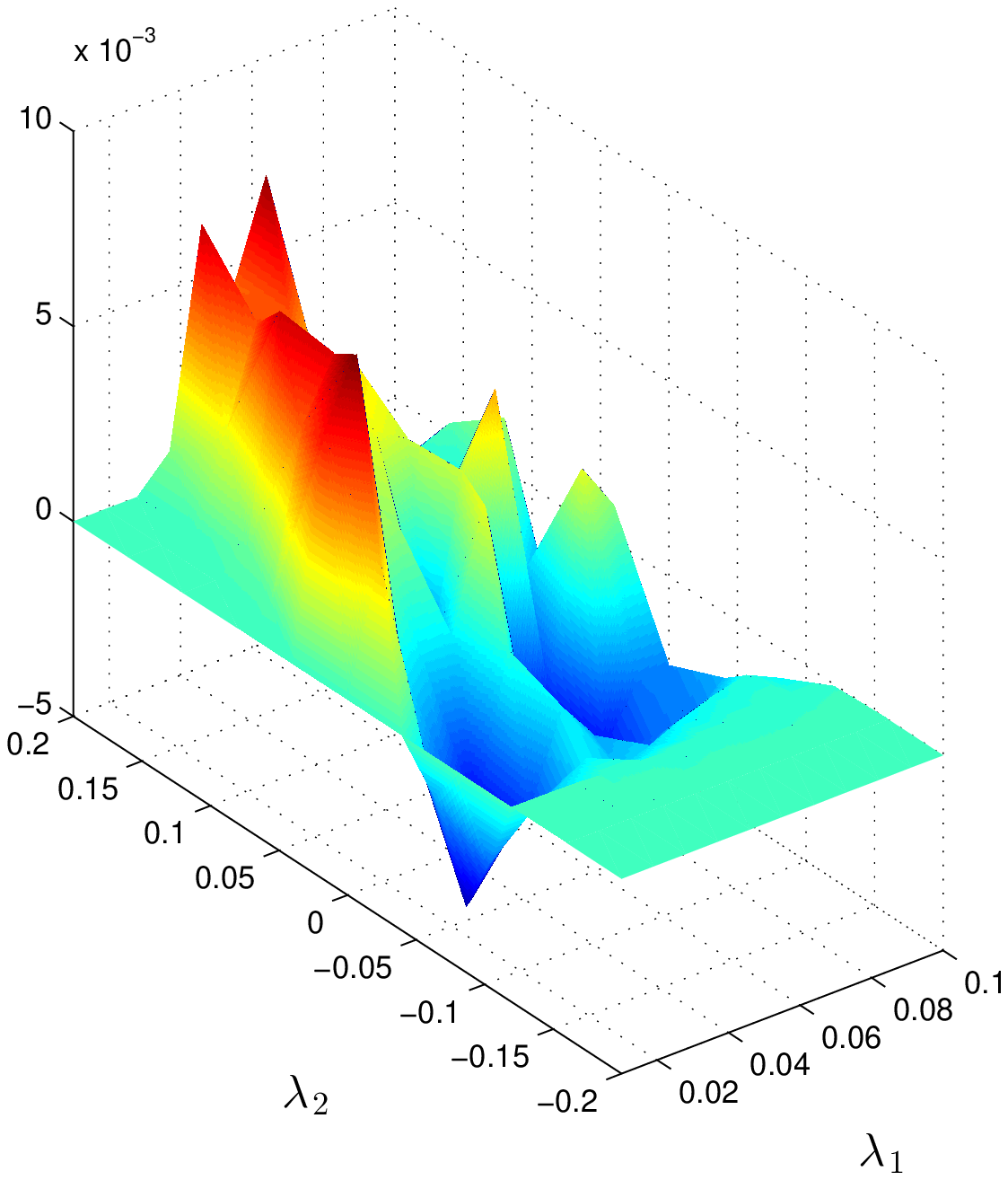}
		\caption{Plots of $\tilde{P}_{\Lambda,15k,h_1}-P_{\Lambda,62500}$ on the generating sets of  $\tilde{\mathcal{B}}_{\Lambda}$. }
	\label{fig:NS_LowRes_Errors}
\end{figure}

\begin{figure}[htb]
	\centering
\includegraphics[height=5cm,bb= 165 205 429 593,clip=true]{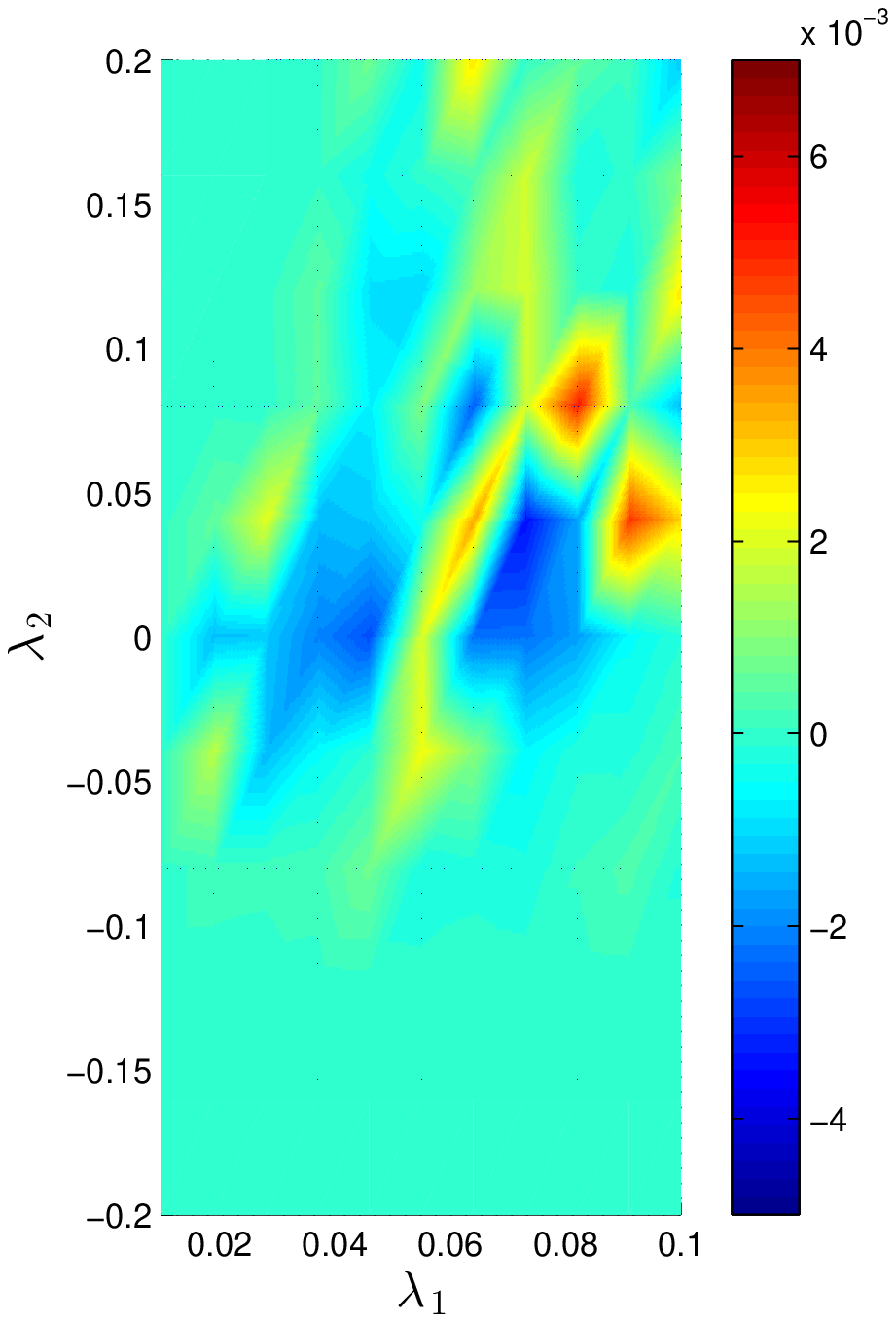}
\includegraphics[height=5cm,bb= 155 200 475 580,clip=true]{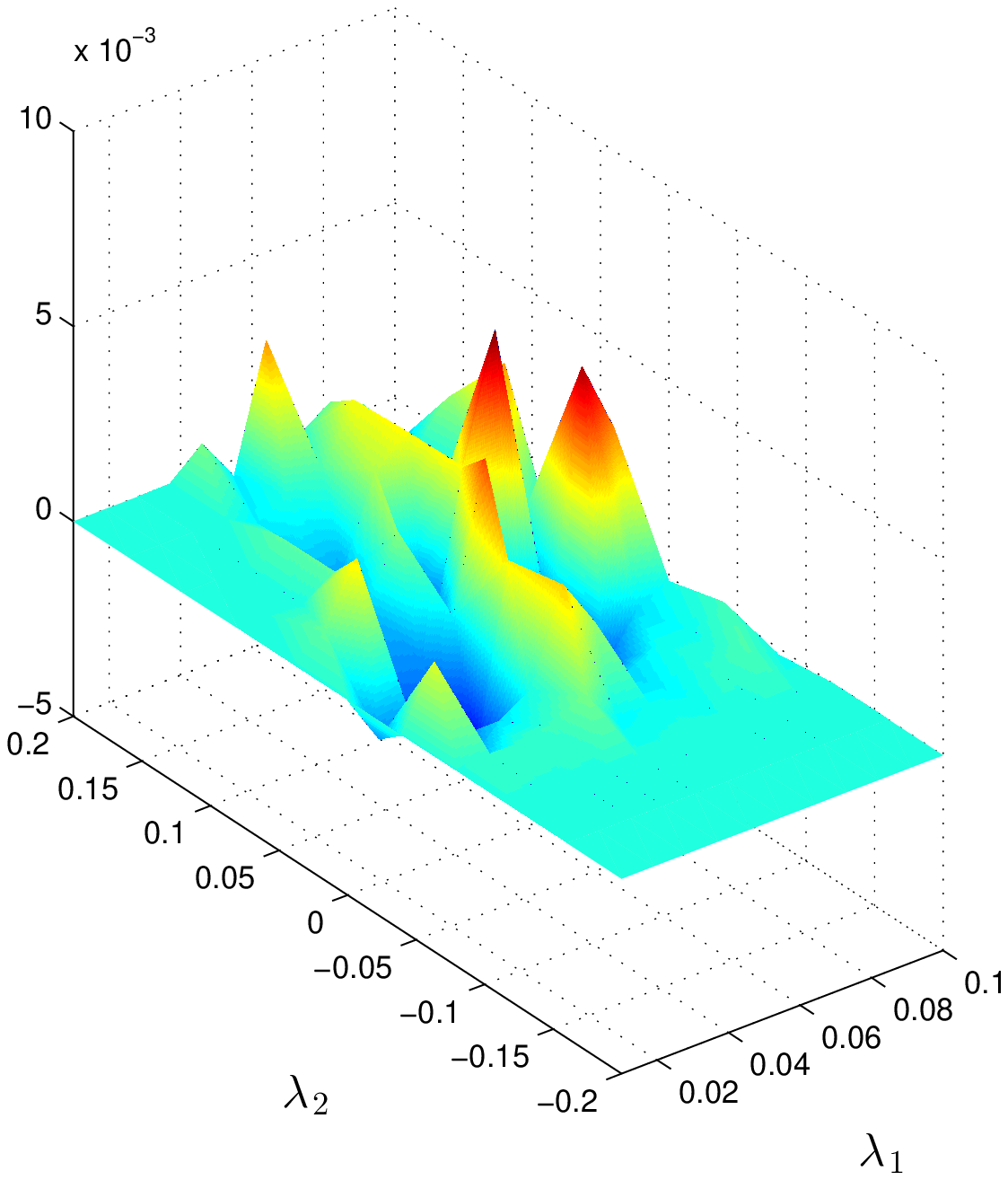}
		\caption{Plots of $\tilde{P}_{\Lambda,15k,h_2}-P_{\Lambda,62500}$ on the generating sets of  $\tilde{\mathcal{B}}_{\Lambda}$. }
	\label{fig:NS_LowRes_ErrCorrected_Errors}
\end{figure}

In Figures~\ref{fig:NS_LowRes_Errors} and \ref{fig:NS_LowRes_ErrCorrected_Errors}, we show the signed errors $\tilde{P}_{\Lambda,15k,h_1}-P_{\Lambda,62500}$ and $\tilde{P}_{\Lambda,15k,h_2}-P_{\Lambda,62500}$, respectively,  on the generating sets of $\tilde{\mathcal{B}}_{\Lambda}$. We observe a general reduction in the magnitude of the errors across large areas of the domain when we compute the counting measure using the a posteriori error estimates to correct for the numerical error.

\section{A higher dimensional example}\label{S:CM_Example}

We apply Algorithm~\ref{Alg_MC} to compute a counting measure for a 15-dimensional parameter space defined by uncertain coefficients and initial conditions of an MSEIRS model. The MSEIRS model is a generalization of the well studied SIR epidemic model that takes into account groups of infants protected by maternal antibodies (the ``M'' class) and groups of exposed and latent infected but not infectious (the ``E'' class) \cite{Hethcote_2000}. The coupled nonlinear system of differential equations defining this model is given by
\begin{equation}\label{eq:MSEIRS}
	\begin{cases}
		\frac{dM}{dt} &= B(S+E+I+R)-(\delta+\mu_M) M \\
		\frac{dS}{dt} &= \delta M - \beta SI - (\mu_G+\iota) S  + fR \\
		\frac{dE}{dt} &= \beta SI - (\epsilon + \mu_G)E \\
		\frac{dI}{dt} &= \epsilon E - (\gamma+\mu_I+\mu_G)I\\
		\frac{dR}{dt} &= \gamma I - (\mu_G+f)R+\iota S
	\end{cases}
\end{equation}
We define the uncertain parameters of this model in Table~\ref{tab:MSEIRS_params} including the bounds (shown without dimension and a single unit of time is 1 week, a unit of population is $1E6$, and birth/death rates are normalized to a population of size 300 million) that we use to define $\Lambda$.
\begin{table}[h]
\centering

 \begin{tabular}{l l}

 \begin{tabular}{| l | l | l | l |}
 \hline
 B & avg.~birth rate & $[ 2.72E-4,3.04E-4 ]$    &  $3.02E-4$  \\
 $\delta$ &  avg.~temporary immunity & $[1/12,1/4]$   & 0.16 \\
 $\mu_M$ & avg.~infant death rate & $[4E-3, 6E-3$  &  $4.5E-3$\\
 $\beta$ & contact/infectivity rate & $ [1.92E-3, 3.85E-3]$  & $3.4E-3$ \\
 $\mu_G$ & avg.~general death rate & $[2.4E-4,2.72E-4]$ &  $2.52E-4$\\
 $1/\epsilon$ & avg.~infection time& $[0.571, 1]$ & 0.7\\
 $\mu_I$ & avg.~infected death rate & $[4.81E-6 , 2.11E-5]$ & $1.75E-5$\\
 $1/\gamma$ & avg.~recovery time& $[0.7 , 2.33]$  & $0.8$\\
 $f$ & avg.~loss of immunity rate & $[0.125, 0.25]$ & 0.18\\
 $\iota$ & avg.~immunization rate & $[0.015,0.0375]$ & 0.026\\
 $M_0$ & initial infants & $[2.5,3.5]$ &  3.25 \\
 $S_0$ & initial susceptibles & $[260,275]$ & 270 \\
 $E_0$ & initial exposed & $[0.01,0.5]$& 0.425\\
 $I_0$ &  initial infected & $[0.1,4]$ & 3.8\\
 $R_0$ &  initial recovered/immunized & $[10,20]$ & 13\\
 \hline
 \end{tabular}

 \end{tabular}
 \caption{Uncertain parameters and initial conditions and their interval bounds in the MSEIRS model given by Eq.~\ref{eq:MSEIRS}. The last column of values are the reference parameter values used to define the mean of the output densities. }
 \label{tab:MSEIRS_params}
\end{table}

We sample the 15-dimensional $\Lambda$ uniformly $1E6$ times. For each parameter sample, we solve the model until a final time of six weeks numerically using the Dormand-Prince method coded within the ode45 function in Matlab, which implements a variable step Runge-Kutta method based on fourth-order error estimates \cite{Dormand_1980}. For the QoI, we take the number of immune infants and the number of infected individuals at the sixth week (denoted by $Q_1$ and $Q_2$, respectively). These quantities are uncertain (e.g., they may be estimated by surveys of hospital data during the disease outbreak), and it is generally impossible to know these values exactly at any particular time during the outbreak. Here, we consider the problem of parameter identification under uncertainty. We seek to determine the most probable configurations (defined by events) of parameters and initial conditions associated with the uncertain output data.

We use a tensor product of marginal probability densities defined by shifted and scaled (Beta) $B(4,5)$ distributions to model the uncertainty in the joint QoI (denoted by $Q=(Q_1,Q_2)$). These densities are defined by first solving the model with the reference parameter and initial conditions shown in Table~\ref{tab:MSEIRS_params}, which produced reference QoI values of approximately $Q_1=1.54$ and $Q_2=2.43$ (in millions). We determined the rescaling and shifting of the Beta distributions so that the range of any sampled QoI value stays within $\pm 0.15$ million people of the reference QoI value. We first invert only the marginal densities and then invert the joint density. To invert the marginal densities using Algorithm~\ref{Alg_MC}, we approximate each marginal using simple function approximations computed on a uniform grid of 200 subintervals on the computed ranges of $Q_1$ and $Q_2$. To invert the marginal density, we use $50\times 50$ uniformly sized rectangles to partition the Cartesian product of the ranges of $Q_1$ and $Q_2$ to compute the simple function approximation.

\begin{figure}[htb]
	\centering
\includegraphics[height=3.75cm,clip=true]{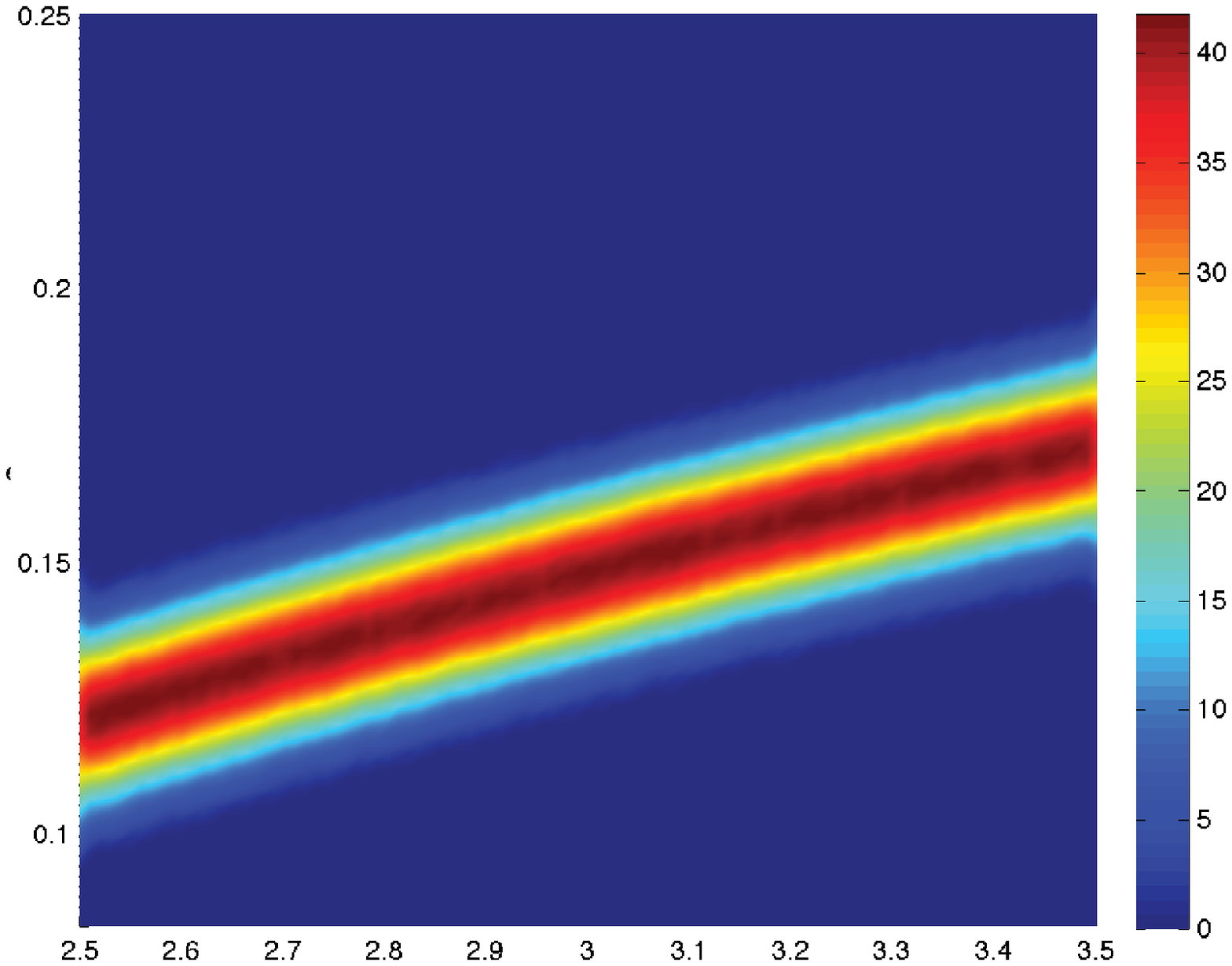}
\includegraphics[height=3.75cm,clip=true]{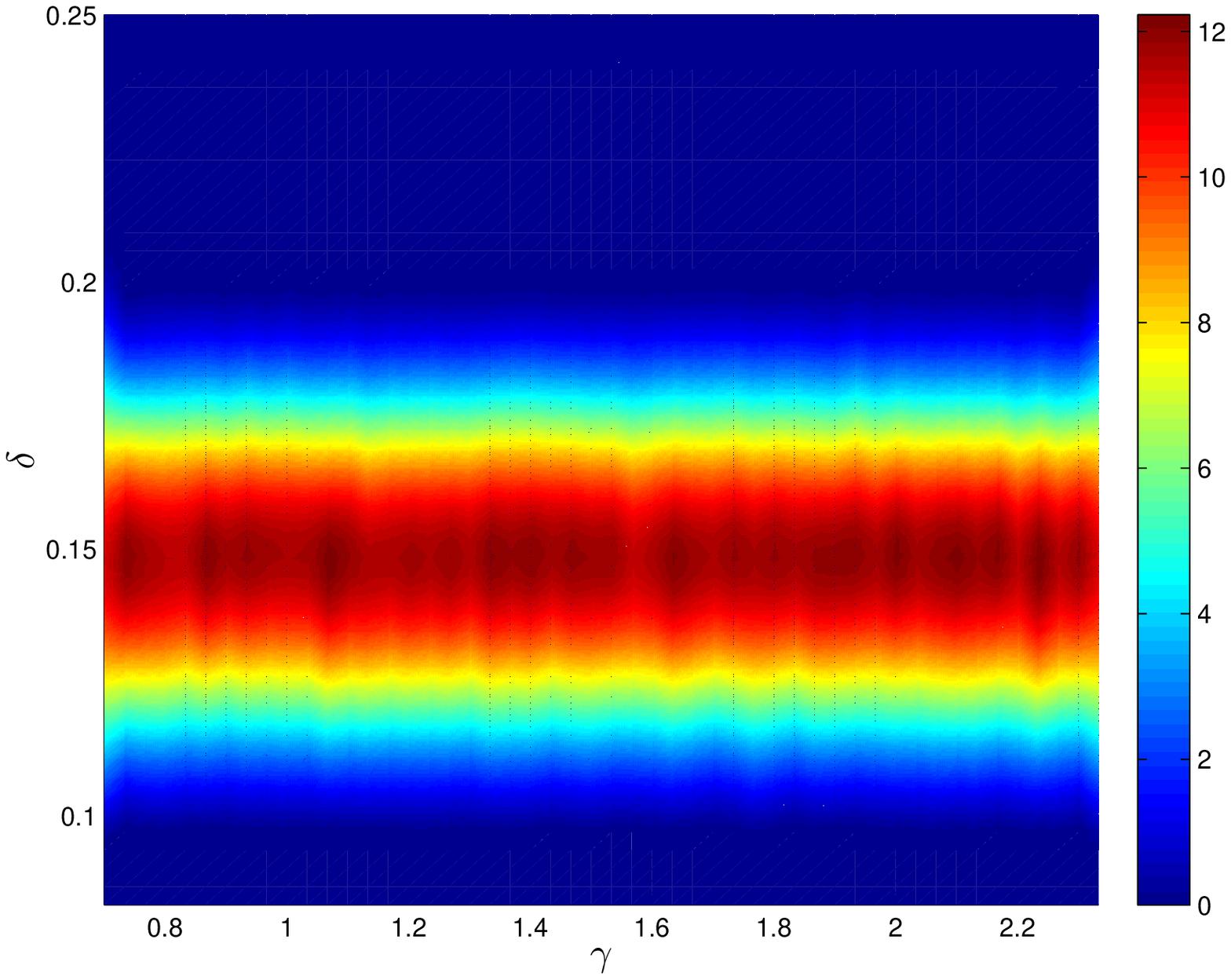}
\includegraphics[height=3.75cm,clip=true]{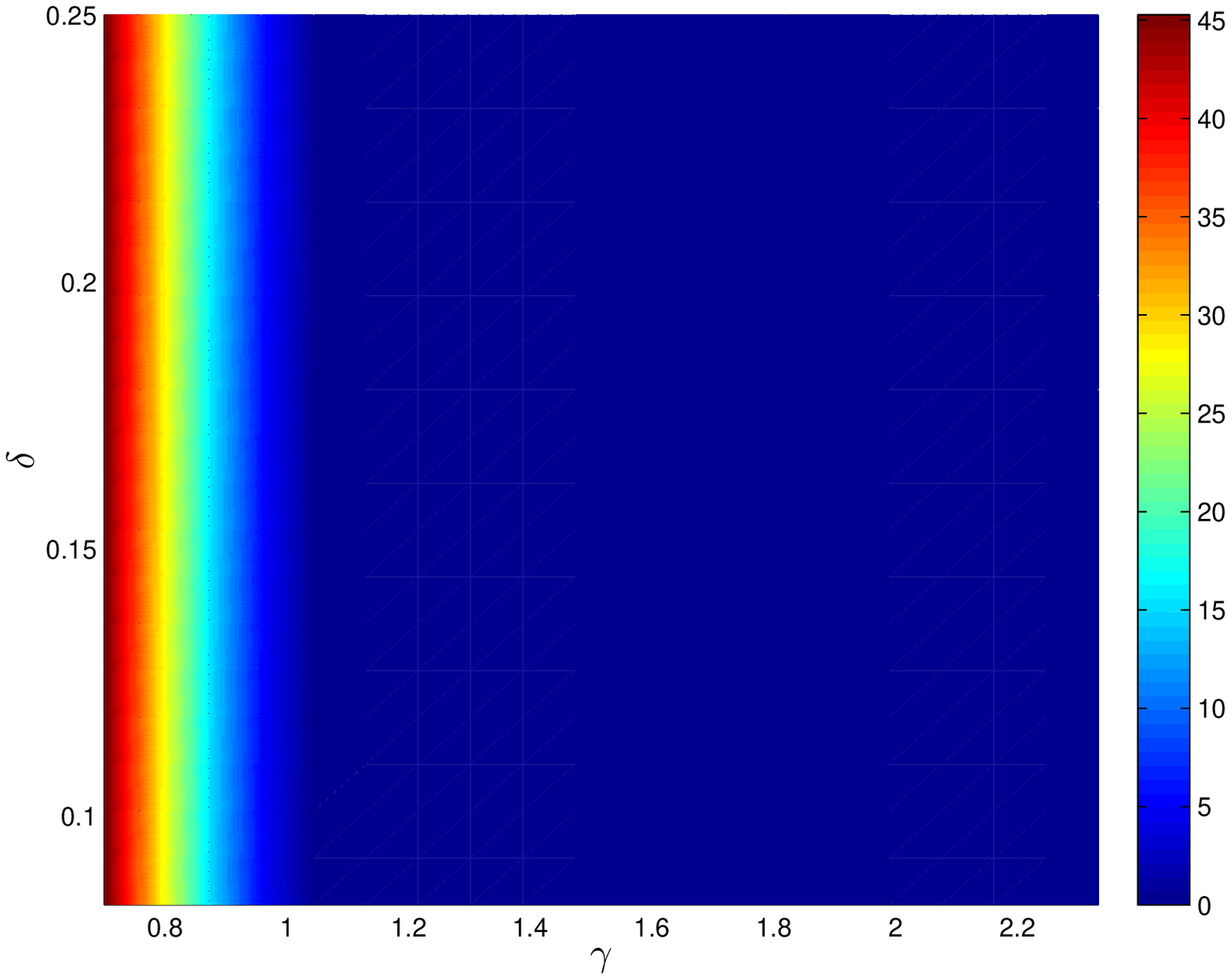}
\includegraphics[height=3.75cm,clip=true]{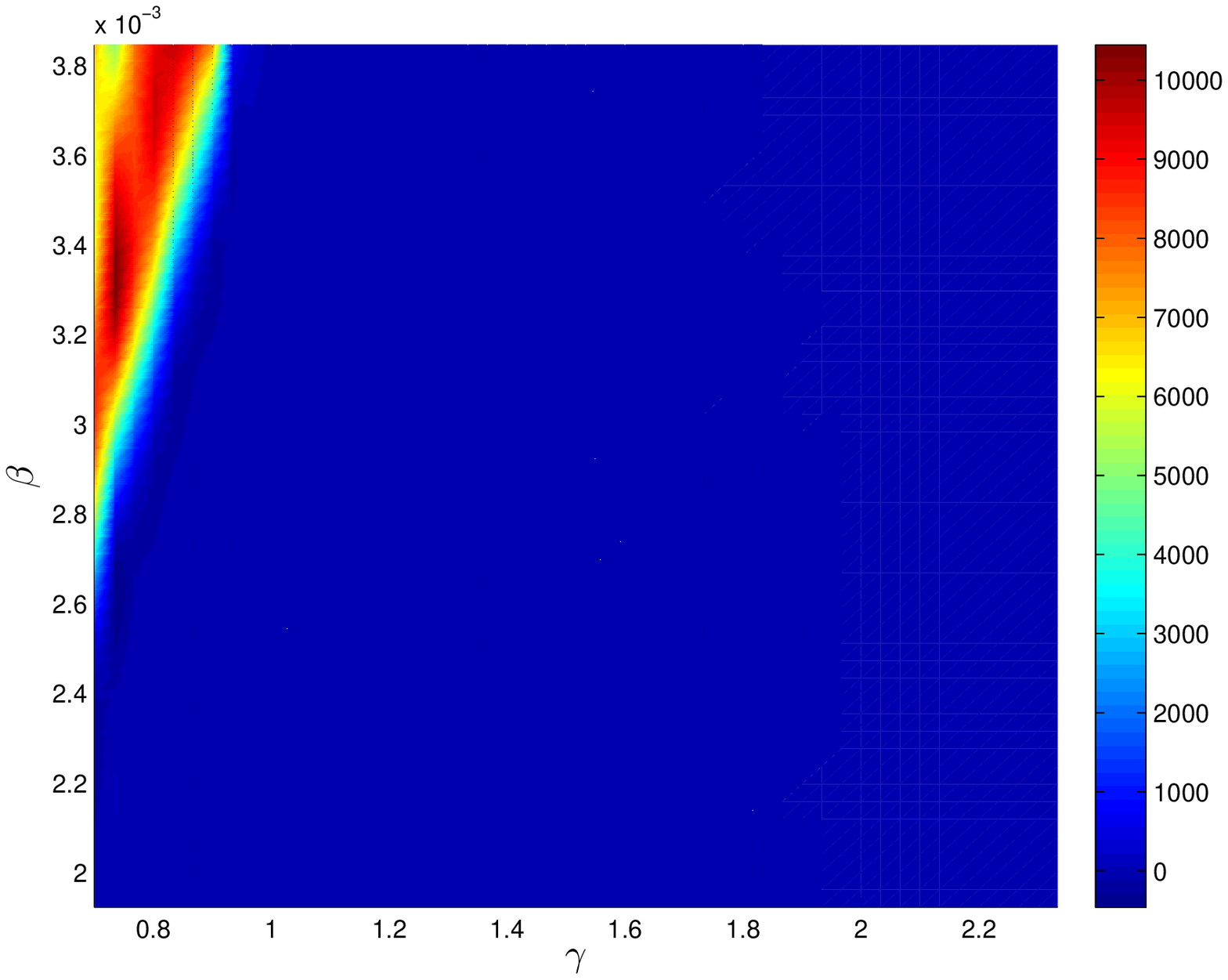}
\includegraphics[height=3.75cm,clip=true]{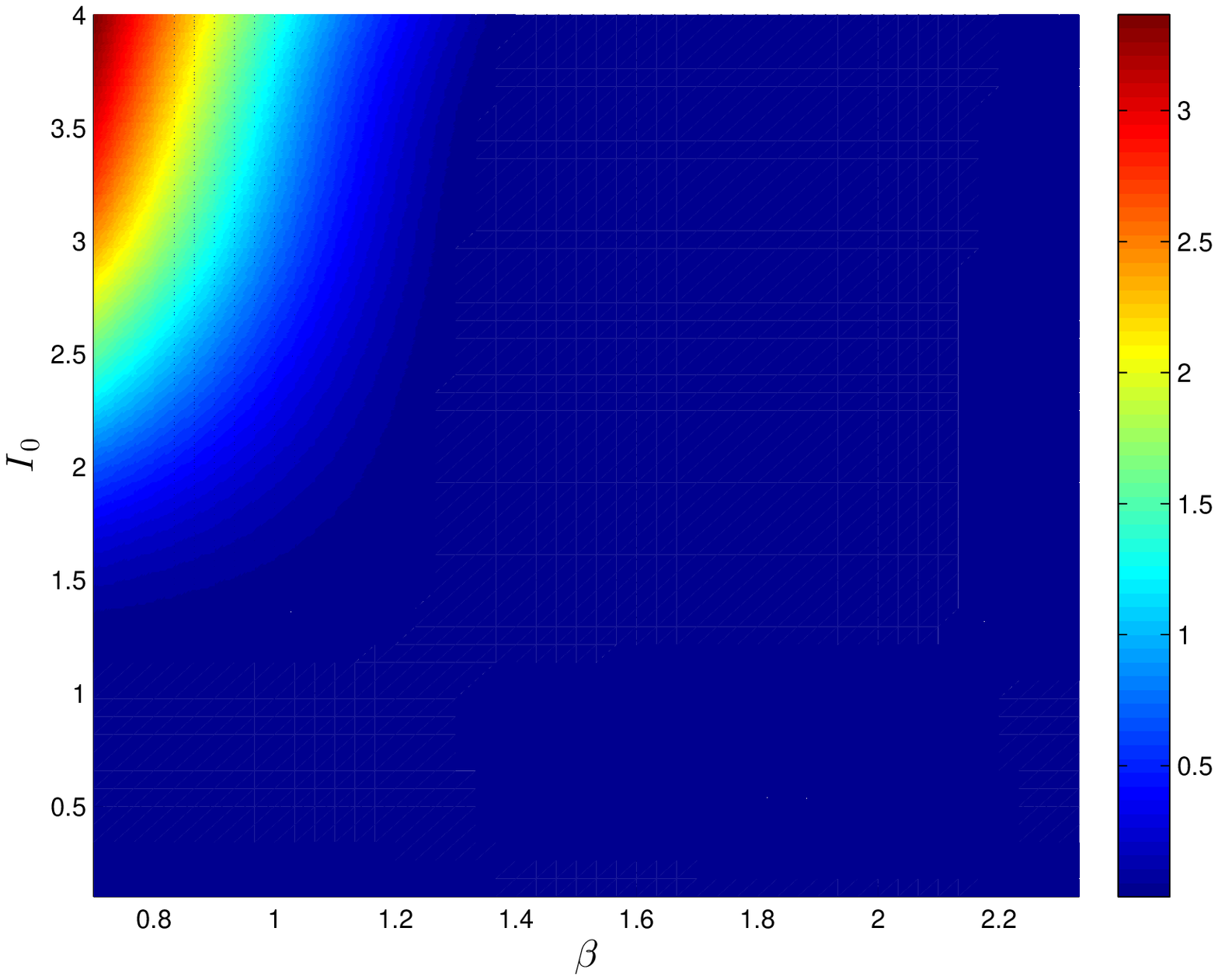}
\includegraphics[height=3.75cm,clip=true]{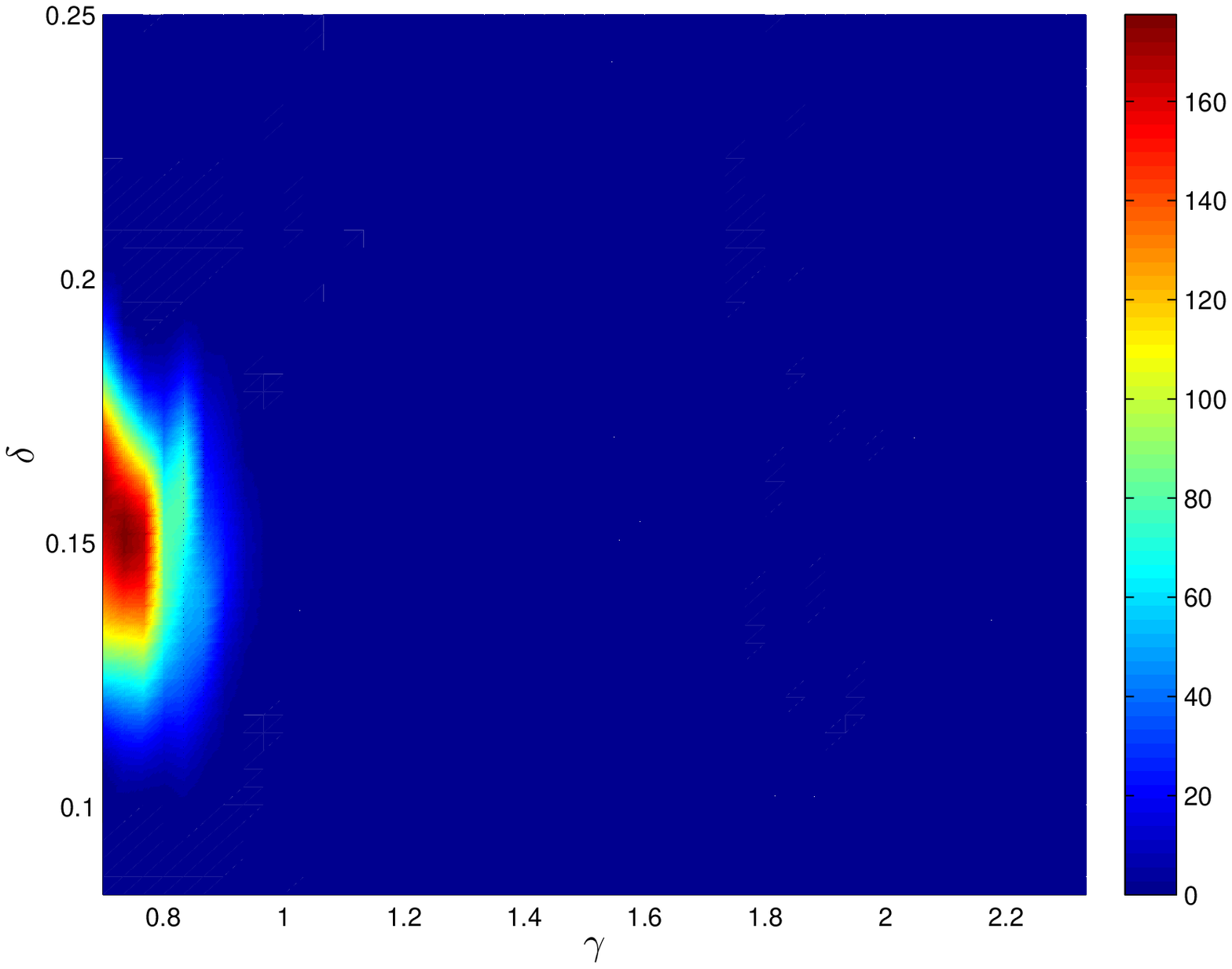}
		\caption{Thin plate spline approximations to some of the counting measure densities for various inverse problems. The top left and middle plots are marginals over $(M_0,\delta)$ and $(\gamma,\delta)$, respectively, computed from inverting the Beta density on $Q_1$. The top right plot is a marginal over $(\gamma,\delta)$ computed from inverting the density on $Q_2$. The bottom left and middle plots are marginals over $(\gamma,\beta)$ and $(\gamma,I_0)$, respectively, computed from inverting the Beta density on $Q_2$.  The bottom right plot is a marginal over $(\gamma,\delta)$ computed from inverting the joint output density on $Q$. }
	\label{fig:MSEIRS_marginals}
\end{figure}

In high dimensional parameter spaces, visualization and analysis of probability measures and/or densities can be difficult. Here, to highlight some of the underlying geometric structure, we show plots of some of the more interesting marginals for which $Q_1$ and $Q_2$ exhibit different sensitivities to the underlying parameters in Figure~\ref{fig:MSEIRS_marginals}. Most of the marginals over pairs of parameters not including those shown in Figure~\ref{fig:MSEIRS_marginals} were approximately uniform due to the lack of sensitivity of the QoI maps with respect to the parameters (i.e., the generalized contours were nearly flat in the directions of these parameters). We see that inverting the density for either $Q_1$, $Q_2$, or $Q$ all determine different marginal distributions on $(\gamma,\delta)$. In general, each QoI can only produce ``good estimates'' of two or three parameters in terms of a particular reference parameter or initial condition being in a relatively small event of high probability. Using both QoI provides good estimates of five of the 15 parameters and initial conditions. By utilizing additional geometrically distinct QoI, we can further improve results in terms of localizing the probability into regions of decreasing volume as described in \cite{Butler2014a}.

\appendix
\section{Proof of Uniqueness in Theorem~\ref{thm:stoch_error}}\label{App:A}

Suppose there is also a $C\in\mathcal{B}_{\Lambda,N}$ satisfying the equality for arbitrary output probability measures and $\mu_{\Lambda}(B\triangle C) >0$. By Lemma~\ref{lemma:nonempty}, there exists a $\lambda^{(k)}$ for some $k\in\set{1,\ldots,N}$ in either $B$ or $C$ but not both. Without loss of generality, let $\mathcal{K}_{C\backslash B}$ denote the set of indices such that $\lambda^{(k)}\in C\backslash B$ for all $k\in \mathcal{K}_{C\backslash B}$.  Let $I_{\mathcal{K}_{C\backslash B}}\subset\mathcal{D}$ denote the measurable event defined by $I_{\mathcal{K}_{C\backslash B}} = Q(\cup_{k\in\mathcal{K}_{C\backslash B}}\mathcal{V}(\lambda^{(k)}))\in\mathcal{B}_{\mathcal{D}}$ and $P_{\mathcal{D}}$ be defined by a uniform probability measure on $I_{\mathcal{K}_{C\backslash B}}$. This implies $\rho_{\mathcal{D}}$ is a simple function and we choose $\rho_{\mathcal{D},M}$ as this exact density. There are two cases to consider: (1) $\mu_{\Lambda}(B\backslash C) = 0$, or (2) $\mu_{\Lambda}(B\backslash C)\neq 0$. For case (1), we immediately have that $\tilde{P}_{\Lambda,N}(C)\neq \tilde{P}_{\Lambda,N}(B)$, which is a contradiction. For case (2), let $\mathcal{K}_{B\backslash C}$ denote the set of indices such that $\lambda^{(k)}\in B\backslash C$ for any $k\in \mathcal{K}_{B\backslash C}$. If $Q(\lambda^{(k)})\notin I_{\mathcal{K}_{C\backslash B}}$ for all $k\in \mathcal{K}_{B\backslash C}$, then we immediately arrive at the contradiction $\tilde{P}_{\Lambda,N}(C)\neq \tilde{P}_{\Lambda,N}(B)$ by construction of $\tilde{P}_{\Lambda,N}$. If there exists any $k\in\mathcal{K}_{B\backslash C}$ such that $Q(\lambda^{(k)})\in I_{\mathcal{K}_{C\backslash B}}$, then the corresponding Voronoi cells are identified as approximating induced regions of generalized contours with non-zero probabilities whose global approximation also includes Voronoi cells indexed by $\mathcal{K}_{C\backslash B}$. We are free to choose any Ansatz with the required modifications to Algorithm~\ref{Alg_MC} and choosing the Ansatz such that probabilities of parts of generalized contours through $B$ not included in $C$ are set to zero leads to a contradiction. $\Box$

\bibliographystyle{siam}
\bibliography{ReferencesBib}

\end{document}